\documentclass{amsart}
\usepackage{amsmath}
\usepackage{amsthm}
\usepackage{amsfonts}
\usepackage{amssymb}
\usepackage{dsfont}
\usepackage{xr}
\usepackage{graphicx}
\usepackage[all]{xy}

\newlength{\abstand}
\catcode`\ä=\active
\def ä{\"a}
\catcode`\ö=\active
\def ö{\"o}
\catcode`\ü=\active
\def ü{\"u}
\catcode`\Ä=\active
\def Ä{\"A}
\catcode`\Ö=\active
\def Ö{\"O}
\catcode`\Ü=\active
\def Ü{\"U}
\catcode`\ß=\active
\def ß{\ss{}}
\catcode`\é=\active
\def é{\'{e}}

\DeclareFontFamily{OT1}{pzc}{}
\DeclareFontShape{OT1}{pzc}{m}{it}{<-> s * [1.150] pzcmi7t}{}
\DeclareMathAlphabet{\mathcal}{OT1}{pzc}{m}{it}

\def\A{{\mathds A}}
\def\F{{\mathds F}}
\def\R{{\mathds R}}
\def\Rge{\R^{\geq 0}}
\def\Rp{\R^{>0}}
\def\C{{\mathds C}}
\def\N{{\mathds N}}
\def\Nn{\N_0}
\def\Np{\N_+}
\def\Q{{\mathds Q}}
\def\Z{{\mathds Z}}
\def\P{{\mathds P}}
\def\Ns{\str{\N}}
\def\Nns{\str{\Nn}}
\def\Nps{\str{\Np}}

\def\Qs{\str{\Q}}

\def\Zs{\str{\Z}}

\def\b{{\mathfrak b}}
\def\m{{\mathfrak m}}
\def\n{{\mathfrak n}}
\def\p{{\mathfrak p}}

\def\cB{\mathcal B}
\def\cC{\mathcal C}

\def\cE{\mathcal E}
\def\cF{\mathcal F}
\def\cG{\mathcal G}
\def\cH{\mathcal H}
\def\cI{\mathcal I}
\def\cJ{\mathcal J}

\def\cO{\mathcal O}

\def\cR{\mathcal R}
\def\cS{\mathcal S}

\def\cU{\mathcal U}

\def\bP{\mathbf{P}}
\def\scB{\str{\cB}}
\def\scR{\str{\cR}}
\def\scS{\str{\cS}}
\def\Rings{\mathcal{Rings}}
\def\Sets{\mathcal{Sets}}
\def\ALG{\mathcal{Alg}}

\def\SCH{\mathcal{Sch}}
\def\SCHFP{\SCH^{\mathrm{fp}}}
\def\MOD{\mathcal{Mod}}
\def\MODFP{\MOD^{\mathrm{fp}}}
\def\COH{\mathcal{Coh}}
\def\QCOH{\mathcal{QCoh}}
\def\RF{\mathrm{R}}
\def\Ks{\str{K}}
\def\Ksf{\Ks^{\mathrm{fin}}}
\def\Ksi{\Ks^{\mathrm{inf}}}
\def\Kh{\hat{K}}
\def\Ksh{\str{\Kh}}
\def\Mi{\MM{1}}
\def\Mii{\MM{2}}
\def\Miii{\MM{3}}
\def\sh{\mathrm{sh}}
\def\eps{\varepsilon}
\def\Guniv{\cG^{\mathrm{univ}}}

\newcommand\spec[1]{{\mathrm{Spec}\,(#1)}}

\newcommand\str[1]{{\mbox{}^*#1}}
\newcommand\Mor[3]{\mathrm{Mor}_{#1}(#2,#3)}
\newcommand\id[1]{\mathds{1}_{#1}}
\newcommand\Alg[1]{\ALG_{#1}}
\newcommand\Sch[1]{\SCH_{#1}}
\newcommand\Algfp[1]{\ALG^{\mathrm{fp}}_{#1}}
\newcommand\Schfp[1]{\SCH^{\mathrm{fp}}_{#1}}
\newcommand\sAlg[1]{\str{\Alg{#1}}}
\newcommand\sAlgfp[1]{\str{\Algfp{#1}}}
\newcommand\sSch[1]{\str{\Sch{#1}}}
\newcommand\sSchfp[1]{\str{\SCH^{\mathrm{fp}}_{#1}}}
\newcommand\op[1]{{#1}^{\mathrm{op}}}
\newcommand\Mod[1]{\MOD_{#1}}
\newcommand\sMod[1]{\str{\Mod{#1}}}
\newcommand\Modfp[1]{\MOD_{#1}^{\mathrm{fp}}}
\newcommand\sModfp[1]{\str{\Modfp{#1}}}
\newcommand\Coh[1]{\COH_{#1}}
\newcommand\sCoh[1]{\str{\COH_{#1}}}
\newcommand\QCoh[1]{\QCOH_{#1}}
\newcommand\sQCoh[1]{\str{\QCOH_{#1}}}

\newcommand\s[1]{N\,#1}
\newcommand\sls[1]{N_*#1}
\newcommand\sm[1]{\sigma_{#1}}
\newcommand\B[1]{S\,#1}
\newcommand\T[1]{T\,#1}
\newcommand\tm[1]{\rho_{#1}}
\newcommand\Ob[1]{\mathrm{Ob}(#1)}
\newcommand\bs[3]{\tau_{#1,#2,#3}}
\newcommand\bsi[3]{\tau^{-1}_{#1,#2,#3}}
\newcommand\HOM[3]{{\mathcal{Hom}}_{\cO_{#1}}(#2,#3)}
\newcommand\Hom[3]{\mathrm{Hom}_{#1}(#2,#3)}
\newcommand\HZar[3]{\mathrm{H}^{#1}(#2,#3)}
\newcommand\Rf[3]{\mathrm{R}^{#1}{#2}_*#3}

\newcommand\Der[1]{\mathrm{D}^+(#1)}
\newcommand\Derb[1]{\mathrm{D}^\mathrm{b}(#1)}
\newcommand\MMod[1]{#1\hbox{-}\mathbf{Mod}}
\newcommand\MM[1]{\textbf{(M#1)}}
\newcommand\ball[3]{U_{#1}(#2,#3)}
\newcommand\Pic[1]{\mathrm{Pic}(#1)}
\newcommand\sPic[1]{\str{\Pic{#1}}}
\newcommand\image[2]{\mathrm{Im}\,\left(#1\rightarrow #2\right)}
\newcommand\FQuot[4]{\mathrm{Quot}^{#1}(#2/#3/#4)}
\newcommand\sFQuot[4]{\str{\FQuot{#1}{#2}{#3}{#4}}}
\newcommand\Quot[4]{\mathrm{Quot}^{#1}_{#2/#3/#4}}
\newcommand\sQuot[4]{\str{\Quot{#1}{#2}{#3}{#4}}}
\newcommand\Hilb[3]{\mathrm{Hilb}^{#1}_{#2/#3}}
\newcommand\sHilb[3]{\str{\Hilb{#1}{#2}{#3}}}
\newcommand\eq[2]{\mathrm{eq}(#1,#2)}
\newcommand\Eq[2]{\mathrm{Eq}(#1,#2)}
\newcommand\Homf[4]{\mathrm{Hom}^{#1}_{#2}(#3,#4)}
\newcommand\sHomf[4]{\str{\Homf{#1}{#2}{#3}{#4}}}
\newcommand\QK[1]{\mathrm{Quot}(#1)}
\newcommand\D[1]{\mathrm{D}(#1)}

\swapnumbers
\theoremstyle{definition}
\newtheorem{defi}{Definition}[section]
\newtheorem{bsp}[defi]{Example}

\newtheorem{satzdefi}[defi]{Proposition/ Definition}
\newtheorem{lemma}[defi]{Lemma}
\newtheorem{lemmadefi}[defi]{Lemma/ Definition}
\newtheorem{bem}[defi]{Remark}
\newtheorem{satz}[defi]{Proposition}
\newtheorem{thm}[defi]{Theorem}
\newtheorem{cor}[defi]{Corollary}

\title{Enlargements of Schemes}
\author{Lars Brünjes, Christian Serpé}
\date{\today}
\thanks{Parts of this work have been written under hospitality of the
  Mathematical Research Institute Oberwolfach}
\address{Universität Regensburg \\ NWF I - Mathematik \\ D-93040 Regensburg \\ Germany}
\email{lars.bruenjes@mathematik.uni-regensburg.de}
\address{
  Christian Serpé \\
  Westfälische Wilhelms-Universität Münster \\
  Mathematisches Institut \\
  Sonderforschungsbereich 478 ``Geometrische Strukturen in der Mathematik'' \\
  Hittorfstr. 27 \\
  D-48149 Münster \\
  Germany}
\email{serpe@math.uni-muenster.de}


\date{\today}

\begin{document}

\begin{abstract}
  In this article we use our constructions from \cite{enlcat}
  to lay down some foundations for the application of A.~Robinson's nonstandard methods
  to modern Algebraic Geometry.
  The main motivation is the search for another tool to transfer results from
  characteristic zero to positive characteristic and vice versa. We give
  applications to the resolution of singularities and weak factorization.
\end{abstract}

\maketitle


\section{Introduction}

\vspace{\abstand}

The difficulty of
many problems about algebraic varieties
depends on the characteristic of the base field.
Resolution of singularities (proved in characteristic zero, open in characteristic $p$)
and Grothendieck's standard conjecture on the rationality of Künneth components
(proved over finite fields, open in characteristic zero) are prominent examples.
This is mostly due to the fact that some tools --- like transcendental methods --- are only available in
characteristic zero while others --- like Frobenius morphisms --- only exist in characteristic $p$.

A link between the apparently so different worlds of characteristic zero and characteristic $p$
is provided by internal fields of \emph{infinite} characteristic,
for example the *finite field $\Zs/P$, where $\Zs$ is an enlargement of $\Z$
and $P\in\Zs$ is an infinite prime:

Let $\Phi$ be a \emph{first order statement} in the \emph{language of fields}. If $\Phi$ is true for all fields of
characteristic zero, it is in particular true for $\Zs/P$ (which externally has characteristic zero),
so by the permanence principle it is true for $\F_p$ for infinitely many \emph{finite} primes $p\in\Z$.
If, on the other hand, $\Phi$ is true for almost all $\F_p$, it is also true for $\Zs/P$, a field of characteristic zero.

Unfortunately, being first order is a strong condition in whose absence the above reasoning fails,
and the language of fields is ill adapted to dealing with schemes, sheaves and cohomology
in Grothendieck's modern language of Algebraic Geometry.

Building on our paper \cite{enlcat},
we therefore use the notion of \emph{enlargement of categories}
to establish a more flexible method of transferring properties from characteristic zero to characteristic $p$
and vice versa in the framework of schemes:

Starting from a category $\cB$ of rings,
we consider the fibred category $\Schfp{\cB}/\cB$ of finitely presented schemes over objects of $\cB$
and enlarge it to get the category of *schemes $\sSchfp{\cB}$, fibred over $\scB$.
Here the main point is the following:

An object $A$ of $\scB$
is also an ordinary ring,
and we can consider the category $\Schfp{A}$ of finitely presented schemes over $A$.
The notion of scheme is \emph{not} first order, so an object $X$ of $\Schfp{A}$ is \emph{not}
an $A$-*scheme. Nevertheless, $X$ is given by finitely many equations in finitely many unknowns,
and these define a *scheme $\s{X}$ over $A$
(in fact, we construct a canonical fibred functor from $\Schfp{\scB}$ to $\sSchfp{\cB}$,
which turns out to be a fibred Kan extension and is therefore unique up to unique isomorphism).
Similarly, any finitely presented $\cO_X$-module $\cF$ defines a *finitely-presented $\cO_{\s{X}}$-*module
given by "the same" presentation.
For modules, there is even a canonical functor $\B{}$ in the opposite direction, sending $\cO_{\s{X}}$-*modules to
$\cO_X$-modules, and the functors $\s{}$ and $\B{}$ turn out to have many nice properties.

The main part of our paper is devoted to proving that many properties of $X$ (like for example being smooth or proper)
translate into corresponding properties of $\s{X}$.
--- Let us stress the fact that this is \emph{not} simply an application of the transfer principle,
because the standard universe does not contain $A$ and $X$ and is thus not applicable.

Especially in the case where $A$ is a \emph{field}, properties of $\s{X}$ often also
imply corresponding properties of $X$ ---
for example, $\s{X}$ is *irreducible respectively *integral if
\emph{and only if} $X$ is irreducible respectively integral.

Furthermore, we can give criteria (mostly of cohomological nature) for whether
a given *scheme or *module lies in the essential image of $\s{}$,
thus enabling us to deduce the existence of schemes and modules with certain properties from the existence
of *schemes and *modules with the corresponding properties
(note that there are many *scheme which do no lie in the essential image of $\s{}$, for example *schemes
of *finite but infinite *dimension and *schemes given by equations of *finite but infinite *degree).

At this point, let us mention Angus Macintyre's ``many sorted'' approach to the application of Model Theory to
Algebraic Geometry in \cite{macintyre},
where he considers ultraproducts of varieties (and algebraic cycles) of fixed complexity.
Though a direct comparison between Macintyre's approach and ours is difficult due to the different languages used,
*schemes respectively *schemes in the essential image of $\s{}$
correspond to ultraproducts of varieties of arbitrary respectively bounded complexity.

The announced method of transfer between characteristic zero and characteristic $p$ now roughly works as follows:
Let $\Phi$ be a statement of schemes.
Assume first that $\Phi$ holds in characteristic zero, and consider a class $C$
of *schemes over *fields which lie in the essential image of $\s{}$
(i.e. a class of ``bounded complexity", for example the class of *projective *schemes
whose *dimension and *degree is bounded by a finite number).
If $k$ is a *field in $\scB$ of infinite *characteristic, $\Phi$ holds for schemes over $k$ (which has characteristic zero as a field),
and using properties of $\s{}$, it will often be possible to show that $\str{\Phi}$ then holds for *schemes in $C$, hence
$\Phi$ holds for (certain) schemes over fields of \emph{finite}
characteristic (by the permanence principle). --- We will give two applications of this method, namely to the problems of
resolution of singularities and of weak factorization in characteristic $p$.

If, on the other hand, $\Phi$ holds for schemes in characteristic $p$,
by transfer $\str{\Phi}$ holds for *schemes over *fields $k$ in $\scB$ of infinite *characteristic,
so if $X$ is a scheme over $k$, $\str{\Phi}$ holds for $\s{X}$. Again, using properties of $\s{}$,
it will often be possible to use this fact to prove that $\Phi$ holds for $X$, a scheme in characteristic zero.
For example, if the (modified) Jacobian conjecture was proven is characteristic $p$,
this method, combined with an easy application of the Lefschetz principle, would imply the Jacobian conjecture over $\Q$.

In subsequent papers, we plan to define similar functors $\s{}$ for $K$-theory, cycles
and étale cohomology,
and even though we demonstrate the usefulness of our method as it stands in the present paper
(and it will not be hard to find other applications along similar lines),
our main motivation for this paper is to lay the ground for that future work,
from which we hope to gain new insights into the theory of algebraic cycles over varieties in characteristic zero
and characteristic $p$.

\vspace{\abstand}

The paper is organized as follows: In the second section we give
basic definitions; in particular we define the fibration
$\Schfp{\cB}/\cB$ of finitely presented schemes over a small category
of rings $\cB$ and consider the enlargement $\sSchfp{\cB}/\scB$.

In the third section we relate schemes and *schemes. For that, we
define a functor $\s{}:\Schfp{\scB}/\scB\rightarrow\sSchfp{\cB}/\scB$ which
extends the canonical functor $\Schfp{\cB}/\cB\rightarrow\sSchfp{\cB}/\scB$.
In particular, for an internal ring $A$, we get a
functor $\s{}:\Schfp{}/A\rightarrow \sSchfp{}/A$.

Section 4 discusses more
properties of the functor $\s{}$ and shows that it respects
many properties of morphism between schemes.

In section 5 we define and investigate an analogous
functor $\s{}$ for coherent modules. That is, for a scheme $X$ of finite
presentation over an internal ring, we define a functor from
coherent modules on $X$ to *coherent modules on $\s{X}$.

Section 6 specializes to the case where the internal ring $A$ is actually an internal
field. Mainly, we apply a theorem of van den Dries and Schmidt to
show --- among other things --- that the functor $\s{}$ on modules is exact and
that the functor $\s{}$ on schemes is compatible with Quot- and Hilbert
schemes.

In section 7 we show that $\s{}$ is compatible with higher direct images of coherent sheaves for proper morphisms
(the proof of this is similar to the proof of the theorem on formal functions in Algebraic Geometry).
One main application of this theorem is that $\s{}$
is fully faithful on coherent modules and induces an injection on
Picard groups.

Section 8 shows that it is possible to define a kind of shadow
map for varieties over an internal valued field with locally compact
completion.

In section 9 finally we give two standard applications of the
theory: First we reprove a result on resolution of
singularities in characteristic $p$ by Eklof,
and second we show a similar result for the factorization of birational morphisms.

\vspace{\abstand}


\section{Basic definitions}

\vspace{\abstand}

Let $\Rings$ be the category of rings,
let $\cB\subset\Rings$ be a \emph{small} (not necessarily full) subcategory,
let $\cR$ be the small full subcategory of $\Rings$ containing every object of $\cB$ and (an isomorphic image of)
every ring finitely presented over $\Z$ or over an object of $\cB$,
and let $\cS$ be the small full subcategory of the category of schemes containing (an isomorphic image of)
every scheme which is finitely presented over an object of $\cR$.

Choose a universe $\cU$ such that $\cS$ is $\cU$-small,
and choose a superstructure $\hat{M}$ containing $\cU$ (such that any $\cU$-small category is $\hat{M}$-small ---
compare \cite[A.3]{enlcat}).

Let $*:\hat{M}\rightarrow\widehat{\str{M}}$ be an enlargement.
Since $\cS$ is $\hat{M}$-small, so are $\cB$ and $\cR$,
and we can consider the enlargements $\scB\subseteq\scR$ and $\scS$, all $\widehat{\str{M}}$-small categories,
where $\scB$ and $\scR$ can be thought of as categories of (internal) rings with
(internal) ring homomorphisms as morphisms (compare \cite[4.7]{enlcat}).

We call objects of $\scR$ \emph{*rings} and objects of $\scS$ \emph{*schemes}.

\vspace{\abstand}

Define $\SCH$ to be the category whose objects
are morphisms $X\rightarrow\spec{S}$ (with $X$ an arbitrary scheme and $S$ an arbitrary ring)
and whose morphisms $[X'\xrightarrow{\pi_{X'}}\spec{S'}]\rightarrow[X\xrightarrow{\pi_X}\spec{S}]$
are pairs $\langle X'\xrightarrow{f}X,\ S\xrightarrow{\varphi}S'\rangle$
such that the following square commutes:
\[
  \xymatrix@C=20mm{
    {X'} \ar[r]^f \ar[d]_{\pi_{X'}} &
    {X} \ar[d]^{\pi_X} \\
    {\spec{S'}} \ar[r]^{\spec{\varphi}} &
    {\spec{S}.} \\
  }
\]
(If the morphism $X\rightarrow\spec{S}$ is understood, we often denote the object $X\rightarrow\spec{S}$
by $X/S$ or --- if $S$ is understood as well - simply by $X$.)

Projection onto the second component defines a functor $\SCH\rightarrow\op{\Rings}$ which is obviously a bifibration:
For a ring homomorphism $\varphi:S\rightarrow S'$, inverse and direct image are given by
\[
  \begin{array}{lclc}
    \varphi^*[X\rightarrow\spec{S}] & = & [X\times_SS'\rightarrow\spec{S'}] & \text{and} \\[2mm]
    \varphi_*[X'\rightarrow\spec{S'}] & = & [X'\rightarrow\spec{S'}\xrightarrow{\spec{\varphi}}\spec{S}].
  \end{array}
\]
The fibre over a ring $S$ is obviously the category $\Sch{S}$ of $S$-schemes.

Let $\SCHFP$ be the full subcategory of $\SCH$ consisting of morphisms $X\rightarrow\spec{S}$
with $X$ a \emph{finitely presented} $S$-scheme.
Then $\SCHFP$ is a subfibration of $\SCH$ over $\Rings$ (but no longer a bifibration,
because for a ring homomorphism $S\rightarrow S'$, not every finitely presented $S'$-scheme will in general
be finitely presented as an $S$-scheme).
Of course, the fibre over a ring $S$ is the category $\Schfp{S}$ of finitely presented $S$-schemes.

\vspace{\abstand}

For an arbitrary subcategory $\cC$ of $\Rings$,
we can form the pullbacks of $\SCH\rightarrow\op{\Rings}$ and $\SCHFP\rightarrow\op{\Rings}$
along $\op{\cC}\rightarrow\op{\Rings}$,
and we denote the resulting bifibration respectively fibration over $\op{\cC}$
by $\Sch{\cC}$ respectively $\Schfp{\cC}$.
\footnote[2]{When we view $\SCH$ and $\SCHFP$ as pseudo-functors from $\op{\Rings}$ to the category of categories,
then $\Sch{\cS}$ and $\Schfp{\cS}$ are just the restrictions of these functors to $\op{\cC}$.}

Applying this to $\cC:=\cB$ and $\cC:=\scB$, we get bifibrations
$\Sch{\cB}\rightarrow\op{\cB}$ and $\Sch{\scB}\rightarrow\op{\scB}$
and fibrations $\Schfp{\cB}\rightarrow\op{\cB}$ and $\Schfp{\scB}\rightarrow\op{\scB}$.

\vspace{\abstand}

Since the fibrations $\Schfp{\cR}\rightarrow\op{\cR}$ and $\Schfp{\cB}\rightarrow\op{\cB}$ are obviously
$\hat{M}$-small,
we can consider their enlargements
\[
  \xymatrix{
    {\sSchfp{\scB}} \ar@{^{(}->}[r] \ar[d] & {\sSchfp{\scR}} \ar[d] \\
    {\op{\scB}} \ar@{^{(}->}[r] & {\op{\scR}} \\
  }
\]
which are again fibrations (compare \cite[7.3]{enlcat}),
whose fibres we denote by $\sSchfp{S}$ for objects $S$ of $\scR$.

\vspace{\abstand}

\begin{defi}\label{defstarscheme}\mbox{}\\[-2mm]
  \begin{enumerate}
    \item
      For a *ring $S$, we call the category $\sAlg{S}:=\scR\backslash S$ of objects under $S$
      the category of \emph{$S$-*algebras}.
    \item
      By transfer we have a functor $\str{\mathrm{Spec}}:\op{\scR}\rightarrow\scS$ from *rings
      to *schemes, and we call *schemes in the essential image of this functor \emph{*affine}.
    \item
      For a *scheme $X$, we call the category $\sSch{X}:=\scS/X$ of objects over $X$
      the category of \emph{$X$-*schemes}
      or --- if $X=\str{\spec{A}}$ is *affine ---
      the category $\sSch{A}$ of \emph{$A$-*schemes}.
    \item
      Let $\bP$ be a property of rings (schemes, morphisms of rings, morphisms of schemes).
      When considering $\bP$ as a predicate on the set of objects of $\cR$ (of objects of $\cS$,\ldots),
      we get a predicate $\str{\bP}$ on the set of objects of $\scR$ (of objects of $\scS$,\ldots),
      i.e. a property of *rings (*schemes, morphisms of *rings, morphisms of *schemes).
  \end{enumerate}
\end{defi}

\vspace{\abstand}

\begin{bem}\label{bemsternringe}
  It follows immediately from transfer that objects of $\sSchfp{\scS}$
  are morphisms of *schemes $X\rightarrow\str{\spec{S}}$,
  where $S$ is a *ring and $X$ is a *scheme.
  Morphisms $[X'\rightarrow\str{\spec{S'}}]\rightarrow[X\rightarrow\str{\spec{S}}]$
  are pairs $\langle f,\varphi\rangle$
  with $f\in\Mor{\scS}{X'}{X}$ and $\varphi\in\Mor{\scR}{S}{S'}$ such that the following square
  commutes in $\scS$:
  \[
    \xymatrix@C=20mm{
      {X'} \ar[r]^f \ar[d] &
      {X} \ar[d] \\
      {\str{\spec{S'}}} \ar[r]^{\str{\spec{\varphi}}} &
      {\str{\spec{S}}.} \\
    }
  \]
  In particular, for a *ring $S$, the fibre $\sSchfp{S}$ is the full subcategory of the category of $S$-*schemes
  defined in \ref{defstarscheme} consisting only of *finitely presented $S$-*schemes.
\end{bem}

\vspace{\abstand}

\begin{defi}\label{defstarpol}
  Consider the functor $\text{Pol}:\N_0\times\cR\rightarrow\cR$
  (where $\N_0$ is the category associated to the partially ordered set $(\N_0,\leq)$),
  sending a pair $(n,S)$ to the polynomial ring $S[X_1,\ldots,X_n]$.
  This is a functor between $\hat{M}$-small categories,
  so we can enlarge it to a functor $\str{\text{Pol}}:\str{\N_0}\times\scR\rightarrow\scR$.
  For a (not necessarily finite) natural number $n\in\str{\N_0}$ and a *ring $S$,
  we denote $\str{\text{Pol}}(n,S)$ by $S\str{[X_1,\ldots,X_n]}$
  and call it the \emph{*polynomial ring over $S$ in $n$ unknowns}.
\end{defi}

\vspace{\abstand}

\begin{bem}\label{bemstarpol}
  Let $(n,S)$ be an object of $\str{\N_0}\times\scR$ as above.
  \begin{enumerate}
    \item
      The morphism $\str{\text{Pol}}(0\leq n,\id{S}):S=S\str{[]}\rightarrow S\str{[X_1,\ldots,X_n]}$
      canonically turns $S\str{[X_1,\ldots,X_n]}$ into an $S$-*algebra.
    \item
      It is easy to see that $S\str{[X_1,\ldots,X_n]}$ has the following explicit description when viewed
      as an internal ring: Elements are internal *finite $S$-linear combinations of \emph{*monomials in $n$ unknowns},
      i.e. of internal products of the form $X_1^{d_1}\cdot\ldots\cdot X_n^{d_n}$ with exponents $d_i\in\str{\N_0}$.
      These elements are added and multiplied in the obvious way.
    \item\label{bempoluniv}
      Transfer immediately shows that $S\str{[X_1,\ldots,X_n]}$ has the following universal property:
      If $T$ is an $S$-*algebra
      and if $(t_1,\ldots,t_n)$ is an internal family of elements of $T$,
      then there is a unique morphism of $S$-*algebras from $S\str{[X_1,\ldots,X_n]}$ to $T$
      which sends $X_i$ to $t_i$ for all $i$.
    \item\label{bempolinj}
      Let $n$ be a \emph{finite} natural number. Then by the universal property of usual polynomial rings,
      we have a canonical morphism of $S$-algebras (\emph{not} $S$-*algebras)
      $S[X_1,\ldots,X_n]\rightarrow S\str{[X_1,\ldots,X_n]}$ which sends $X_i$ to $X_i$.
      This map is easily seen to be injective,
      but is is (for $n\geq 1$) \emph{not} bijective:
      For example, for an infinite $h\in\str{\N_0}$,
      the monomial $X_1^h$ is obviously not contained in the image.
  \end{enumerate}
\end{bem}

\vspace{\abstand}

\begin{defi}\label{defA}
  Let $X$ be a *scheme, and let $n$ be a *natural number.
  We define the \emph{$n$-dimensional *affine space over $X$}
  as the $X$-*scheme $X\str{\times}_{\str{\Z}}\str{\Z}\str{[X_1,\ldots,X_n]}$
  (note that the fibre product exists by transfer).
\end{defi}

\vspace{\abstand}

\begin{bem}\label{bemP}
  For every scheme $X$ and every natural number $n\in\N_0$,
  we have the finitely presented $X$-scheme $\P^n_X=\P^n_{\Z}\times_{\Z}X$,
  the $n$-dimensional projective space over $X$,
  which is covered by $(n+1)$ copies of $\A^n_X$, glued together by certain universal morphisms.

  By transfer, for every *scheme $X$ and every *natural number $n\in\str{\N_0}$,
  we get a *finitely presented $X$-scheme $\str{\P^n_X}$,
  covered by $(n+1)$ copies of $\str{\A^n_X}$, the \emph{$n$-dimensional *projective space over $X$}.

  If $n$ is \emph{finite}, then these *affine spaces are glued together by the enlargements of
  the corresponding morphisms from the standard world.
\end{bem}

\vspace{\abstand}

\begin{defi}\label{defideal}
  If $S$ is a ring in $\cR$,
  and if $\{f_1,\ldots,f_m\}\subseteq S$ is a finite set of elements,
  then the category of $S$-algebras $A\in\Ob{\cR}$ with $f_1=\ldots=f_m=0\in A$ has an initial object,
  namely the $S$-algebra $S/(f_1,\ldots,f_m)$ (which is obviously finitely presented).

  It follows by transfer that for every *ring $S$ and any *finite internal subset
  $\{f_1,\ldots,f_m\}\subseteq S$,
  there is a $S$-*algebra $S/\str{(f_1,\ldots,f_m)}$ which is initial in the category of $S$-*algebras in which
  the $f_i$ are mapped to zero.
  --- We call $S/\str{(f_1,\ldots,f_m)}$
  the \emph{*factor ring of $S$ with respect to the *ideal $\str{(f_1,\ldots,f_m)}$}.
  \footnote[2]{By transfer, it is obvious that a *ideal of a *ring $S$ is in particular an ideal of $S$.}
\end{defi}

\vspace{\abstand}

\begin{bem}\label{bemideal}
  Let $S$ be a *ring, and let $(f_1,\ldots,f_m)$ be an ideal of $S$ with $m$ \emph{finite}.
  Then it follows by easy transfer that
  \[
    \str{(f_1,\ldots,f_m)}
    =(f_1,\ldots f_m)\cdot S
    \subseteq S.
  \]
\end{bem}

\vspace{\abstand}


\section{Relating schemes and *schemes}

\vspace{\abstand}

Let $A$ be a *ring in $\scB$.
On the one hand, when considering $A$ simply as a ring, we have the category $\Schfp{A}$
of finitely presented $A$-schemes.
On the other hand, we have the category $\sSchfp{A}$ of *finitely presented *schemes over $A$.

Intuitively, every finitely presented $A$-scheme determines
a *finitely presented $A$-*scheme which is "defined by the same relations",
and every morphism between finitely presented $A$-schemes gives a morphism between the associated $A$-*schemes.

In this section, we want to make this intuition precise by defining a morphism
$\s{}:\Schfp{\scB}\rightarrow\sSchfp{\cB}$ of fibrations over $\op{\scB}$.
In particular, by restricting to the fibre over $A$, this then gives us
the desired functor $\Schfp{A}\rightarrow\sSchfp{A}$.

\vspace{\abstand}

\begin{lemma}\label{lemmastarsquare}
  Let $\varphi:A\rightarrow B$ be a ring homomorphism in $\cR$.
  Then the diagram
  \begin{equation}\label{eqstsq}
    \xymatrix{
      {A} \ar[r]^{\varphi} \ar@{^{(}->}[d]_{*} & {B} \ar@{^{(}->}[d]^{*} \\
      {\str{A}} \ar[r]_{\str{\varphi}} & {\str{B}} \\
    }
  \end{equation}
  commutes in $\Rings$.
\end{lemma}

\vspace{\abstand}

\begin{proof}
  This follows immediately from elementary properties of enlargements.
\end{proof}

\vspace{\abstand}

\begin{satzdefi}\label{starfinpres}
  Let $A$ be an object of $\cR$.
  \begin{enumerate}
    \item\label{starfinpresi}
      Let $B=A[X_1,\ldots,X_n]/(f_1,\ldots,f_m)$ be a finitely presented $A$-algebra.
      Then
      \[
        \str{B}=\str{A}\str{[X_1,\ldots,X_n]}/\str{(f_1,\ldots,f_m)}.
      \]
    \item\label{starfinpresii}
      Let $\Alg{A}$ respectively $\Algfp{A}$ denote the category of $A$-algebras respectively
      finitely presented $A$-algebras.
      The canonical functors
      \[
        \begin{array}{ccccc}
          \op{\bigl(\Algfp{A}\bigr)} & \times & \sAlg{\str{A}} & \longrightarrow & \Sets \\[3mm]
          (B & , & C) & \mapsto &
          \left\{\begin{array}{l}
            \Mor{\sAlg{\str{A}}}{\str{B}}{C} \\[3mm]
            \Mor{\Alg{A}}{B}{C},
          \end{array}\right.
        \end{array}
      \]
      induced by
      $*:\Algfp{A}\rightarrow\sAlg{\str{A}}$ and the forgetful functor $\sAlg{\str{A}}\rightarrow\Alg{A}$,
      are canonically isomorphic via
      \[
        \bs{A}{B}{C}:\Mor{\sAlg{\str{A}}}{\str{B}}{C}
        \longrightarrow
        \Mor{\Alg{A}}{B}{C},\;\;
        \bigl[\str{B}\xrightarrow{\varphi}C\bigr]\mapsto
        \bigl[B\xrightarrow{*}\str{B}\xrightarrow{\varphi}C\bigr].
      \]
  \end{enumerate}
\end{satzdefi}

\vspace{\abstand}

\begin{proof}
  By transfer, \ref{bemstarpol}\ref{bempoluniv} and \ref{defideal},
  both $\str{B}$ and $\str{A}\str{[X_1,\ldots,X_n]}/\str{(f_1,\ldots,f_m)}$ have the same universal property
  in the category of $\str{A}$-algebras, which proves \ref{starfinpresi}.

  To show \ref{starfinpresii},
  we must first check that $\bs{A}{B}{C}$ is indeed functorial in the arguments $B$ and $C$.
  For argument $C$ this is trivial,
  and for argument $B$ it follows immediately from \ref{lemmastarsquare}.

  To see that $\bs{A}{B}{C}$ is a bijection,
  let $B=A[X_1,\ldots,X_n]/(f_1,\ldots,f_m)$. Then
  \begin{multline*}
    \Mor{\sAlg{\str{A}}}{\str{B}}{C}
    \stackrel{\ref{starfinpresi}}{=}
    \Mor{\sAlg{\str{A}}}{\str{A}\str{[X_1,\ldots,X_n]}/\str{(f_1,\ldots,f_m)}}{C} \\
    \stackrel{\ref{bemstarpol}\ref{bempoluniv}, \ref{defideal}}{=}
    \Bigl\{(c_1,\ldots,c_n)\in C^n\Bigl\vert\forall i\in\{1,\ldots,m\}:\, f_i(c_1,\ldots,c_n)=0\in C\Bigr\}
    =\Mor{\Alg{A}}{B}{C},
  \end{multline*}
  where this identification of the two sets is obviously just given by $\bs{A}{B}{C}$.
\end{proof}

\vspace{\abstand}

\begin{defi}\label{defT}
  For every ring $A$ in $\cR$,
  base change along the (external) ring homomorphism $*:A\rightarrow\str{A}$
  defines a functor $\T{}:\Sch{A}\rightarrow\Sch{\str{A}}$ (which respects schemes of finite presentation),
  and if $\varphi:A\rightarrow A'$ is a ring homomorphism,
  the diagram
  \[
    \xymatrix{
      {\Sch{A'}} \ar[d]_{\T{}} & {\Sch{A}} \ar[l]_{\varphi^*} \ar[d]^{\T{}} \\
      {\Sch{\str{A'}}} & {\Sch{\str{A}}} \ar[l]^{{(\str{\varphi})}^*} \\
    }
  \]
  commutes because of \ref{lemmastarsquare}.
  Consequently, we get "base change"-functors $\T{}$ of fibrations
  \[
    \xymatrix{
      {\Schfp{\cR}} \ar[r]^{\T{}} \ar@{^{(}->}[d] & {\Schfp{\scR}} \ar@{^{(}->}[d] \\
      {\Sch{\cR}} \ar[r]^{\T{}} \ar[d] & {\Sch{\scR}} \ar[d] \\
      {\op{\cR}} \ar[r]_{*} & {\op{\scR}.} \\
    }
  \]
  For every ring $A$ in $\cR$,
  base change along $\spec{\str{A}}\xrightarrow{\spec{*}}\spec{A}$
  defines for every $A$-scheme $X$ a morphism $\tm{X}:\T{X}\rightarrow X$ of schemes
  which is clearly functorial, i.e. the $\tm{X}$ define a 2-morphism $\tm{}$ of fibrations as follows:
  \begin{center}
    \includegraphics[scale=1.5]{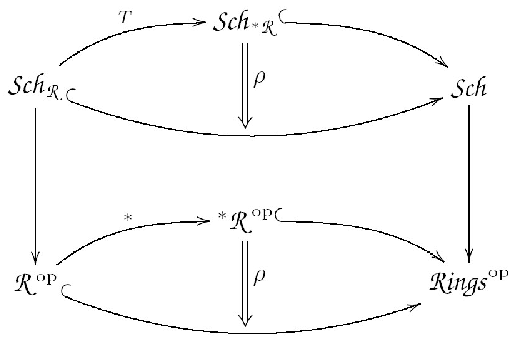}
  \end{center}
\end{defi}

\vspace{\abstand}

\begin{thm}\label{thmS}
  There is an essentially unique functor $\s{}:\Schfp{\scR}\rightarrow\sSchfp{\cR}$
  of fibrations over $\op{\scR}$ such that the following diagram of fibrations commutes:
  \begin{equation}\label{eq1}
    \includegraphics[scale=0.33]{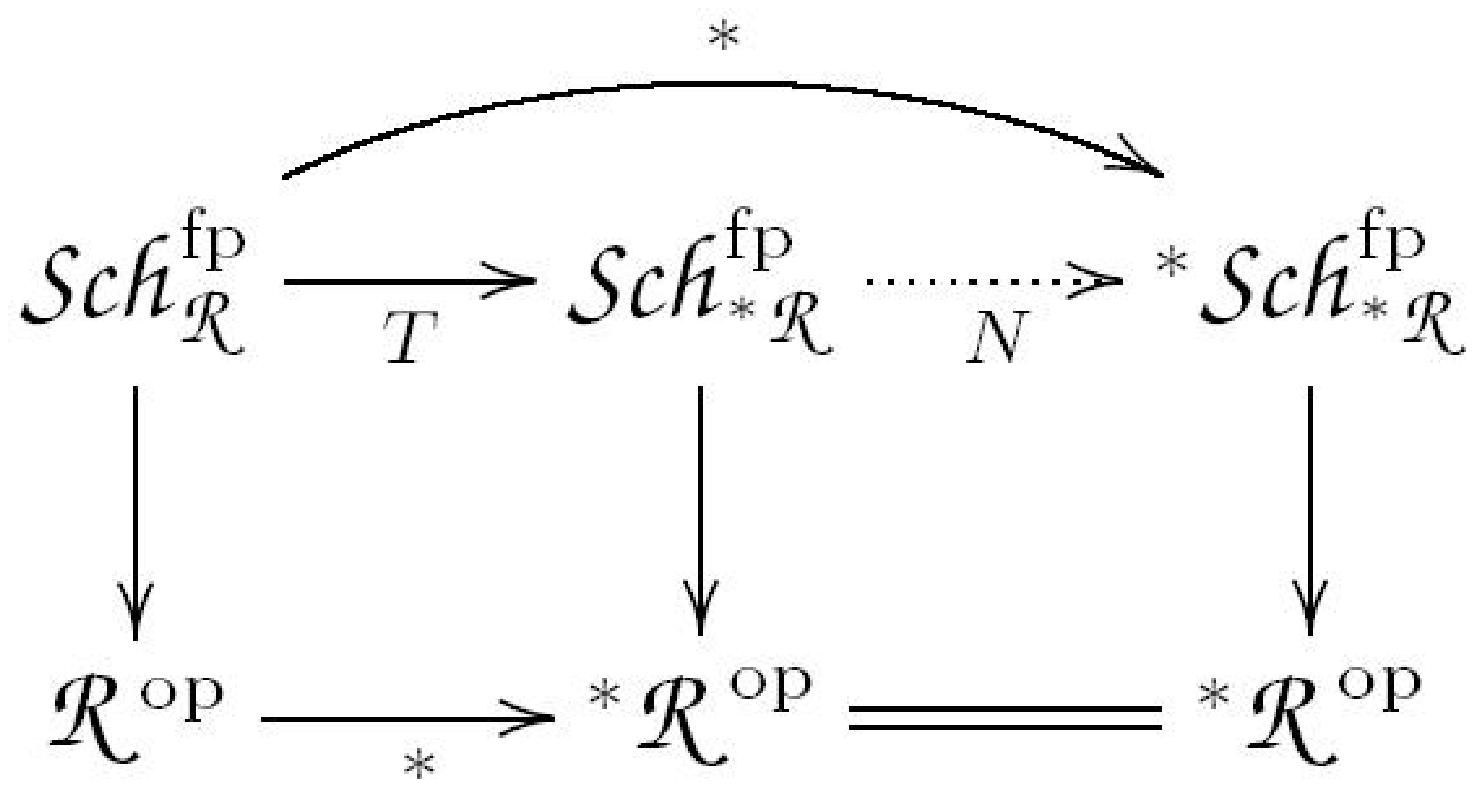}
  \end{equation}
  In particular, by restriction to $\scB$, we get a canonical functor
  $\s{}:\Schfp{\scB}\rightarrow\sSchfp{\cB}$ of fibrations over $\op{\scB}$.
\end{thm}

\vspace{\abstand}

\begin{proof}
  Let $A$ be a *ring,
  and let $X$ be a scheme of finite presentation over $A$.
  According to \cite[8.9.1]{ega43},
  there exist a subring $A_0\subseteq A$, finitely generated over $\Z$,
  and a finitely generated (and hence finitely presented) $A_0$-scheme $X_0$,
  such that $X_0\times_{A_0}A$ is isomorphic to $X$ over $A$.

  So $A_0$ is an object of $\cR$, and $X_0/A_0$ is an object of $\Schfp{A_0}$.
  According to \ref{starfinpres}\ref{starfinpresii},
  we get the following cartesian diagram of schemes:
  \begin{equation}\label{eq2}
    \includegraphics[scale=0.5]{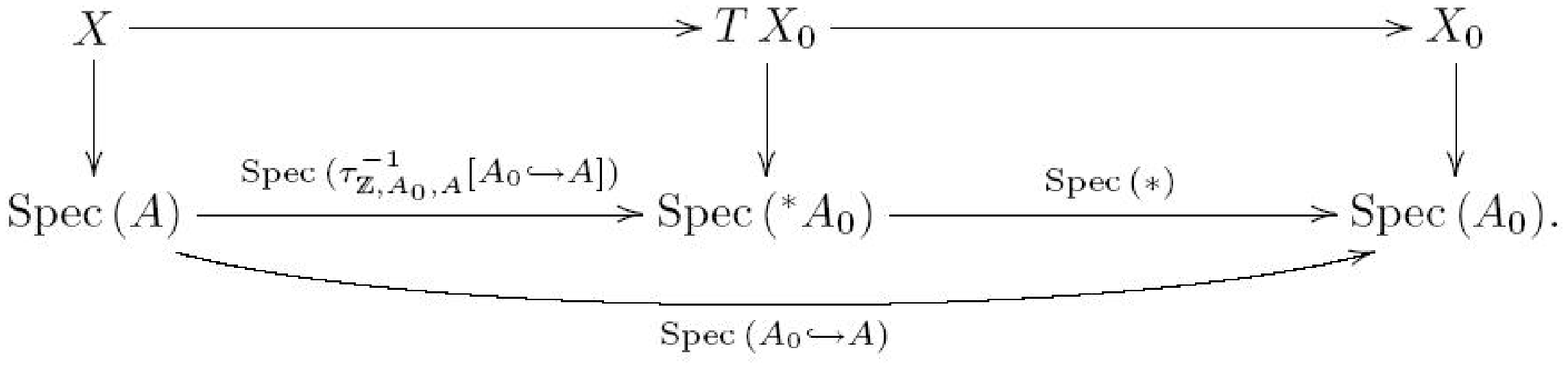}
  \end{equation}
  Therefore, in order to get a morphism of fibrations that makes \eqref{eq1} commute,
  we must define
  \[
    \s{X}
    \stackrel{\eqref{eq2}}{:=}
    {\Bigl(\bsi{\Z}{A_0}{A}[A_0\hookrightarrow A]\Bigr)}^*(\s{\T{X_0}})
    \stackrel{\eqref{eq1}}{:=}
    {\Bigl(\bsi{\Z}{A_0}{A}[A_0\hookrightarrow A]\Bigr)}^*(\str{X_0}).
  \]
  Now let $Y/S$ be another scheme of finite presentation,
  and let $f:X\rightarrow Y$ be an $S$-morphism.
  As before, there is a finitely generated ring $B_0\subseteq A$ and a finitely presented $B_0$-scheme $Y_0$
  such that $Y\cong Y_0\times_{B_0}A$.

  Let $\cI$ be the partially ordered set of finitely generated subrings of $A$ containing both $A_0$ and $B_0$,
  and put $X_C:=X_0\times_{A_0}C$ and $Y_C:=Y_0\times_{B_0}C$ for $C\in\cI$.

  Then $A=\varinjlim_{C\in\cI}C$, $X=\varprojlim_{C\in\cI}X_C$ and $Y=\varprojlim_{C\in\cI}Y_C$,
  and by \cite[8.8.2]{ega43} we have
  \begin{equation}\label{eq3}
    \varinjlim_{C\in\cI}\Mor{\Sch{C}}{X_C}{Y_C}=\Mor{\Sch{A}}{X}{Y}.
  \end{equation}
  In particular, there exists a $C_0\in\cI$ and a $C_0$-morphism $f_0:X_{C_0}\rightarrow Y_{C_0}$
  such that $f=f_0\times 1_A$.
  Therefore we get the following cartesian diagram of schemes
  \begin{center}
    \includegraphics[scale=0.5]{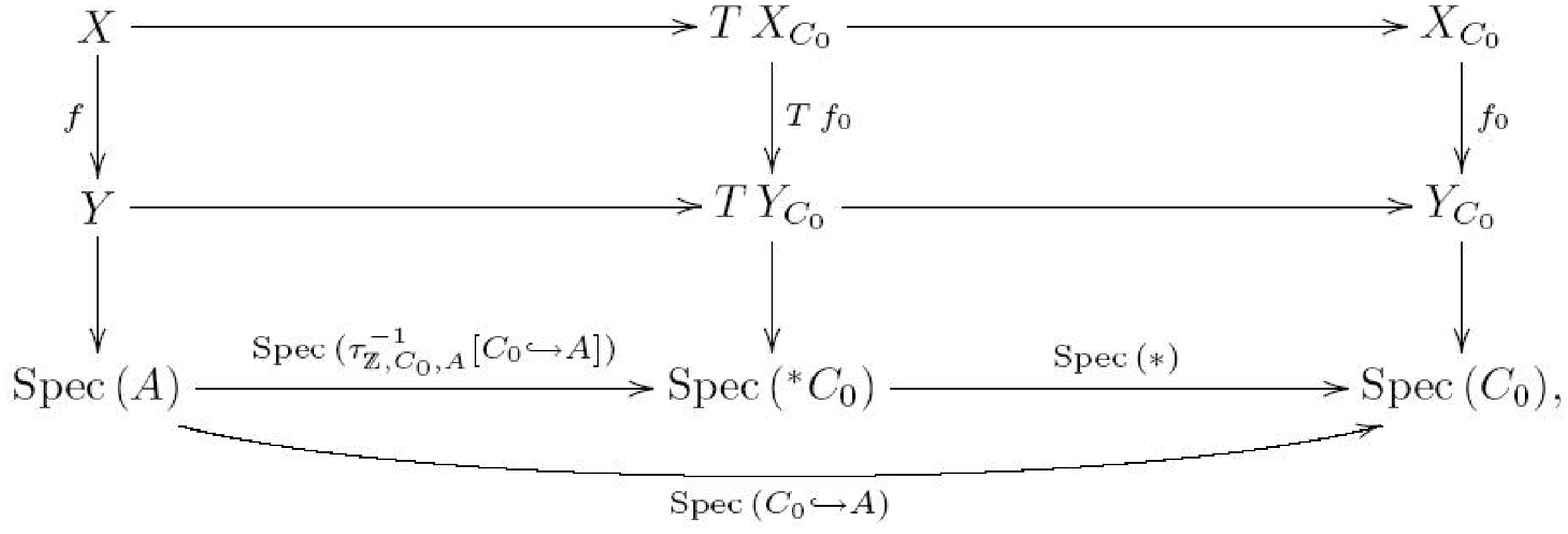}
  \end{center}
  and we are forced to set
  \[
    \s{f}
    \stackrel{\eqref{eq2}}{:=}
    {\Bigl(\bsi{\Z}{C_0}{A}[C_0\hookrightarrow A]\Bigr)}^*(\s{\T{f_0}})
    \stackrel{\eqref{eq1}}{:=}
    {\Bigl(\bsi{\Z}{C_0}{A}[C_0\hookrightarrow A]\Bigr)}^*(\str{f_0}).
  \]
  To check that this is well defined, let
  $C_1\in\cI$ be another subring of $A$ that admits a $C_1$-morphism $f_1:X_{C_1}\rightarrow Y_{C_1}$
  with $f=f_1\times 1_A$.

  Using \eqref{eq3} again, we find a subring $C_2$ of $A$ containing both $C_0$ and $C_1$
  with $f_0\times 1_{C_2}=f_1\times 1_{C_2}:X_{C_2}\rightarrow Y_{C_2}$,
  and \ref{starfinpres}\ref{starfinpresii} implies that the diagram
  \begin{center}
    \includegraphics[scale=0.5]{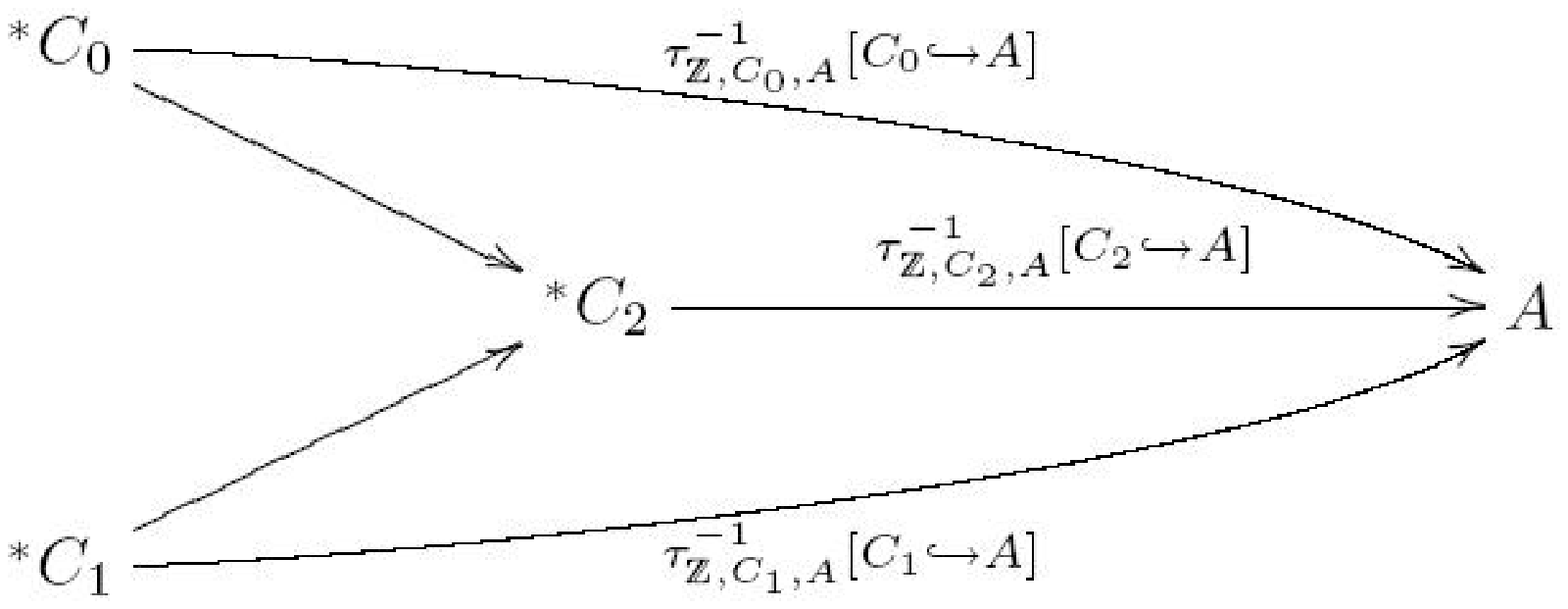}
  \end{center}
  commutes. Therefore we have
  \begin{multline*}
    {\Bigl(\bsi{\Z}{C_0}{A}[C_0\hookrightarrow A]\Bigr)}^*(\str{f_0})
    ={\Bigl(\bsi{\Z}{C_2}{A}[C_2\hookrightarrow A]\Bigr)}^*(\str{f_0\times 1_{C_2}}) \\
    ={\Bigl(\bsi{\Z}{C_2}{A}[C_2\hookrightarrow A]\Bigr)}^*(\str{f_1\times 1_{C_2}})
    ={\Bigl(\bsi{\Z}{C_1}{A}[C_1\hookrightarrow A]\Bigr)}^*(\str{f_1}).
  \end{multline*}
  Thus $\s{f}$ is well defined, and since this definition is obviously functorial,
  we get a functor $\s{}:\Schfp{A}\rightarrow\sSchfp{A}$
  which furthermore is uniquely determined (up to isomorphism) by the conditions stated in the theorem.

  It remains to show that this functor is compatible with inverse images and hence defines
  a morphism of fibrations $\s{}:\Schfp{\scR}\rightarrow\sSchfp{\cR}$ as claimed:
  If $\varphi:A\rightarrow A'$ is any morphism of *rings,
  we have to show that $\s{\varphi^* f}=\varphi^*\s{f}$ (for $f:X\rightarrow Y$ as above).
  With $D_0:=\varphi(C_0)\subseteq A'$ we have
  \[
    \varphi^*f=(f_0\times 1_A)\times 1_{A'}
    =f_0\times 1_{A'}
    =(f_0\times 1_{D_0})\times 1_{A'},
  \]
  so
  \begin{multline*}
    \s{\varphi^*f}
    ={\Bigl(\bsi{\Z}{D_0}{A'}[D_0\hookrightarrow A']\Bigr)}^*\str{(f_0\times 1_{D_0})}
    ={\Bigl(\bsi{\Z}{D_0}{A'}[D_0\hookrightarrow A']\Bigr)}^*\Bigl(\str{C_0}\rightarrow\str{D_0}\Bigr)^*(\str{f_0}) \\
    =\Bigl(\str{C_0}\rightarrow\str{D_0}\xrightarrow{\bsi{Z}{D_0}{A'}[D_0\hookrightarrow A']}A'\Bigr)^*(\str{f_0})
    \stackrel{\ref{starfinpres}\ref{starfinpresii}}{=}
    \Bigl(\bsi{Z}{C_0}{A'}[C_0\rightarrow D_0\xrightarrow{\varphi}A']\Bigr)^*(\str{f_0}) \\
    =\Bigl(\bsi{Z}{C_0}{A'}[C_0\hookrightarrow A\xrightarrow{\varphi}A']\Bigr)^*(\str{f_0})
    \stackrel{\ref{starfinpres}\ref{starfinpresii}}{=}
    \Bigl(\str{C_0}\xrightarrow{\bsi{Z}{C_0}{A}[C_0\hookrightarrow A]}A\xrightarrow{\varphi}A'\Bigr)^*(\str{f_0}) \\
    =\varphi^*\Bigl(\bsi{Z}{C_0}{A}[C_0\hookrightarrow A]\Bigr)^*(\str{f_0})
    =\varphi^*\s{f}.
  \end{multline*}
\end{proof}

\vspace{\abstand}

\begin{bem}\label{bemKan}
  The uniqueness of $\s{}$ in \ref{thmS} can be made precise as follows:
  It is easy to see that $\s{}$ is a right Kan extension of * along $\T{}$
  in the 2-category of fibrations (compare \cite[XII.4]{maclane_cat}),
  therefore enjoys a universal property
  and consequently is uniquely determined up to a canonical 2-isomorphism between morphisms of fibrations.
\end{bem}

\vspace{\abstand}

\begin{bsp}\label{exS}
  Let $A$ be a *ring,
  and let $B=A[X_1,\ldots,X_n]/(f_1,\ldots,f_m)$ be a finitely presented $A$-algebra.
  Let $A_0$ be the subring of $A$ generated by the (finitely many) coefficients of the $f_i$.
  Then we can consider the $f_i$ as elements of $A_0[X_1,\ldots,X_n]$,
  and we have $B=A_0[X_1,\ldots,X_n]/(f_1,\ldots,f_m)\otimes_{A_0}A$.
  Hence
  \begin{multline}\label{eqSaffine}
    \s{\spec{B}}
    =\bigl(\bsi{\Z}{A_0}{A}[A_0\subseteq A]\bigr)^*
      \Bigl[\str{\Bigl(\spec{A_0[X_1,\ldots,X_n]/(f_1,\ldots,f_m)}\Bigr)}\Bigr] \\
    \stackrel{\ref{starfinpres}\ref{starfinpresi}}{=}
      \bigl(\bsi{\Z}{A_0}{A}[A_0\subseteq A]\bigr)^*
      \Bigl[\str{\spec{\str{A_0}\str{[X_1,\ldots,X_n]}/\str{(f_1,\ldots,f_m)}}}\Bigr] \\
    \stackrel{\text{\tiny{transfer}}}{=}
      \str{\spec{A\str{[X_1,\ldots,X_n]}/\str{(f_1,\ldots,f_m)}}}.
  \end{multline}
  In particular, for $n\in\N_0$ we get $\s{\A^n_A}=\str{\A^n_A}$
  and --- taking $n=0$ --- $\s{\spec{A}}=\str{\spec{A}}$.
\end{bsp}

\vspace{\abstand}

\begin{satz}\label{satzProj}
  Let $A$ be a *ring,
  let $X$ be a finitely presented $A$-scheme,
  and let $n\in\N_0$ be a natural number.
  Then
  \[
    \s{\bigl(\P^n_X\xrightarrow{\text{\tiny can}}X\bigr)}
    =\str{\P^n_{\s{X}}}\xrightarrow{\text{\tiny can}}\s{X}.
  \]
\end{satz}

\vspace{\abstand}

\begin{proof}
  We know from the proof of \ref{thmS} that there exist
  a finitely generated subring $A_0$ of $A$ and a finitely presented $A_0$-scheme $X_0$
  with $X=X_0\times_{A_0}A$.
  Then
  \[
    \bigl(\P^n_X\xrightarrow{\text{\tiny can}}X\bigr)
    =\bigl(\P^n_{X_0}\xrightarrow{\text{\tiny can}}X_0\bigr)\times 1_A,
  \]
  and
  \[
    \s{\bigl(\P^n_X\xrightarrow{\text{\tiny can}}X\bigr)}
    ={\Bigl(\bsi{\Z}{A_0}{A}[A_0\hookrightarrow A]\Bigr)}^*
      \bigl[\str{\P^n_{\str{X_0}}}\xrightarrow{\text{\tiny can}}\str{X_0}\bigr]
    =\str{\P^n_{\s{X}}}\xrightarrow{\text{\tiny can}}\s{X}.
  \]
\end{proof}

\vspace{\abstand}


\section{Properties of the functor $\s{}$}

\vspace{\abstand}

\noindent
Let $A$ be a *ring in $\scB$.

\vspace{\abstand}

\begin{satz}\label{satzfin}
  The functor $\s{}:\Schfp{A}\rightarrow\sSchfp{A}$
  \begin{enumerate}
    \item\label{finleftex}
      is left exact, i.e. commutes with finite limits;
    \item\label{finglue}
      commutes with finite gluing data, i.e. if $\cI$ is a finite set,
      if $\coprod_{i,j\in\cI}U_{ij}\rightrightarrows\coprod_{i\in\cI}U_i$ with $U_{ij}$, $U_i$ finitely presented
      $A$-schemes is gluing data for an $A$-scheme $X$,
      then $\coprod_{i,j\in\cI}\s{U_{ij}}\rightrightarrows\coprod_{i\in\cI}\s{U_i}$ is gluing data for $\s{X}$;
    \item\label{finempty}
      sends the empty scheme to the empty *scheme;
    \item\label{finsum}
      commutes with finite sums.
  \end{enumerate}
\end{satz}

\vspace{\abstand}

\begin{proof}
  Let $\cI$ be a \emph{finite} category,
  and let $F:\cI\rightarrow\Schfp{A}$, $i\mapsto X^i$ be an arbitrary functor.
  According to \cite[8.8.3]{ega43},
  there exist a finitely generated subring $A_0$ of $A$
  and a functor $F_0:\cI\rightarrow\Schfp{A_0}$, $i\mapsto X^i_0$,
  such that $(\varprojlim_{i\in\cI}X^i_0)\times_{A_0}A=\varprojlim_{i\in\cI}X_i$.
  Since * is exact by \cite{enlcat}, and since inverse image functors in $\sSchfp{\cR}$ are left exact
  by transfer, we get $\s{\left(\varprojlim_{i\in\cI}X^i\right)}=\varprojlim_{i\in\cI}\s{X^i}$ by \ref{thmS}.
  Therefore \ref{finleftex} holds.

  Now let $\cI$ be a finite set, and let
  $\coprod_{i,j\in\cI}U_{ij}\rightrightarrows\coprod_{i\in\cI}U_i$ and $X$ be as in \ref{finglue}.
  By \cite[8.8.2, 8.10.5]{ega43},
  there are a finitely generated subring $A_0$ of $A$ and
  gluing data $\coprod_{i,j\in\cI}V_{ij}\rightrightarrows\coprod_{i\in\cI}V_i$,
  where the $V_{ij}$ and $V_i$ are finitely presented $A_0$-schemes
  and where base change with $A_0\hookrightarrow A$ gives back the original gluing data over $A$
  --- let $X_0$ be the finitely presented $A_0$-scheme defined gluing the $V_i$ along the $V_{ij}$.

  It follows from the construction of fibre products in \cite[3.2.6.3]{ega1}
  that base changes in the category of schemes respect gluing data.
  This implies firstly that $X_0\times_{A_0}A=X$
  and secondly (by transfer) that inverse image functors in $\sSchfp{\cR}$ commute with gluing data as well.
  Combining this with the exactness of * (note that "commuting with gluing data" means
  commuting with certain finite colimits) completes the proof of \ref{finglue} using the same reasoning as for
  \ref{finleftex}.

  Let $0$ denote the trivial ring, and let $\emptyset=\spec{0}$ be the empty (finitely presented) $A$-scheme.
  Then $\emptyset=[\str{\Z}\rightarrow A]^*\spec{0}$, so
  \[
    \s{\emptyset}
    \stackrel{\ref{thmS}}{=}[\str{\Z}\rightarrow A]^*\bigl(\s{\spec{0}}\bigr)
    \stackrel{\ref{exS}}{=}[\str{\Z}\rightarrow A]^*\bigl(\str{\spec{0}}\bigr)
    \stackrel{\text{\tiny{transfer}}}{=}\str{\spec{0}},
  \]
  which is the empty *scheme.

  Finally, \ref{finsum} is just the special case of \ref{finglue} where all the $U_{ij}$ are empty,
  and combining \ref{finglue} with \ref{finempty} immediately finishes the proof.
\end{proof}

\vspace{\abstand}

\begin{bem}
  Combining \ref{exS} with \ref{satzfin}\ref{finglue} provides us with an alternative description of the functor
  $\s{}$, at least when we restrict our attention to \emph{separated} $A$-schemes of finite presentation:

  Every finitely presented $A$-scheme $X$ admits a finite open affine covering $X=\bigcup_{i\in\cI}U_i$,
  and if $X/A$ is separated, the intersections $U_{ij}:=U_i\cap U_j$ are affine as well
  by \cite[5.5.6]{ega1}.
  So in this case, we can compute the $\s{U_{ij}}$ and $\s{U_i}$ using \ref{exS},
  and we know from \ref{satzfin}\ref{finglue} that $\s{X}$ is obtained by glueing the $\s{U_i}$ along
  the $\s{U_{ij}}$.
\end{bem}

\vspace{\abstand}

\begin{cor}\label{corSgroup}
  Let $G$ be a finitely presented (commutative) $A$-group scheme.
  Then $\s{G}$ is a *finitely presented (commutative) $A$-*group *scheme,
  i.e. a (commutative) group object in $\sSchfp{A}$.
\end{cor}

\vspace{\abstand}

\begin{proof}
  The data defining a (commutative) group scheme structure on $G$ can be expressed with diagrams involving
  only $A$, $G$, $G\times_AG$ and $G\times_AG\times_AG$,
  and these products are respected by $\s{}$ according to \ref{satzfin}\ref{finleftex}.
\end{proof}

\vspace{\abstand}

\begin{satz}\label{satzP}
  Let $f:X\rightarrow Y$ be a morphism of finitely presented $A$-schemes,
  and let $\bP$ be one of the following properties of morphisms of schemes:
  \begin{itemize}
    \item
      isomorphism,
    \item
      monomorphism,
    \item
      immersion,
    \item
      open immersion,
    \item
      closed immersion,
    \item
      separated,
    \item
      surjective,
    \item
      radicial,
    \item
      affine,
    \item
      quasi-affine,
    \item
      finite,
    \item
      quasi-finite,
    \item
      proper,
    \item
      projective,
    \item
      quasi-projective.
  \end{itemize}
  If $f$ has property $\bP$, then $\s{f}:\s{X}\rightarrow\s{Y}$ has property $\str{\bP}$.
\end{satz}

\vspace{\abstand}

\begin{proof}
  Let $\bP$ be one of the above properties.
  By \cite[8.8.2, 8.10.5]{ega43},
  there exist a finitely generated ring $A_0\subseteq A$
  and a morphism $f_0:X_0\rightarrow Y_0$ of finitely presented $A_0$-schemes
  such that $X_0\times_{A_0}A=X$, $Y_0\times_{A_0}A=Y$, $f_0\times 1_A=f$
  and such that $f_0$ has property $\bP$.

  Then $\str{f_0}:\str{X_0}\rightarrow\str{Y_0}$ has property $\str{\bP}$,
  and since property $\str{\bP}$ is stable under base change
  (by transfer, because $\bP$ is stable under base change),
  we see that
  $\s{f}\stackrel{\ref{thmS}}{=}\bigl(\bsi{\Z}{A_0}{A}[A_0\subseteq A]\bigr)^*(\str{f_0})$
  has property $\str{\bP}$ as well.
\end{proof}

\vspace{\abstand}

\begin{bem}\label{bemQC}
  Let $X$ be a finitely presented $A$-scheme, and let $U\subseteq X$ be an open subscheme.
  According to \cite[1.6.2(i),(v)]{ega41}, $U$ is a finitely presented $A$-scheme if and only if $U$ is
  quasi-compact. It follows that $\s{U}$ is defined (and then a *open *subscheme of $\s{X}$ by \ref{satzP})
  if and only if $U$ is quasi-compact.

  Note that the quasi-compact open subsets of $X$ form a basis for the Zariski topology (since affine
  open sets are quasi-compact), so that there will be no harm in restricting our attention to
  quasi-compact open subschemes.
\end{bem}

\vspace{\abstand}

\begin{cor}\label{corinvim}
  Let $f:X\rightarrow Y$ be a morphism of finitely presented $A$-schemes,
  and let $U\subseteq Y$ be a quasi-compact open subscheme of $Y$.
  Then $\s{U}$ is an open *subscheme of $\s{Y}$, and
  \[
    \s{\bigl(f|_{f^{-1}(U)}\bigr)}
    =(\s{f})|_{(\s{f})^{-1}(\s{U})}
    \in\Mor{\sSchfp{A}}{\s{(f^{-1}(U))}}{\s{U}}
  \]
\end{cor}

\vspace{\abstand}

\begin{proof}
  This follows immediately from the fact that $\s{}$ is left exact by \ref{satzfin}\ref{finleftex}
  and respects open immersions by \ref{satzP},
  applied to the cartesian diagram
  \[
    \xymatrix@C=20mm{
      {f^{-1}(U)} \ar@{^{(}->}[d] \ar[r]^{f|_{f^{-1}(U)}} & {U} \ar@{^{(}->}[d] \\
      {X} \ar[r]_{f} & {Y.} \\
    }
  \]
\end{proof}

\vspace{\abstand}

\begin{cor}\label{corcov}
  Let $X=\bigcup_{i\in\cI}U_i$ be a \emph{finite} (affine) covering by quasi-compact open subschemes.
  Then $\s{X}=\bigcup_{i\in\cI}\s{U_i}$ is a *open (*affine) *covering in $\sSchfp{A}$.
\end{cor}

\vspace{\abstand}

\begin{proof}
  If the $U_i$ are affine, the $\s{U_i}$ are *affine by example \ref{exS}.
  The $\s{U_i}$ are open subschemes of $\s{X}$ by \ref{satzP},
  and since
  \[
    \coprod_{i\in\cI}\s{U_i}
    \stackrel{\ref{satzfin}}{=}\s{\bigl(\coprod_{i\in\cI}U_i\bigr)}
    \stackrel{\ref{satzP}}{\twoheadrightarrow}
    \s{X},
  \]
  they cover $\s{X}$.
\end{proof}

\vspace{\abstand}

\begin{lemma}\label{lemmaComplement}
  Let $X$ be a finitely presented $A$-scheme,
  let $Y\subseteq X$ be a closed, finitely presented subscheme,
  and assume that the open complement $U:=X\setminus Y$ is quasi-compact.
  Then $\s{U}$ is $[\s{X}]\setminus[\s{Y}]$, the *complement of $\s{Y}$ in $\s{X}$.
\end{lemma}

\vspace{\abstand}

\begin{proof}
  Since the diagram
  \[
    \xymatrix{
      {\emptyset} \ar@{^{(}->}[r] \ar@{^{(}->}[d] \ar@{}[rd]|\Box &
      U \ar@{^{(}->}[d] \\
      Y \ar@{^{(}->}[r] &
      X \\
    }
  \]
  is cartesian, \ref{satzfin}\ref{finleftex}, \ref{finempty} imply that
  \[
    \xymatrix{
      {\emptyset} \ar@{^{(}->}[r] \ar@{^{(}->}[d] \ar@{}[rd]|\Box &
      {\s{U}} \ar@{^{(}->}[d] \\
      {\s{Y}} \ar@{^{(}->}[r] &
      {\s{X}} \\
    }
  \]
  is also cartesian, i.e. $\s{U}$ lies in $[\s{X}]\setminus[\s{Y}]$.
  For the other inclusion,
  note that the surjectivity of $Y\coprod U\rightarrow X$ implies the *surjectivity of
  $[\s{Y}]\coprod[\s{U}]\rightarrow\s{X}$ by \ref{satzfin}\ref{finsum} and \ref{satzP}.
\end{proof}

\vspace{\abstand}

Let $\varphi:R\rightarrow S$ be a ring homomorphism,
let $X$ be an $R$-scheme, and let $Y$ be an $S$-scheme.
Then it is common practice to simply write $X(Y)$ for the set of those morphisms $f:Y\rightarrow X$ of schemes
that make the diagram
\[
  \xymatrix@C=15mm{
    {Y} \ar[r]^{f} \ar[d] & {X} \ar[d] \\
    {\spec{S}} \ar[r]_{\spec{\varphi}} & {\spec{R}} \\
  }
\]
commute, thus dropping $R$, $S$ and $\varphi$ from the notation.
In other words,
when $R$, $S$ and $\varphi$ are understood,
$X(Y)$ denotes the subset of those morphisms in $\SCH$ which project to $\varphi$
in the bifibration $\SCH\rightarrow\op{\Rings}$.

In analogy to this practice, we make the following definition:

\vspace{\abstand}

\begin{defi}\label{defpoint}
  Let $X$ be a *scheme in $\sSchfp{A}$,
  let $\varphi:A\rightarrow B$ be a morphism of *rings,
  and let $Y$ be a *scheme in $\sSchfp{B}$.

  Then we denote the set of those morphisms in $\Mor{\sSchfp{\cR}}{Y/B}{X/A}$
  which are projected to $\varphi$ under $\sSchfp{\cR}\rightarrow\op{\scR}$ by $X(Y)$
  and call it the set of \emph{$Y$-valued points of $X$}
  (where we assume that $A$, $B$ and $\varphi$ are understood).

  In the special case $Y=\str{\spec{B}}$, we put $X(B):=X(Y)$ and call $X(B)$
  the set of \emph{$B$-valued points of $X$}.
\end{defi}

\vspace{\abstand}

\begin{bem}\label{bempoint}
  Let $X$ be a finitely presented $A$-scheme,
  let $\varphi:A\rightarrow B$ be a morphism of *rings,
  and let $T$ be a finitely presented $B$-scheme.

  Then the functor $\s{}$ induces a canonical map
  \[
    \xymatrix@C=40mm{
      {X(T)} \ar@{.>}[r]^{\s{}} \ar@{^{(}->}[d] & {(\s{X})(\s{T})} \ar@{^{(}->}[d] \\
      {\Mor{\Schfp{\cB}}{T/B}{X/A}} \ar[r]_{\s{}} & {\Mor{\sSchfp{\cB}}{\s{T}/B}{\s{X}/A}} \\
    }
  \]
  (note that $\s{}$, restricted to $X(T)$, factorizes over $(\s{X})(\s{T})$,
  because $\s{}$ is a morphism of fibrations and hence in particular a morphism of categories over $\op{\scR}$).

  Since $\s{\spec{B}}=\str{\spec{B}}$ by \ref{exS},
  we in particular get a map $\s{}:X(B)\rightarrow(\s{X})(B)$
  from $B$-valued points of $X$ to $B$ valued points of $\s{X}$.
\end{bem}

\vspace{\abstand}

\begin{defi}\label{defSaff}
  As we have seen in \ref{exS}, the functor $\s{}:\Schfp{A}\rightarrow\sSchfp{A}$
  sends affine schemes to *affine schemes and thus induces a functor $\Algfp{A}\rightarrow\sAlgfp{A}$
  --- which we want to denote by $\s{}$ as well ---
  satisfying
  \begin{equation}\label{eqSaff}
    \forall B\in\Ob{\Algfp{A}}:\
    \s{\spec{B}}
    =\str{\spec{\s{B}}}.
  \end{equation}
  If $B=A[X_i]/(f_j)$,
  then we have calculated in \ref{exS} that
  $\s{B}=A\str{[X_i]}/\str{(f_j)}$.
  It follows from \ref{bemstarpol}\ref{bempoluniv} and \ref{defideal}
  that sending $X_i$ to $X_i$ defines a canonical morphism of $A$-algebras $\sm{B}:B\rightarrow\s{B}$,
  which is obviously functorial:
  If $\varphi:B\rightarrow C$ is a morphism of $A$-algebras,
  then
  \begin{equation}\label{eqSfunct}
    \xymatrix{
      {B} \ar[r]^{\sm{B}} \ar[d]_{\varphi} & {\s{B}} \ar[d]^{\s{\varphi}} \\
      {C} \ar[r]_{\sm{C}} & {\s{C}} \\
    }
  \end{equation}
  commutes in the category of $A$-algebras.
\end{defi}

\vspace{\abstand}

\begin{lemma}\label{lemmaBSaff}
  Let $k$ be an $A$-*algebra,
  and let $B$ be a finitely presented $A$-algebra.
  Then the canonical map
  \[
    (\sm{B})_*:\Mor{\sAlg{A}}{\s{B}}{k}\longrightarrow\Mor{\Alg{A}}{B}{k},\;\;\;
    [\s{B}\xrightarrow{\varphi}k]\mapsto
    [B\xrightarrow{\sm{B}}\s{B}\xrightarrow{\varphi}k]
  \]
  is bijective.
\end{lemma}

\vspace{\abstand}

\begin{proof}
  Let $B=A[X_i]/(f_j)$.
  We can argue as in the proof of \ref{starfinpres}\ref{starfinpresii}:
  A morphism $\varphi:\s{B}\rightarrow k$ in $\sAlg{A}$
  is precisely given by a tuple $(x_1,\ldots,x_n)\in k^n$ satisfying
  $f_j(x_1,\ldots,x_n)=0\in k$ for all $j$,
  and the exact same data defines a morphism $\varphi':B\rightarrow k$ of $A$-algebras.
  --- It is clear that this identification between the two sets of morphisms is
  just the one given in the lemma.
\end{proof}

\vspace{\abstand}

\begin{thm}\label{thmBS}
  Let $k$ be a *artinian $A$-*algebra,
  and let $X$ be a finitely presented $A$-scheme.
  Then the canonical map $\s{}:X(k)\rightarrow(\s{X})(k)$ is bijective.
\end{thm}

\vspace{\abstand}

\begin{proof}
  We choose a finite affine open covering
  $X=\bigcup_{i\in\cI}U_i$,
  so that $\s{X}=\bigcup_{i\in\cI}\s{U_i}$ is a *open *affine *covering of $\s{X}$ by \ref{corcov}.

  To prove surjectivity, let $f:\str{\spec{k}}\rightarrow X$ be an arbitrary $k$-valued point of $X$.
  By transfer, since $k$ is *artinian, $f$ factorizes over one of the $\s{U_i}$,
  so without loss of generality, we can assume that $X=\spec{B}$ is affine.

  Then $\s{X}\stackrel{\eqref{eqSaff}}{=}\str{\spec{\s{B}}}$,
  and $f$ corresponds to a morphism
  $\varphi:\s{B}\rightarrow k$ of $A$-*algebras
  which induces a morphism $\varphi':=\varphi\sm{B}:B\rightarrow k$ of $A$-algebras as in \ref{lemmaBSaff},
  hence a $k$-valued point $f':=\spec{\varphi'}$ of $X$.
  It is clear that $\s{f'}=f$, so $\s{}$ is indeed surjective.

  For injectivity, let $f,g\in X(k)$ be two $k$-valued points of $X$ with
  $\s{f}=\s{g}\in(\s{X})(k)$.

  If $X_k$ denotes the inverse image of $X$ under $A\rightarrow k$,
  then the canonical map $X_k(k)\rightarrow X(k)$ is a bijection, so that we can assume $A=k$ without
  loss of generality.
  As above, it follows that $f$ factorizes over one of the $U_i$, say over $U_{i_0}$ ---
  then $\s{f}$ factorizes over $\s{U_{i_0}}$.
  Let us assume that $g$ does not factorize over $U_{i_0}$.
  This would imply that the following diagram of finitely presented $k$-schemes is cartesian:
  \[
    \xymatrix@C=20mm{
      {\emptyset} \ar[r] \ar[d] & {U_{i_0}} \ar@{^{(}->}[d] \\
      {\spec{k}} \ar[r]_g & {X.} \\
    }
  \]
  Then \ref{satzfin}\ref{finleftex} and \ref{finempty} imply that
  \[
    \xymatrix@C=20mm{
      {\str{\emptyset}} \ar[r] \ar[d] & {\s{U_{i_0}}} \ar@{^{(}->}[d] \\
      {\str{\spec{k}}} \ar[r]_{\s{g}=\s{f}} & {\s{X}} \\
    }
  \]
  is cartesian as well, a contradiction to the fact that $\s{f}$ factorizes over $\s{U_{i_0}}$.

  Therefore both $f$ and $g$ factorize over $U_{i_0}$,
  and we can again assume that $X=\spec{B}$ is affine.
  But then $f$ and $g$ correspond to $k$-algebra morphisms $\varphi,\psi:B\rightarrow k$,
  and $\s{f}=\s{g}$ means that the induced morphisms of $k$-*algebras
  $\varphi',\psi':\s{B}\rightarrow k$ are the same.
  But then $\varphi$ and $\psi$ must be the same as well according to \ref{lemmaBSaff}.
\end{proof}

\vspace{\abstand}

\begin{satz}\label{satzSmooth}
  Let $f:X\rightarrow Y$ be a morphism of finitely presented $A$-schemes.
  If $f$ is étale (unramified, smooth), then $\s{f}:\s{X}\rightarrow\s{Y}$ is
  *étale (*unramified, *smooth).
\end{satz}

\vspace{\abstand}

\begin{proof}
  First consider the case where $f:X\rightarrow Y$ is unramified.
  By \cite[17.4.2]{ega44}, a morphism $f:X\rightarrow Y$ of (locally) finite presentation is unramified
  if and only if the diagonal
  $\Delta_{X/Y}:X\xrightarrow{(f,f)}X\times_YX$ is an open immersion.
  So in our case, $\Delta_{X/Y}$ is an open immersion, and
  \ref{satzfin}\ref{finleftex} and \ref{satzP} show that the *diagonal
  $\Delta_{\s{X}/\s{Y}}:\s{X}\xrightarrow{(\s{f},\s{f})}\s{X}\times_{\s{Y}}\s{X}$
  is a *open immersion, hence transferring \cite[17.4.2]{ega44} proves that $\s{f}$ is *unramified
  (since it is *finitely presented by construction).

  Now let $f:X\rightarrow Y$ be étale.
  By \cite[17.1.6]{ega44}, \ref{corinvim} and \ref{corcov},
  we can assume without loss of generality that $X$ and $Y$ are affine
  and that $f$ is given by a morphism $\varphi:B\rightarrow C$ of finitely presented $A$-algebras.
  Furthermore, by \cite[I.3.16]{milne},
  we can assume that $C=B[T_1,\ldots,T_n]/(P_1,\ldots,P_n)$ with
  $d:=\det(\partial P_i/\partial T_j)\in C^\times$ and that $\varphi$ is the canonical morphism,
  and we have to show that $\s{\varphi}:\s{B}\rightarrow\s{C}=(\s{B})\str{[T_i]}/\str{(P_j)}$
  is *étale.
  By transfer of \cite[I.3.16]{milne},
  for this it suffices to show that $d':=\str{\det}(\str{\partial P_i}/\str{\partial T_j})$
  is a *unit in $\s{C}$.

  Since partial derivatives of polynomials and determinants of matrices are given
  by universal polynomials in the coefficients,
  it follows easily that the diagrams
  \[
    \xymatrix@C=13mm{
      {B[T_1,\ldots,T_n]} \ar[r]^{\partial/\partial T_j} \ar[d]_{\sm{B[T_i]}} &
      {B[T_1,\ldots,T_n]} \ar[d]^{\sm{B[T_i]}} & &
      {B[T_1,\ldots,T_n]^{n\times n}} \ar[r]^{\det} \ar[d]_{\sigma_{B[T_i]}^{n\times n}} &
      {B[T_1,\ldots,T_n]} \ar[d]^{\sm{B[T_i]}} \\
      {B\str{[T_1,\ldots,T_n]}} \ar[r]_{\str{\partial}/\str{\partial T_j}} &
      {B\str{[T_1,\ldots,T_n]}} & &
      {B\str{[T_1,\ldots,T_n]}^{n\times n}} \ar[r]_{\str{\det}} &
      {B\str{[T_1,\ldots,T_n]}} \\
    }
  \]
  commute, which implies $d'=\sm{C}(d)\in\s{C}$.
  Since $\sm{C}$ is a ring homomorphism, it maps units to units, so $d'$ is a unit in $\s{C}$.
  But being a unit is obviously a first order property, so units and *units are the same thing,
  and we are done in the case where $f$ is étale.

  Finally, let $f:X\rightarrow Y$ be smooth.
  By \cite[3.24]{milne},
  this is equivalent to the existence of a (finite) open affine covering $U_i$ of $X$,
  such that for every $i$ the restriction $f|_{U_i}$ factorizes as
  \[
    \xymatrix{
      {U_i} \ar[r]^{f|_{U_i}} \ar[d]_{g_i} & {Y} \\
      {\A^n_{V_i}} \ar[r]_{\text{\tiny{can}}} & {V_i} \ar@{^{(}->}[u] \\
    }
  \]
  with $g_i$ étale and $n\in\N_0$. Since the functor $\s{}$ respects open affine coverings by \ref{corcov},
  restrictions by \ref{corinvim}, open immersions by \ref{satzP},
  affine spaces (over affine bases) by \ref{exS} and étale morphisms by the second part of the proof,
  transfer of \cite[3.24]{milne} shows that $\s{f}$ is indeed *smooth.
\end{proof}

\vspace{\abstand}

\begin{lemma}\label{lemmaSS}
  Let $B$ be a finitely presented $A$-algebra,
  and let $C=B[Y_1,\ldots,Y_k]/\cJ$ be a finitely presented $B$-algebra.
  Then $\s{C}=(\s{B})\str{[Y_j]}/\str{\cJ}$.
\end{lemma}

\vspace{\abstand}

\begin{proof}
  Let $B=A[X_1,\ldots,X_n]/\cI$ be a finite presentation of $B$ as an $A$-algebra.
  Then
  \[
    \s{C}
    =\s{\Bigl(A[X_i,Y_j]/(\cI+\cJ)\Bigr)}
    \stackrel{\ref{exS}}{=}A\str{[X_i,Y_j]}/\str{(\cI+\cJ)} \\
    \stackrel{\text{\tiny transfer}}{=}\Bigl(A\str{[X_i]}/\str{\cI}\Bigr)\str{[Y_j]}/\str{\cJ}
    \stackrel{\ref{exS}}{=}(\s{B})[Y_j]/\str{\cJ}.
  \]
\end{proof}

\vspace{\abstand}

\begin{satz}\label{satzfinite}
  Let $B$ be a finitely presented $A$-algebra,
  and let $C$ be a \emph{finite} $B$-algebra.
  Then the canonical ring homomorphism $C\otimes_B\s{B}\longrightarrow\s{C}$
  induced by \eqref{eqSfunct} is an isomorphism.
\end{satz}

\vspace{\abstand}

\begin{proof}
  First consider the case where $C=B/\cI$ is a quotient of $B$. Then
  \[
    C\otimes_B\s{B}
    =(\s{B})/\cI\cdot\s{B}
    =(\s{B})/\str{\cI}
    \stackrel{\ref{lemmaSS}}{=}\s{C}.
  \]
  Next let $C=B[c]/(c^n+b_{n-1}c^{n-1}+\ldots+b_0)$ with $n\in\N_+$ and $b_0,\ldots,b_{n-1}\in B$.
  Consider the following true statement in $\hat{M}$:
  \begin{quote}
    For every object $R$ of $\cR$ and for every tuple $(r_0,\ldots,r_{n-1})\in R^n$,
    sending $e_i$ to $\bar{X}^{i-1}$ defines an isomorphism of $R$-modules
    $R^n\xrightarrow{\sim}R[X]/(X^n+r_{n-1}X^{n-1}+\ldots+r_0)$.
  \end{quote}
  By transfer and the fact that an isomorphism of *modules is in particular an isomorphism of modules,
  we get:
  \begin{quote}
    For every *ring $R$ and for every tuple $(r_0,\ldots,r_{n-1})\in R^n$,
    sending $e_i$ to $\bar{X}^{i-1}$ defines an isomorphism of $R$-modules
    $R^n\xrightarrow{\sim}R\str{[X]}/\str{(X^n+r_{n-1}X^{n-1}+\ldots+r_0)}$.
  \end{quote}
  By \ref{lemmaSS}, we have $\s{C}=(\s{B})\str{[c]}/\str{(c^n+b_{n-1}c^{n-1}+\ldots+b_0)}$,
  so we get the following commutative diagram of $\s{B}$-modules:
  \[
    \xymatrix{
      {\bar{c}^{i-1}\otimes 1} \ar@{}[r]|{\in} & {C\otimes_B\s{B}} \ar@{.>}[r] &
        {\s{C}} \ar@{}[r]|{\ni} & {\bar{c}^{i-1}} \\
      {e_i} \ar@{|->}[u] \ar@{}[r]|{\in} & {\s{B}^n} \ar[u]^{\wr} \ar[r]_{\sim} &
        {\s{B}^n} \ar[u]_{\wr} \ar@{}[r]|{\ni} & {e_i,} \ar@{|->}[u] \\
    }
  \]
  and we are done in this case as well.

  Now let $C=B[c]/\cI$. Then the element $\bar{c}$ of $C$ is integral over $B$, because $C/B$ is finite,
  so there is a relation $\bar{c}^n+b_{n-1}\bar{c}^{n-1}+\ldots+b_0=0$ in $C$, which means that
  $B\rightarrow C$ factorizes as
  \[
    B\rightarrow
    \underbrace{B[c]/(c^n+b_{n-1}c^{n-1}+\ldots+b_0)}_{=:C'}
    \stackrel{c\mapsto\bar{c}}{\twoheadrightarrow}C,
  \]
  and we get
  \[
    \s{C}
    \stackrel{\text{\tiny 1. case}}{=}C\otimes_{C'}\s{C'}
    \stackrel{\text{\tiny 2. case}}{=}C\otimes_{C'}C'\otimes_B\s{B}
    =C\otimes_B\s{B}.
  \]
  Finally, in the general case, let $C=B[X_1,\ldots,X_n]/\cI$ for an $n\in\N_+$.
  We prove the proposition by induction on $n$:
  The case $n=1$ has been proven above, so let $C=B[X_1,\ldots,X_{n+1}]/\cI$ for $n\geq 1$.
  Let $C'$ be the subring of $C$ generated by $\bar{X}_1,\ldots,\bar{X}_n$ as a $B$-algebra.
  Then $C=C'[X_{n+1}]/\cJ$, and
  \[
    \s{C}
    \stackrel{\text{\tiny 3. case}}{=}C\otimes_{C'}\s{C'}
    \stackrel{\text{\tiny induction}}{=}C\otimes_{C'}C'\otimes_B\s{B}
    =C\otimes_B\s{B}.
  \]
\end{proof}

\vspace{\abstand}


\section{*Modules over *schemes}

\vspace{\abstand}

Let $\MOD$ be the category whose objects are pairs
$\langle\cF,X/A\rangle$,
consisting of an $A$-scheme $X$ and an $\cO_X$-module $\cG$,
and whose morphisms from $\langle\cF,X/A\rangle$ to $\langle\cG,Y/B\rangle$
are pairs $\langle\alpha,\langle f,\varphi\rangle\rangle$
with $\langle f,\varphi\rangle$ a morphism from $X/A$ to $Y/B$ in $\SCH$
and $\alpha:f^*\cG\rightarrow\cF$ a morphism of $\cO_X$-modules.

Projection onto the second component defines an abelian bifibration $\MOD\rightarrow\SCH$
(or $\MOD\rightarrow\op{\Rings}$ after composing with $\SCH\rightarrow\op{\Rings}$):
For a morphism $\langle f,\varphi\rangle :X/A\rightarrow Y/B$,
direct and inverse image functor are given by
$\langle f,\varphi\rangle_*\langle\cF,X/A\rangle=\langle f_*\cF,Y/B\rangle$ respectively
$\langle f,\varphi\rangle^*\langle\cG,Y/B\rangle=\langle f^*\cG,X/A\rangle$,
and the fibre over an object $X/A$ is \emph{the opposite} of the category $\Mod{X}$
of $\cO_X$-modules.

Let $\Mod{\cR}^\cU\rightarrow\Schfp{\cR}$
be the full subcategory of the pullback of this fibration along $\Schfp{\cR}\hookrightarrow\SCH$
consisting of $\cU$-sheaves, i.e. sheaves in our chosen universe $\cU$.
--- this is an abelian, $\hat{M}$-small bifibration
where the opposite of each fibre has enough injective objects.

\vspace{\abstand}

For a scheme $X$,
denote the category of quasi-projective (respectively finitely presented) $\cO_X$-modules
by $\QCoh{X}$ (respectively $\Modfp{X}$).
Recall from
\cite[5.2.5]{ega1} that an $\cO_X$-module $\cF$ is called \emph{finitely presented}
if for every $x\in X$, there is an open neighborhood $U\subseteq X$ of $x$ and an exact sequence
$\cO_U^m\rightarrow\cO_U^n\rightarrow\cF|_U\rightarrow 0$ of $\cO_U$-modules
with natural numbers $m$ and $n$.
If $X$ is locally noetherian,
this is equivalent to $\cF$ being a coherent $\cO_X$-module.

Let $\QCOH$ (respectively $\MODFP$) be the full subcategory of $\MOD$
whose fibre over $X/A$ is the opposite of $\QCoh{X}$ (respectively of $\Modfp{X}$).

Pulling back along $\Schfp{\cR}\rightarrow\SCH$ and restricting to $\cU$-sheaves,
we get $\hat{M}$-small fibrations $\QCoh{\cR}^\cU$ and $\Modfp{\cR}$ over $\Schfp{\cR}$
(note that any finitely presented $\cO_X$-module for $X$ in $\cS$ is automatically a $\cU$-sheaf).
We sum up the situation in the following diagram of additive fibrations:
\[
  \xymatrix{
    {\Modfp{\cR}} \ar@{^{(}->}[r] \ar[d] &
    {\QCoh{\cR}^\cU} \ar@{^{(}->}[r] \ar[d] &
    {\Mod{\cR}^\cU} \ar@{^{(}->}[r] \ar[d] &
    {\MOD} \ar[d] \\
    {\Schfp{\cR}} \ar@{=}[r] \ar[d] &
    {\Schfp{\cR}} \ar@{=}[r] \ar[d] &
    {\Schfp{\cR}} \ar@{^{(}->}[r] \ar[d] &
    {\SCH} \ar[d] \\
    {\op{\cR}} \ar@{=}[r] &
    {\op{\cR}} \ar@{=}[r] &
    {\op{\cR}} \ar@{=}[r] &
    {\op{\Rings}} \\
  }
\]
The first three columns in this diagram are $\hat{M}$-small, and we enlarge them to get
an additive fibration $\sModfp{\cR}/\scR$,
an abelian fibration $\sQCoh{\cR}/\scR$ and
an abelian bifibration $\sMod{\cR}/\scR$.

For a *scheme $X$, we denote the \emph{opposite} of the fibre of $\sMod{\cR}$
(respectively $\sModfp{\cR}$, respectively $\sQCoh{\cR}$)
over $X$ by $\sMod{X}$ (respectively $\sModfp{X}$, respectively $\sQCoh{X}$),
and we call the objects of this fibre \emph{$\cO_X$-*modules}
(respectively \emph{*finitely presented} $\cO_X$-*modules,
respectively \emph{*quasi-coherent} $\cO_X$-*modules).

If $X$ is *locally noetherian,
we also say \emph{*coherent} instead of *finitely presented,
and $\sCoh{X}:=\sModfp{X}$ is an abelian category.

\vspace{\abstand}

\begin{lemmadefi}\label{defTcoh}
  Sending $\langle\cF,X/A\rangle$ to $\langle\tm{X}^*\cF,\T{X}/\str{A}\rangle$
  induces a canonical morphism of additive fibrations $\T{}:\Modfp{\cR}\rightarrow\Modfp{\scR}$:
  \[
    \xymatrix{
      {\Modfp{\cR}} \ar@{.>}[r]^{\T{}} \ar[d] & {\Modfp{\scR}} \ar[d] \\
      {\Schfp{\cR}} \ar[r]_{\T{}} \ar[d] & {\Schfp{\scR}} \ar[d] \\
      {\op{\cR}} \ar[r]_{*} & {\op{\scR}} \\
    }
  \]
\end{lemmadefi}

\vspace{\abstand}

\begin{proof}
  This is obvious.
\end{proof}

\vspace{\abstand}

\begin{thm}\label{thmScoh}
  There is an (essentially) unique morphism $\s{}:\Modfp{\scR}\rightarrow\sModfp{\cR}$
  of additive fibrations over $\sSch{\cR}$ that makes the following diagram commute:
  \[
    \xymatrix{
      {\Modfp{\cR}} \ar@(ur,ul)[rr]^{*} \ar[r]_{\T{}} \ar[d] & {\Modfp{\scR}} \ar@{.>}[r]_{\s{}} \ar[d] &
        {\sModfp{\cR}} \ar[d] \\
      {\Schfp{\cR}} \ar[r]_{\T{}} \ar[d] & {\Schfp{\scR}} \ar[r]_{\s{}} \ar[d] & {\sSchfp{\cR}} \ar[d] \\
      {\op{\cR}} \ar[r]_{*} & {\op{\scR}} \ar@{=}[r] & {\op{\scR}} \\
    }
  \]
  In particular, for every *ring $A$ and every finitely presented $A$-scheme $X$,
  we get a canonical additive functor $\s{}:\Modfp{X}\rightarrow\sModfp{\s{X}}$.
\end{thm}

\vspace{\abstand}

\begin{proof}
  This follows from \cite[8.5.2]{ega43} in the same way
  as \ref{thmS} follows from \cite[8.8.2]{ega43}.
\end{proof}

\vspace{\abstand}

\noindent
From now on for the rest of this section,
let $A$ be a *ring, and let $X$ be a finitely presented $A$-scheme.

\vspace{\abstand}

\begin{satz}\label{satzsre}
  Let $\cF\rightarrow\cG\rightarrow\cH\rightarrow 0$ be a sequence in $\Modfp{X}$
  which is exact in $\Mod{X}$.
  Then the sequence $\s{\cF}\rightarrow\s{\cG}\rightarrow\s{\cH}\rightarrow 0$
  of *finitely-presented $\cO_{\s{X}}$-*modules is exact in $\sMod{\s{X}}$.

  In particular, if $A$ is noetherian (for example a *field),
  then the functor $\s$ from coherent $\cO_X$-modules to *coherent $\cO_{\s{X}}$-*modules is right exact.
\end{satz}

\vspace{\abstand}

\begin{proof}
  This follows from \cite[8.5.6]{ega43} and the construction of $\s{}$.
\end{proof}

\vspace{\abstand}

\begin{satz}\label{satzox}
  For $n\in\N_0$, we have $\s{\cO_X^n}=\cO_{\s{X}}^n$.
\end{satz}

\vspace{\abstand}

\begin{proof}
  Since $\s{}$ is additive, we only have to consider the case $n=1$.
  Because $A$ is a *ring, we have a canonical morphism of *rings $\str{\Z}\rightarrow A$
  and hence a canonical morphism $f:X/A\rightarrow\spec{\str{\Z}}/\str{\Z}$ in $\Schfp{\scR}$.
  Then $\cO_X=f^*\cO_{\spec{\str{\Z}}}$, so
  \[
    \s{\cO_X}
    \stackrel{\ref{thmScoh}}{=}(\s{f})^*\s{\cO_{\spec{\str{\Z}}}}
    =(\s{f}^*)\s{\T{\cO_{\spec{\Z}}}}
    =(\s{f})^*\cO_{\str{\spec{\Z}}}
    =\cO_{\s{X}}.
  \]
\end{proof}

\vspace{\abstand}

\begin{cor}\label{corVB}
  Let $\cE$ be a vector bundle of rank $n\in\N_0$ on $X$.
  Then $\s{\cE}$ is a *vector bundle of rank $n$ on $\s{X}$.
\end{cor}

\vspace{\abstand}

\begin{proof}
  This follows immediately from \ref{satzox}.
\end{proof}

\vspace{\abstand}

\begin{lemmadefi}\label{lemmadefSls}
  For an $\cO_{\s{X}}$-*module $\cF$, sending a quasi-compact open subscheme $U$ of $X$ to
  $\cF(\s{U})$ defines an abelian sheaf $\sls{\cF}$ on $X$.
  In this way, we get an additive functor $\sls{}$ from $\sMod{\s{X}}$ to the category of abelian sheaves on $X$.
\end{lemmadefi}

\vspace{\abstand}

\begin{proof}
  First of all, $\sls{\cF}$ is clearly an abelian presheaf on the category of quasi-compact open subsets of $X$,
  because $\s{}$ is a functor from that category to the category of *open *subschemes of $\s{X}$.
  By \ref{bemQC}, such a presheaf defines a sheaf on $X$,
  provided the sheaf-condition with respect to finite, quasi-compact, open coverings is satisfied.

  So let $U\subseteq X$ be quasi-compact and open, and
  let $U=U_1\cup\ldots\cup U_n$ be a finite, quasi-compact, open covering of $U$.
  Then by \ref{corcov}, $[\s{U}]=[\s{U_1}]\cup\ldots\cup[\s{U_n}]$ is a *open covering of $\s{U}$,
  which is \emph{internal} because it is finite.
  By transfer, since $\cF$ is a $\cO_{\s{X}}$-*module, we get the following exact sequence
  (of abelian *groups):
  \[
    0\longrightarrow\cF(\s{U})\longrightarrow\str{\prod_{i=1}^n}\cF(\s{U_i})
    \longrightarrow\str{\prod_{i,j=1}^n}\cF(\s{U_i}\cap\s{U_j}).
  \]
  But $n$ is finite, and finite *products are simply products, so we get the following sequence of abelian groups
  \[
    0\longrightarrow[\sls{\cF}](U)\longrightarrow\prod_{i=1}^n[\sls{\cF}](U_i)
    \longrightarrow\prod_{i,j=1}^n[\sls{\cF}](U_i\cap U_j),
  \]
  which is just the sheaf condition we wanted to prove, so $\sls{\cF}$ is indeed an abelian sheaf on $X$.

  Finally, since $\s{}$ is a functor, we really get an additive functor $\sls{}$ as desired.
\end{proof}

\vspace{\abstand}

\begin{defi}\label{defXhat}
  Since $\sls{\cO_{\s{X}}}$ is a sheaf of rings on $X$ by \ref{lemmadefSls}, we get a ringed space
  \[
    \hat{X}:=(X,\cO_{\hat{X}}):=(X,\sls{\cO_{\s{X}}}),
  \]
  and from now on, we want to consider $\sls{}$ as an additive functor from $\sMod{\s{X}}$
  to $\Mod{\hat{X}}$.

  If $U=\spec{B}$ is an affine, open subscheme of $X$, then we have a canonical morphism of $A$-algebras
  \[
    \cO_X(U)=B\xrightarrow{\sm{B}}\s{B}=\cO_{\s{X}}(\s{U})=\cO_{\hat{X}}(U),
  \]
  which is functorial in $U$ by \eqref{eqSfunct}, i.e. we get a morphism of sheaves of rings
  $\sigma:\cO_X\rightarrow\cO_{\hat{X}}$ on $X$ and hence a canonical morphism of ringed spaces
  $\sigma^*:\hat{X}\rightarrow{X}$,
  which in turn defines a canonical additive functor $\sigma_*:\Mod{\hat{X}}\rightarrow\Mod{X}$.
  We denote the composition
  \[
    \sMod{\s{X}}\xrightarrow{\sls{}}\Mod{\hat{X}}\xrightarrow{\sigma_*}\Mod{X}
  \]
  by $\B{}$.
\end{defi}

\vspace{\abstand}

\begin{satz}\label{satzslsBex}
  The functors $\sls{}:\sMod{\s{X}}\rightarrow\Mod{\hat{X}}$
  and $\B{}:\sMod{\s{X}}\rightarrow\Mod{X}$ are left exact,
  and their restrictions to $\sQCoh{\s{X}}$ are exact and faithful.
\end{satz}

\vspace{\abstand}

\begin{proof}
  The functor
  $\sigma_*:\Mod{\hat{X}}\rightarrow\Mod{X}$ is exact and faithful,
  because it is the identity functor on the underlying abelian sheaves,
  so if $\B{}$ is left exact respectively exact and faithful,
  so is $\sls{}$.

  Let
  $0\rightarrow\cF'\rightarrow\cF\rightarrow\cF''\rightarrow 0$
  be an exact sequence of $\cO_{\s{X}}$-*modules.
  If $U$ is a *open *subscheme of $\s{X}$, then by transfer the sequence
  \[
    0\longrightarrow\cF'(U)\longrightarrow\cF(U)\longrightarrow\cF''(U)
  \]
  is exact (in the category of internal $\cO_{\s{X}}(U)$-modules and hence in particular in the category of
  abelian groups),
  which proves that $\B{}$ is left exact.

  Now let $\cF'$, $\cF$ and $\cF''$ be *quasi-coherent.
  Let $x\in X$ be an arbitrary point, and let $t_x\in[\B{\cF''}]_x$ be an arbitrary element in the stalk.
  There is an affine open subscheme $U$ of $X$ with a local section $t_U\in[\B{\cF''}](U)$
  which represents $t_x$.
  Since $\cF'$, $\cF$ and $\cF'$ are *quasi-coherent and since $\s{U}$ is *affine,
  it follows by transfer that
  \[
    0\longrightarrow\cF'(\s{U})\longrightarrow\cF(\s{U})\longrightarrow\cF''(\s{U})\longrightarrow 0
  \]
  is exact, so that there is a preimage
  $s_U\in\cF(\s{U})=[\B{\cF}](U)$ of $t_U$ which then represents a preimage $s_x\in[\B{\cF}]_x$ of $t_x$.
  This shows that $\B{}$ is also right exact and hence exact.

  Now let $\cF\xrightarrow{\varphi}\cG$ be a morphism of *quasi-coherent $\cO_{\s{X}}$-*modules
  with $\B{\varphi}=0$. For faithfulness, we have to show $\varphi=0$.
  Choose a finite affine open covering $X=U_1\cup\ldots\cup U_n$ of $X$.
  Then $[\s{X}]=[\s{U_1}]\cup\ldots\cup[\s{U_n}]$ is a finite *affine *open covering by \ref{corcov},
  and it suffices to show $\varphi\vert_{\s{U_i}}=0$ for all $i\in\{1,\ldots,n\}$
  or equivalently --- because $\cF$ and $\cG$ are *quasi-coherent --- $\varphi_{\s{U_i}}=0$ for all $i$.
  But $\varphi_{\s{U_i}}=[\B{\varphi}]_{U_i}=0$, and we are done.
\end{proof}

\vspace{\abstand}

Let $\cF$ and $\cG$ be $\cO_{\s{X}}$-*modules. Then $\s{}$ induces a canonical morphism
\begin{equation}\label{eqMorHom}
  \sls{\HOM{\s{X}}{\cF}{\cG}}\rightarrow\HOM{\hat{X}}{\sls{\cF}}{\sls{\cG}}
\end{equation}
of $\cO_{\hat{X}}$ modules
by
\begin{multline*}
  \bigl[\sls{\HOM{\s{X}}{\cF}{\cG}}\bigr](U)
  =\Hom{\cO_{\s{U}}}{\cF\vert_{\s{U}}}{\cG\vert_{\s{U}}} \\
  \xrightarrow{\sls{}}
  \Hom{\cO_{\hat{U}}}{\sls\cF\vert_U}{\sls\cG\vert_U}
  =\bigl[\HOM{\hat{X}}{\sls{\cF}}{\sls{\cG}}\bigr](U)
\end{multline*}
for quasi-compact, open subschemes $U$ of $X$.

\vspace{\abstand}

\begin{satz}\label{satzslsff}
  Let $\cF$ be a finitely presented $\cO_X$-module,
  and let $\cG$ be an $\cO_{\s{X}}$-*module.
  Then the canonical morphism \eqref{eqMorHom} (for $\s{\cF}$ and $\cG$)
  \[
    \sls{\HOM{\s{X}}{\s{\cF}}{\cG}}
    \longrightarrow
    \HOM{\hat{X}}{\sls{\s{\cF}}}{\sls{G}}
  \]
  of $\cO_{\hat{X}}$-modules is an isomorphism.
  Taking global sections, this in particular implies that
  \[
    \Hom{\cO_{\s{X}}}{\s{\cF}}{\cG}
    \xrightarrow{\sls{}}
    \Hom{\cO_{\hat{X}}}{\sls{\s{\cF}}}{\sls{\cG}}
  \]
  is an isomorphism.
\end{satz}

\vspace{\abstand}

\begin{proof}
  The question whether a given morphism of sheaves on $X$ is an isomorphism is local on $X$,
  so we can assume that $X$ is affine.
  If $\cF=\cO_X^n$, then $\s{\cF}=\cO_{\s{X}}^n$ and $\sls{\s{\cF}}=\cO_{\hat{X}}^n$,
  i.e. $\HOM{\s{X}}{\s{\cF}}{\cG}$ is canonically isomorphic
  to $\cG^n$ (by transfer),
  and $\HOM{\hat{X}}{\sls{\s{\cF}}}{\sls{G}}$ is canonically isomorphic to $\sls{\cG}^n$,
  so that the statement is obviously true in this case.

  In the general case --- since $X$ is affine --- there is a finite presentation
  \[
    \cO_X^m\longrightarrow\cO_X^n\longrightarrow\cF\longrightarrow 0
  \]
  of $\cF$,
  which (by \ref{satzsre} and \ref{satzox}) induces an exact sequence
  \[
    \cO_{\s{X}}^m\longrightarrow\cO_{\s{X}}^n\longrightarrow\s{\cF}\longrightarrow 0
  \]
  of $\cO_{\s{X}}$-*modules
  and (by \ref{satzslsBex}) an exact sequence
  \[
    \cO_{\hat{X}}^m\longrightarrow\cO_{\hat{X}}^n\longrightarrow\sls{\s{\cF}}\longrightarrow 0
  \]
  of $\cO_{\hat{X}}$-modules.
  Since the functors
  \[
    \begin{array}{lccccl}
      \HOM{\s{X}}{\_}{\cG} & : & \sMod{\s{X}} & \rightarrow & \sMod{\s{X}}, & \\[2mm]
      \HOM{\hat{X}}{\_}{\sls{\cG}} & : & \Mod{\hat{X}} & \rightarrow & \Mod{\hat{X}} & \text{and} \\[2mm]
      \sls{} & : & \sMod{\s{X}} & \rightarrow & \Mod{\hat{X}} & \\
    \end{array}
  \]
  are left exact,
  we get the following commutative diagram of $\cO_{\hat{X}}$-modules with exact rows:
  \[
    \xymatrix{
      0 \ar[r] &
      {\sls{\HOM{\s{X}}{\s{\cF}}{\cG}}} \ar[r] \ar@{.>}[d]_\alpha &
      {\sls{\HOM{\s{X}}{\cO_{\s{X}}^n}{\cG}}} \ar[r] \ar[d]_\beta &
      {\sls{\HOM{\s{X}}{\cO_{\s{X}}^m}{\cG}}} \ar[d]_\gamma \\
      0 \ar[r] &
      {\HOM{\hat{X}}{\sls{\s{\cF}}}{\sls{\cG}}} \ar[r] &
      {\HOM{\hat{X}}{\cO_{\hat{X}}^n}{\sls{\cG}}} \ar[r] &
      {\HOM{\hat{X}}{\cO_{\hat{X}}^m}{\sls{\cG}}.} \\
    }
  \]
  According to the first case, $\beta$ and $\gamma$ are isomorphisms.
  But then $\alpha$ must be an isomorphism as well, and we are done.
\end{proof}

\vspace{\abstand}

Let $\cF$ be a finitely $\cO_X$-module. Choose a subring $A_0$ of $A$ of finite type over $\Z$,
a scheme $X_0$ of finite type over $A_0$ and a finitely presented $\cO_{X_0}$-module $\cF_0$
such that $\langle\cF,X/A\rangle$ is the pullback of $\langle\cF_0,X_0/A_0\rangle$ along
$\varphi:=A_0\hookrightarrow A$.

By \ref{thmS} and \ref{thmScoh}, we get the following diagram (where we put
$\bar{\varphi}:=\bsi{\Z}{A_0}{A}[\varphi]$):
\[
  \xymatrix@C=10mm@R=10mm{
    & {\cF} \ar[dl] \ar@{-}[dr] & & \ar@{.}[ddd] & & {\s{\cF}} \ar@{-}[dl] \ar[dr] \\
    {\cF_0} \ar@{-}[dr] & \Box & X \ar[dl]^f \ar[dr] & &
      {\s{X}} \ar[dl] \ar[dr]_{\bar{f}} & \Box & {\str{\cF_0}} \ar@{-}[dl] \\
    & {X_0} \ar[dr] & \Box & A \ar[dl]^{\op{\varphi}} \ar[dr]_{\op{\bar{\varphi}}} & \Box & {\str{X_0}} \ar[dl] \\
    & & {A_0} \ar[rr]_{*} & & {\str{A_0}.} \\
  }
\]
The squares are cartesian (in $\SCHFP$ respectively $\MODFP$ on the left,
in $\sSchfp{\scR}$ respectively $\sModfp{\scR}$ on the right),
and we have isomorphisms $f^*\cF_0\xrightarrow{\sim}\cF$ and $\bar{f}^*(\str{\cF_0})\xrightarrow{\sim}\s{\cF}$
and their adjoints $\cF_0\rightarrow f_*\cF$ and $\str{\cF_0}\rightarrow\bar{f}_*\s{\cF}$.

Now let $U_0$ be an open subscheme of $X_0$, and put $U:=U_0\times_{X_0}X$.
We get an $\cO_{X_0}(U_0)$-linear map
\[
  F_0(U_0)\xrightarrow{*}
  [\str{\cF_0}](\str{U_0})\longrightarrow
  [\bar{f}_*\s{\cF}](\str{U_0})
  =[\s{\cF}](\s{U})
  =[\B{\s{\cF}}](U)
  =[f_*\B{\s{\cF}}](U_0)
\]
which is clearly functorial in $U_0$, so that we get a morphism of $\cO_{X_0}$-modules
$\cF_0\longrightarrow f_*\B{\s{\cF}}$ and hence --- by adjunction --- a canonical morphism of
$\cO_X$-modules $\cF\longrightarrow\B{\s{\cF}}$.

This morphism is clearly functorial in $\cF$, so that we get a canonical morphism of functors
\begin{equation}\label{eqIdSN}
  \xymatrix{
    {\Modfp{X}} \ar@{}[r]|{\displaystyle\Downarrow} \ar@(ur,ul)[r]^{\mathrm{can}} \ar@(dr,dl)[r]_{\B{\s{}}} &
    {\Mod{X}}
  }
\end{equation}
and --- again taking adjoints --- a canonical morphism of functors
\begin{equation}\label{eqTNlsN}
  \xymatrix{
    {\Modfp{X}} \ar@{}[r]|{\displaystyle\Downarrow} \ar@(ur,ul)[r]^{\_\otimes_{\cO_X}\cO_{\hat{X}}}
      \ar@(dr,dl)[r]_{\sls{\s{}}} &
    {\Mod{\hat{X}}.}
  }
\end{equation}

\vspace{\abstand}

\begin{satz}\label{satzTNlsNiso}
  The canonical morphism of functors \eqref{eqTNlsN} is an isomorphism.
\end{satz}

\vspace{\abstand}

\begin{proof}
  Let $\cF$ be a finitely presented $\cO_X$-module. We claim that the canonical morphism
  $\cF\otimes_{\cO_X}\cO_{\hat{X}}\longrightarrow\sls{\s{\cF}}$ of $\cO_{\hat{X}}$-modules
  (or of abelian sheaves on $X$) is an isomorphism. This claim is local in $X$,
  so we can assume that $X$ is affine and choose a finite presentation
  \[
    \cO_X^m\longrightarrow\cO_X^n\longrightarrow\cF\longrightarrow 0.
  \]
  By \ref{satzsre}, \ref{satzox} and \ref{satzslsBex},
  we get the following commutative diagram of $\cO_{\hat{X}}$-modules with exact rows:
  \[
    \xymatrix{
      {\cO_{\hat{X}}^m} \ar[r] \ar[d]^{\wr} &
      {\cO_{\hat{X}}^n} \ar[r] \ar[d]^{\wr} &
      {\cF\otimes_{\cO_X}\cO_{\hat{X}}} \ar[r] \ar@{.>}[d] &
      0 \\
      {\cO_{\hat{X}}^m} \ar[r] &
      {\cO_{\hat{X}}^n} \ar[r] &
      {\sls{\s{\cF}}} \ar[r] &
      {0.}
    }
  \]
  The first two vertical morphisms are obviously simply the identity, so the third vertical morphism must be
  an isomorphism.
\end{proof}

\vspace{\abstand}

\begin{satz}\label{satzaff}
  For any affine open subscheme $U=\spec{B}$ of $X$,
  there is a canonical isomorphism of functors
  \begin{equation}\label{eqaff}
    \xymatrix@C=40mm{
      {\Modfp{X}}
        \ar@(ur,ul)[r]^{\Gamma_U(\_)\otimes_B\s{B}}
        \ar@(dr,dl)[r]_{\Gamma_{\s{U}}\circ\s{}}
        \ar@{}[r]|{\displaystyle\Downarrow\wr} &
      {\MMod{[\s{B}]}.} \\
    }
  \end{equation}
\end{satz}

\vspace{\abstand}

\begin{proof}
  Using \eqref{eqIdSN}, composed with $\Gamma_U$, defines a canonical morphism of functors
  \[
    \xymatrix{
      {\Modfp{X}} \ar[rr]^{\mathrm{can}} \ar[rd]_{\s{}} &
      \ar@{}[d]|{{\displaystyle\Downarrow}\eqref{eqIdSN}} &
      {\Mod{X}} \ar[rr]^{\Gamma_U} \ar@{}[rd]|{\displaystyle =} & &
      {\MMod{B}} \\
      & {\sModfp{\s{X}}} \ar[ru]_{\B{}} \ar[rr]_{\Gamma_{\s{U}}} & &
      {\MMod{[\s{B}}]} \ar[ru]_{\mathrm{can}} & \\
    }
  \]
  and thus by adjunction the morphism of functors \eqref{eqaff}.
  To see that this is an isomorphism, let $\cF$ be a finitely presented $\cO_X$-module,
  and choose a finite presentation
  \[
    B^m\longrightarrow
    B^n\longrightarrow
    \cF(U)\longrightarrow
    0.
  \]
  Taking associated sheaves and applying $\s{}$, we get an exact sequence of *finitely presented
  $\cO_{\s{U}}$-*modules
  \[
    \cO_{\s{U}}^m\longrightarrow
    \cO_{\s{U}}^n\longrightarrow
    [\s{\cF}]\vert_{\s{U}}\longrightarrow
    0.
  \]
  By transfer, $\Gamma_{\s{U}}:\Modfp{\s{U}}\longrightarrow\MMod{[\s{B}]}$ is exact,
  so we get the exact sequence of $\s{B}$-modules
  \[
    \s{B}^m\longrightarrow
    \s{B}^n\longrightarrow
    [\s{\cF}](\s{U})\longrightarrow
    0
  \]
  and hence the following commutative diagram of $\s{B}$-modules  with exact rows:
  \[
    \xymatrix{
      {B^m\otimes_B\s{B}} \ar[r] \ar[d]^{\alpha} &
      {B^n\otimes_B\s{B}} \ar[r] \ar[d]^{\beta} &
      {\cF(U)\otimes_B\s{B}} \ar[r] \ar@{.>}[d]^{\gamma} &
      0 \\
      {\s{B}^m} \ar[r] &
      {\s{B}^n} \ar[r] &
      {[\s{\cF}](\s{U})} \ar[r] &
      {0.} \\
    }
  \]
  Since $\alpha$ and $\beta$ are clearly isomorphisms, so is $\gamma$, and we are done.
\end{proof}

\vspace{\abstand}

\begin{cor}\label{corbsmod}
  The canonical functors
  \[
    \begin{array}{ccccc}
      \op{\bigl(\Modfp{X}\bigr)} & \times & \sMod{\s{X}} & \longrightarrow & \Sets \\[3mm]
      (\cF & , & \cG) & \mapsto &
      \left\{\begin{array}{l}
        \Hom{\cO_{\s{X}}}{\s{\cF}}{\cG} \\[3mm]
        \Hom{\cO_X}{\cF}{\B{\cG}}
      \end{array}\right.
    \end{array}
  \]
  are canonically isomorphic via
  \[
    \bs{X}{\cF}{\cG}:\Hom{\cO_{\s{X}}}{\s{\cF}}{\cG}
    \longrightarrow
    \Hom{\cO_X}{\cF}{\B{\cG}},\;\;
    \bigl[\s{\cF}\xrightarrow{\varphi}\cG\bigr]\mapsto
    \bigl[\cF\xrightarrow{\eqref{eqIdSN}}\B{\s{\cF}}\xrightarrow{\B{\varphi}}\B{\cG}\bigr].
  \]
\end{cor}

\vspace{\abstand}

\begin{proof}
  Let $\cF$ be a finitely presented $\cO_X$-module, and let $\cG$ be an $\cO_{\s{X}}$-*module.
  Then
  \[
    \Hom{\cO_{\s{X}}}{\s{\cF}}{\cG}
    \stackrel{\ref{satzslsff}}{\cong}
    \Hom{\cO_{\hat{X}}}{\sls{\s{\cF}}}{\sls{\cG}}
    \stackrel{\ref{satzTNlsNiso}}{\cong}
    \Hom{\cO_{\hat{X}}}{\cF\otimes_{\cO_X}\cO_{\hat{X}}}{\sls{\cG}}
    \stackrel{\mathrm{adj.}}{\cong}
    \Hom{\cO_X}{\cF}{\B{\cG}},
  \]
  and it is clear that the composition of these canonical isomorphisms is just $\bs{X}{\cF}{\cG}$.
\end{proof}

\vspace{\abstand}

\begin{bem}\label{bemtaunat}
  Note that the functoriality of the isomorphism from \ref{corbsmod} in particular implies the following:
  If $\cF\xrightarrow{\varphi}\cF'$ is a morphism of finitely presented $\cO_X$-modules,
  if $\cG\xrightarrow{\psi}\cG'$ is a morphism of $\cO_{\s{X}}$-*modules,
  and if
  \begin{equation}\label{eqtaunat}
    \xymatrix{
      {\s{\cF}} \ar[r]^{\s{\varphi}} \ar[d]_{f} & {\s{\cF'}} \ar[d]^{g} \\
      {\cG} \ar[r]_{\psi} & {\cG'} \\
    }
  \end{equation}
  is a diagram of $\cO_{\s{X}}$-*modules,
  then \eqref{eqtaunat} commutes if and only if the corresponding diagram
  \begin{equation}
    \xymatrix{
      {\cF} \ar[r]^{\varphi} \ar[d]_{\bs{X}{\cF}{\cG}(f)} & {\cF'} \ar[d]^{\bs{X}{\cF'}{\cG'}(g)} \\
      {\B{\cG}} \ar[r]_{\B{\psi}} & {\B{\cG'}} \\
    }
  \end{equation}
  of $\cO_X$-modules commutes.
\end{bem}

\vspace{\abstand}

\begin{cor}\label{cortensor}
  Let $\cF$ and $\cG$ be two finitely presented $\cO_X$-modules.
  There is a canonical isomorphism of *finitely presented $\cO_{\s{X}}$-modules
  \[
    \s{\bigl(\cF\otimes_{\cO_X}\cG\bigr)}
    \stackrel{\sim}{\longrightarrow}
    \s{\cF}\otimes_{\cO_{\s{X}}}\s{\cF}.
  \]
\end{cor}

\vspace{\abstand}

\begin{proof}
  For a quasi-compact open subscheme $U$ of $X$, we have a canonical $\cO_X(U)$-linear map
  \begin{multline*}
    \cF(U)\otimes_{\cO_X(U)}\cG(U)
    \xrightarrow{\eqref{eqIdSN}}
    [\B{\s{\cF}}](U)\otimes_{[\B{\s{\cO_X}}](U)}[\B{\s{\cG}}](U) \\
    =
    [\s{\cF}](\s{U})\otimes_{\cO_{\s{X}}(\s{U})}[\s{\cG}](\s{U})
    \xrightarrow{\mathrm{can}}
    \bigl[\s{\cF}\otimes_{\cO_{\s{X}}}\s{\cG}\bigr](\s{U})
    =\B{\bigl[\s{\cF}\otimes_{\cO_{\s{X}}}\s{\cG}\bigr]}(U),
  \end{multline*}
  which is clearly functorial in $U$ and consequently defines
  a functorial morphism of \emph{presheaves} of $\cO_X$-modules
  \[
    \Bigl[U\mapsto\cF(U)\otimes_{\cO_X(U)}\cG(U)\Bigr]
    \longrightarrow
    \B{\bigl[\s{\cF}\otimes_{\cO_{\s{X}}}\s{\cG}\bigr]}
  \]
  and then, by the universal property of the associated sheaf,
  a functorial morphism of $\cO_X$-modules
  \[
    \cF\otimes_{\cO_X}\cG
    \longrightarrow
    \B{\bigl[\s{\cF}\otimes_{\cO_{\s{X}}}\s{\cG}\bigr]},
  \]
  which by \ref{corbsmod} and \ref{bemtaunat} corresponds to a functorial morphism of $\cO_{\s{X}}$-*modules
  \begin{equation}\label{eqtensor}
    \s{\bigl[\cF\otimes_{\cO_X}\cG\bigr]}
    \longrightarrow
    \s{\cF}\otimes_{\cO_{\s{X}}}\s{\cG}.
  \end{equation}
  To prove that \eqref{eqtensor} is an isomorphism, choose a quadruple $\langle A_0,X_0,\cF_0,\cG_0\rangle$,
  where $A_0\stackrel{\varphi}{\hookrightarrow}A$ is a finitely generated subring of $A$,
  $X_0$ is an $A_0$-scheme of finite type with $X\cong\varphi^*X_0$
  and $\cF_0$ and $\cG_0$ are coherent sheaves on $X_0$ with
  $\cF\cong\varphi^*\cF_0$ and $\cG\cong\varphi^*\cG_0$.
  Then of course we also have $\varphi^*[\cF_0\otimes_{\cO_{X_0}}\cG_0]\cong\cF\otimes_{\cO_X}\cG$
  and therefore (with $\bar{\varphi}:=\bsi{\Z}{A_0}{A}[\varphi]$)
  \begin{multline*}
    \s[\cF\otimes_{\cO_X}\cG]
    \cong
    \bar{\varphi}^*\Bigl(\str{\bigl[\cF_0\otimes_{\cO_{X_0}}\cG_0\bigr]}\Bigr)
    =
    \bar{\varphi}^*\Bigl(\str{\cF_0}\otimes_{\cO_{\str{X_0}}}\str{\cG_0}\Bigr) \\
    \cong
    \bar{\varphi}^*\bigl(\str{\cF_0}\bigr)\otimes_{\cO_{\s{X}}}\bar{\varphi}^*\bigl(\str{\cG_0}\bigr)
    \cong\s{\cF}\otimes_{\cO_{\s{X}}}\s{\cG}
  \end{multline*}
  by construction of $\s{}$.
\end{proof}

\vspace{\abstand}

\begin{cor}\label{corpic}
  The functor $\s{}:\Modfp{X}\longrightarrow\sModfp{\s{X}}$ induces a canonical group homomorphism
  $\s{}:\Pic{X}\longrightarrow\sPic{\s{X}}$ between the Picard group of $X$
  and the *Picard group of $\s{X}$.
\end{cor}

\vspace{\abstand}

\begin{proof}
  By \ref{corVB}, $\s{}$ sends line bundles to line bundles,
  so we get a map $\s{}:\Pic{X}\longrightarrow\sPic{\s{X}}$.
  This map is a group homomorphism by \ref{cortensor}.
\end{proof}

\vspace{\abstand}

\begin{cor}\label{corhom}
  and let $\cF$ and $\cG$ be two finitely presented $\cO_X$-modules
  with the property that the $\cO_X$-module $\HOM{X}{\cF}{\cG}$ is also finitely-presented,
  which is for example the case if
  \begin{itemize}
    \item
      $\cF$ is a vector bundle or
    \item
      $\cF$ and $\cG$ are coherent.
  \end{itemize}
  Then there is a canonical morphism of *finitely presented $\cO_{\s{X}}$-*modules
  \begin{equation}\label{eqhom}
    \s{\HOM{X}{\cF}{\cG}}
    \longrightarrow
    \HOM{\s{X}}{\s{\cF}}{\s{\cG}}
  \end{equation}
  which is an isomorphism if $\cF$ is a vector bundle.
\end{cor}

\vspace{\abstand}

\begin{proof}
  Look at the following canonical map of sets of morphisms:
  \begin{multline*}
    \Mor{\Modfp{X}}{\HOM{X}{\cF}{\cG}}{\HOM{X}{\cF}{\cG}}
    \cong\Mor{\Modfp{X}}{\HOM{X}{\cF}{\cG}\otimes_{\cO_X}\cF}{\cG} \\
    \xrightarrow{\s{}}\Mor{\sModfp{\s{X}}}{\s{[\HOM{X}{\cF}{\cG}\otimes_{\cO_X}\cF]}}{\s{\cG}} \\
    \stackrel{\ref{cortensor}}{\cong}
      \Mor{\sModfp{\s{X}}}{\s{\HOM{X}{\cF}{\cG}}\otimes_{\cO_{\s{X}}}\s{\cF}}{\s{\cG}} \\
    \cong\Mor{\sModfp{\s{X}}}{\s{\HOM{X}{\cF}{\cG}}}{\HOM{\s{X}}{\s{\cF}}{\s{\cG}}},
  \end{multline*}
  and take the identity's image under this map to get \eqref{eqhom}.

  Now let $\cF$ be a vector bundle.
  Since the question whether \eqref{eqhom} is an isomorphism is local,
  we can assume that $\cF=\cO_X^n$ is trivial,
  and we have
  \[
    \s{\HOM{X}{\cF}{\cG}}
    \cong
    \s{\cG^n}
    \stackrel{\mathrm{transfer}}{\cong}
    \Hom{\s{X}}{\cO_{\s{X}}^n}{\s{\cG}}
    \stackrel{\ref{satzox}}{\cong}
    \Hom{\s{X}}{\s{\cO_X^n}}{\s{\cG}}
  \]
  as desired.
\end{proof}

\vspace{\abstand}

\begin{cor}\label{cordual}
  For a vector bundle $\cE$ on $X$,
  there is a canonical isomorphism
  $\s{(\cE^\vee)}\cong(\s{\cE})^\vee$.
\end{cor}

\vspace{\abstand}

\begin{proof}
  This follows immediately from \ref{satzox} and \ref{corhom}:
  \[
    \s{(\cE^\vee)}
    =\s\HOM{X}{\cE}{\cO_X}
    \stackrel{\ref{corhom}}{\cong}\HOM{\s{X}}{\s{\cE}}{\s{\cO_X}}
    \stackrel{\ref{satzox}}{=}\HOM{\s{X}}{\s{\cE}}{\cO_{\s{X}}}
    =(\s{\cE})^\vee.
  \]
\end{proof}

\vspace{\abstand}

\begin{satz}\label{satzO1}
  Let $n\in\N_+$ be a natural number,
  and let $k\in\Z$ be an integer.
  \begin{enumerate}
    \item\label{O1proj}
      Under the functor $\s{}$, the invertible $\cO_{\P^n_X}$-module $\cO_{\P^n_X}(k)$ is mapped
      to the *invertible $\cO_{\str{\P^n_{\s{X}}}}$-*module $\cO_{\str{\P^n_{\s{X}}}}(k)$.
    \item\label{O1twist}
      Let $i:Y\hookrightarrow\P^n_X$ be a closed immersion of finitely presented $A$-schemes,
      and let $\cF$ be a finitely presented $\cO_Y$-module.
      Then $\s{[\cF(k)]}=[\s{\cF}](k)$,
      where the twists are taken with respect to $i$ respectively $\s{i}$.
  \end{enumerate}
\end{satz}

\vspace{\abstand}

\begin{proof}
  Choose a finitely generated subring $A_0\stackrel{\varphi}{\hookrightarrow}A$ of $A$
  and an $A_0$-scheme $X_0$ of finite type with $\varphi^*X_0=X$.
  Then $\varphi^*\P^n_{X_0}=\P^n_X$, so
  \[
    \s{(\cO_{\P^n_A}(k))}
    =\Bigl[\bsi{\Z}{A_0}{A}[\varphi]\Bigr]^*\cO_{\str{\P^n_{X_0}}}(k)
    =\cO_{\s{(\P^n_X)}}(k)
    \stackrel{\ref{satzProj}}{=}\cO_{\str{\P^n_{\s{X}}}}(k),
  \]
  and this is *invertible by \ref{corVB} (or by transfer), so we have \ref{O1proj}.

  For \ref{O1twist} we get:
  \begin{multline*}
    \s{[\cF(k)]}
    =\s{\bigl[\cF\otimes_{\cO_Y}i^*\cO_{\P^n_X}(k)\bigr]}
    \stackrel{\ref{cortensor}}{=}[\s{\cF}]\str{\otimes}_{\cO_{\s{Y}}}\Bigl[\s{\bigl[i^*\cO_{\P^n_X}(k)\bigr]}\Bigr] \\
    =[\s{\cF}]\str{\otimes}_{\cO_{\s{Y}}}\Bigl[(\s{i})^*\bigl[\s{\cO_{\P^n_X}(k)}\bigr]\Bigr]
    \stackrel{\ref{O1proj}}{=}
      [\s{\cF}]\str{\otimes}_{\cO_{\s{Y}}}\Bigl[(\s{i})^*\bigl[\cO_{\str{\P^n_{\s{X}}}}(k)\bigr]\Bigr]
    =[\s{\cF}](k).
  \end{multline*}
\end{proof}

\vspace{\abstand}

If $Z$ is a finitely presented closed subscheme of $X$, given by a finitely presented sheaf of ideals $\cI$ on $X$,
then we know from \ref{satzP} that $\s{Z}$ is a *closed *subscheme of $\s{X}$.
As final result in this section, we want to determine the relationship between $\s{\cI}$ and
the *ideal on $\s{X}$ defining $\s{Z}$:

\begin{satz}\label{satzNNcompat}
  Let $Z$ be a finitely presented closed subscheme of $X$, given by a finitely presented sheaf of ideals $\cI$.
  Then the *closed *subscheme $\s{Z}$ of $\s{X}$ is given by the *ideal
  $\image{\s{\cI}}{\cO_{\s{X}}}$.
\end{satz}

\vspace{\abstand}

\begin{proof}
  Let $\s{Z}$ be given by the *ideal $\cJ$ on $\s{X}$.
  If $U\subseteq X$ is a quasi-compact open subscheme of $X$,
  then the *closed *subscheme $\s{[Z\cap U]}$ of $\s{U}$ is given by $\cJ\vert_{\s{U}}$,
  so we can assume without loss of generality that $X$ is affine, say $X=\spec{B}$ for a finitely presented
  $A$-algebra $B$.
  Then $Z=\spec{B/\b}$ for a finitely presented ideal $\b$ of $B$,
  and $\cI=\tilde{\b}$.

  Using \ref{lemmaSS}, we have
  \[
    \s{[B/\b]}=[\s{B}]/\str{\b}=[\s{B}]\,/\,\b\cdot[\s{B}],
  \]
  so that $\s{Z}$ is given by the *ideal $\cJ=\widetilde{\b\cdot[\s{B}]}$ of $\s{X}$,
  and we have to prove that this *ideal equals $\image{\s{\cI}}{\cO_{\s{X}}}$
  or --- equivalently --- that the global sections of these two *ideals agree (as ideals of $\s{B}$).
  Using \ref{satzaff}, this is easy:
  \[
    \Gamma_{\s{X}}\Bigl[\image{\s{\cI}}{\cO_{\s{X}}}\Bigr]
    =\image{\Gamma_{\s{X}}[\s{\cI}]}{\s{B}}
    \stackrel{\eqref{eqaff}}{=}\image{\b\otimes_B[\s{B}]}{\s{B}}
    =\b\cdot\s{B}
    =\Gamma_{\s{X}}[\cJ].
  \]
\end{proof}

\vspace{\abstand}


\section{The case of varieties}

\vspace{\abstand}

Let $k$ be a *field in $\scR$, i.e. a *ring which is an (internal) field.
Then $k$ is of course a noetherian ring, so that a $k$-scheme $X$ is finitely presented
if and only if it is of finite type,
and an $\cO_X$-module $\cF$ is finitely presented if and only if it is coherent.

\vspace{\abstand}

\begin{defi}\label{defdim}
  We can consider "dimension" as a function
  $
    \dim:
    \{\text{schemes}\}
    \longrightarrow
    \{-\infty\}\amalg\Nn\amalg\{\infty\},
  $
  so by restriction to $\Ob{\cS}$ and enlarging we get an induced function
  \[
    \str{\dim}:
    \{\text{*schemes}\}
    \longrightarrow
    \{-\infty\}\amalg\Nns\amalg\{\infty\}.
  \]
  For a *scheme $X$, we call $\str{\dim X}$ the \emph{*dimension of $X$}.
\end{defi}

\vspace{\abstand}

For the proof of theorem \ref{thmvariety} below, we will need the following results
of van den Dries and Schmidt which we state here --- in our notation --- for the convenience of the reader:

\begin{thm}\textbf{(Lou van den Dries, K. Schmidt)}\label{thmvdDries}\\
  Let $\cI\subseteq k[X_1,\ldots,X_n]$ be an ideal. Then
  \begin{enumerate}
    \item\label{vdDriesff}
      The ring homomorphism $k[X_i]\rightarrow k\str{[X_i]}$ is faithfully flat.
    \item\label{vdDriesprime}
      $\cI$ is prime if and only if $\str{\cI}\subseteq k\str{[X_i]}$ is *prime or --- what amounts to the same,
      since for an ideal being prime is clearly a first order property --- prime.
    \item\label{vdDriesminprime}
      If $\p_1,\ldots,\p_m$ are the distinct minimal primes of $\cI$,
      then $\str{\p_1},\ldots,\str{\p_m}$ are the distinct minimal primes of $\str{\cI}$
      (in particular, all minimal primes of $\str{\cI}$ are *ideals, hence the notions of
      "minimal prime ideal of $\str{\cI}$"
      and "minimal prime *ideal of $\str{\cI}$" coincide).
    \item\label{vdDriesrad}
      $\sqrt[*]{\str{\cI}}=\str{\bigl[\sqrt{\cI}\bigr]}$.
  \end{enumerate}
\end{thm}

\vspace{\abstand}

\begin{proof}
  Part \ref{vdDriesff} is \cite[1.8]{vandendries},
  part \ref{vdDriesprime} is \cite[2.5]{vandendries},
  and parts \ref{vdDriesminprime} and \ref{vdDriesrad} are \cite[2.7]{vandendries};
  that van den Dries and Schmidt's formulation agrees with the one given here follows immediately from \ref{bemideal}.
\end{proof}

\vspace{\abstand}

\begin{cor}\label{corff}
  Let $A$ be a $k$-algebra of finite type.
  Then $\sm{A}:A\rightarrow\s{A}$ is faithfully flat.
\end{cor}

\vspace{\abstand}

\begin{proof}
  Let $A=k[X_1,\ldots,X_N]/\cI$. Then
  \[
    \sm{A}=
    \sm{k[X_i]}\otimes_{k[X_i]}k[X_i]/\cI
    \;\;:\;\;
    A=k[X_i]/\cI\longrightarrow
    k\str{[X_i]}/\cI\cdot k\str{[X_i]}
    \stackrel{\ref{bemideal}}{=}k\str{[X_i]}/\str{\cI}=\s{A},
  \]
  and $\sm{k[X_i]}$ is faithfully flat by \ref{thmvdDries}\ref{vdDriesff},
  so $\sm{A}$ --- as a base change of $\sm{k[X_i]}$ --- must be faithfully flat as well.
\end{proof}

\vspace{\abstand}

\begin{thm}\label{thmvariety}
  Let $X$ be a scheme of finite type over $k$.
  \begin{enumerate}
    \item\label{varietyzero}
      $X$ is the empty scheme if and only if $\s{X}$ is the *empty scheme.
    \item\label{varietydim}
      $\str{\dim\s{X}}=\dim X$.
    \item\label{varietyinteger}
      $X$ is reduced (irreducible, integer) if and only if $\s{X}$ is *reduced (*irreducible, *integer).
    \item\label{varietycohex}
      The functor $\s{}$ from coherent $\cO_X$-modules to (the abelian category of) *coherent $\cO_{\s{X}}$-modules
      is faithful and exact.
  \end{enumerate}
\end{thm}

\vspace{\abstand}

\begin{proof}
  If $X=\emptyset$, then $\s{X}=\str{\emptyset}$ by \ref{satzfin}\ref{finempty},
  so let $\s{X}=\str{\emptyset}$.
  Let us assume that $X\neq\emptyset$.
  Then $X$ contains a $K$-valued point for a finite field extension $K/k$,
  and applying $\s{}$ gives us an $\s{K}$-valued point of $\s{X}$.
  If $K=k[X_1,\ldots,X_n]/\cI$ is a finite presentation of $K$,
  it follows from \ref{thmvdDries}\ref{vdDriesff} that $\s{K}=k\str{[X_i]}/\str{\cI}$
  is not zero, so the existence of an $\s{K}$-valued point of $\s{X}$ proves the existence of
  a *topological point of $\s{X}$, a contradiction to $\s{X}=\str{\emptyset}$.

  Having settled \ref{varietyzero}, for \ref{varietydim} and \ref{varietyinteger}
  we can assume that $X\neq\emptyset$.
  For \ref{varietydim}, we use \cite[4.1.2]{ega42},
  according to which $\dim X=n$ is equivalent to the existence of a diagram
  \[
    \xymatrix{
      {U} \ar@{^{(}->}[r]^{j} \ar@{->>}[d]_{f} & {X} \\
      {\A^n_k}
    }
  \]
  of $k$-schemes of finite type with an open immersion $j$ and a finite and surjective $f$.
  But then \ref{varietydim} follows from \ref{satzP} and from the transfer of \cite[4.1.2]{ega42}.

  For \ref{varietyinteger},
  note that we only have to prove the claim for "reduced" and "irreducible",
  since "integer" is just the conjunction of those two.

  Let us first consider the case where $X=\spec{k[X_1,\ldots,X_n]/\cI}$
  is affine. We have
  \begin{multline*}
    \text{$X$ reduced}
    \Longleftrightarrow\sqrt{\cI}=\cI
    \Longleftrightarrow\sqrt{\cI}/\cI=(0)
    \stackrel{\ref{thmvdDries}\ref{vdDriesff}}{\Longleftrightarrow}
      \sqrt{\cI}\cdot k\str{[X_i]}/\cI\cdot k\str{[X_i]}=(0) \\
    \stackrel{\ref{bemideal}}{\Longleftrightarrow}\str{\bigl[\sqrt{\cI}\bigr]}=\str{\cI}
    \stackrel{\ref{thmvdDries}\ref{vdDriesrad}}{\Longleftrightarrow}\sqrt[*]{\str{\cI}}=\str{\cI}
    \Longleftrightarrow\text{$\s{X}$ *reduced}
  \end{multline*}
  and
  \begin{multline*}
    \text{$X$ irreducible}
    \Longleftrightarrow
    \text{$\cI$ has exactly one minimal prime ideal} \\
    \stackrel{\ref{thmvdDries}\ref{vdDriesminprime}}{\Longleftrightarrow}
      \text{$\str{\cI}$ has exactly one minimal prime *ideal}
    \Longleftrightarrow\text{$\s{X}$ *irreducible}.
  \end{multline*}
  In the general case, let $(U_j)_{j\in\cJ}$ be a finite open covering of $X$
  by affine schemes $U_j$ which are not empty.
  The scheme $X$ is reduced if and only if
  the $U_j$ are reduced, which we have just proven to be equivalent to the $\s{U_j}$ being *reduced,
  which in turn is equivalent to $\s{X}$ being *reduced by $\ref{corcov}$ and transfer.

  Let $X$ be irreducible. Then all $U_j$ are irreducible, and their intersection is an open non-empty subscheme
  of $X$. Then by \ref{varietyzero}, the *scheme
  \[
    \s{\bigcap_{j\in\cJ}U_j}\stackrel{\ref{satzfin}\ref{finleftex}}{=}\bigcap_{j\in\cJ}\s{U_j}
  \]
  is not *empty.
  Since we already know that the $\s{U_j}$ are *irreducible and therefore *connected,
  this implies that $\s{X}$ is *connected.

  Assume that $\s{X}$ is *reducible. Since $\s{X}$ is *connected,
  there must be a *topological point of $\s{X}$ where two *irreducible components of $\s{X}$ intersect,
  and since the $\s{U_j}$ cover $\s{X}$, this *topological point lies in one of the $U_j$ which consequently
  can not be *irreducible, a contradiction.

  Now let $\s{X}$ be *irreducible, and assume that $X$ is not irreducible.
  Since $\s{X}$ is *irreducible, the $\s{U}$ are *irreducible, and their intersection is not *empty,
  so by \ref{varietyzero} and \ref{satzfin}\ref{finleftex}, the scheme $X$ is connected.
  Reasoning as above, we see this implies that one of the $U_j$ is reducible,
  which contradicts the fact that the $\s{U_j}$ are *irreducible.

  For \ref{varietycohex}, we have to show that
  a short sequence of coherent $\cO_X$-modules
  \[
    0\longrightarrow
    \cF'\longrightarrow
    \cF\longrightarrow
    \cF''\longrightarrow
    0
  \]
  is exact if and only if the induced sequence of *coherent $\cO_{\s{X}}$-modules
  \[
    0\longrightarrow
    \s{\cF'}\longrightarrow
    \s{\cF}\longrightarrow
    \s{\cF''}\longrightarrow
    0
  \]
  is exact,
  which by \ref{satzslsBex} and \ref{satzTNlsNiso} is equivalent to the exactness of
  \[
    0\longrightarrow
    \cF'\otimes_{\cO_X}\cO_{\hat{X}}\longrightarrow
    \cF\otimes_{\cO_X}\cO_{\hat{X}}\longrightarrow
    \cF''\otimes_{\cO_X}\cO_{\hat{X}}\longrightarrow
    0.
  \]
  Taking stalks,
  it is enough to show that for every point $x\in X$,
  \[
    0\longrightarrow
    \cF'_x\longrightarrow
    \cF_x\longrightarrow
    \cF''_x\longrightarrow
    0
  \]
  is exact if and only if
  \[
    0\longrightarrow
    \cF'_x\otimes_{\cO_{X,x}}\cO_{\hat{X},x}\longrightarrow
    \cF\otimes_{\cO_{X,x}}\cO_{\hat{X},x}\longrightarrow
    \cF''\otimes_{\cO_{X,x}}\cO_{\hat{X},x}\longrightarrow
    0
  \]
  is exact.
  But since
  \[
    \cO_{\hat{X},x}=
    \varinjlim_{x\in U\subseteq X}\cO_{\hat{X}}(U)=
    \varinjlim_{x\in U\subseteq X}\s{[\cO_X(U)]},
  \]
  where the limit is taken over all affine neighborhoods of $x$ in $X$,
  we see from \ref{corff} that
  $\cO_{X,x}\longrightarrow\cO_{\hat{X},x}$ is faithfully flat,
  and the claim follows.
\end{proof}

\vspace{\abstand}

\begin{cor}\label{corkhom}
  Let $X$ be a $k$-scheme of finite type,
  and let $\cF$ and $\cG$ be coherent $\cO_X$-modules.
  Then the canonical morphism \eqref{eqhom} is an isomorphism:
  \[
    \s{\HOM{X}{\cF}{\cG}}
    \stackrel{\sim}{\longrightarrow}
    \HOM{\s{X}}{\s{\cF}}{\s{\cG}}.
  \]
\end{cor}

\vspace{\abstand}

\begin{proof}
  Since the question is local,
  we can assume that there exists a global presentation
  \[
    \cO_X^m\longrightarrow
    \cO_X^n\longrightarrow
    \cF\longrightarrow
    0,
  \]
  and since $\s{}$ is exact by \ref{thmvariety}\ref{varietycohex},
  the functors $\s{\HOM{X}{\_}{\cG}}$ and $\HOM{\s{X}}{\s{\_}}{\s{\cG}}$ from
  $\Coh{X}$ to $\sCoh{\s{X}}$ are both left exact,
  so that we get the following commutative diagram with exact rows:
  \[
    \xymatrix{
      0 \ar[r] &
      {\s{\HOM{X}{\cF}{\cG}}} \ar[r] \ar@{.>}[d]_{\alpha} &
      {\s{\HOM{X}{\cO_X^n}{\cG}}} \ar[r] \ar[d]^{\wr}_{\beta} &
      {\s{\HOM{X}{\cO_X^m}{\cG}}} \ar[d]^{\wr}_{\gamma} \\
      0 \ar[r] &
      {\HOM{\s{X}}{\s{\cF}}{\s{\cG}}} \ar[r] &
      {\HOM{\s{X}}{\s{\cO_X^n}}{\s{\cG}}} \ar[r] &
      {\HOM{\s{X}}{\s{\cO_X^m}}{\s{\cG}}} \\
    }
  \]
  with the vertical morphisms given by \eqref{eqhom}.
  By \ref{corhom}, both $\beta$ and $\gamma$ are isomorphisms,
  so $\alpha$ must be an isomorphism as well.
\end{proof}

\vspace{\abstand}

\begin{cor}\label{corkNNcompat}
  Let $X$ be a $k$-scheme of finite type,
  and let $Z$ be a closed subscheme of $X$ corresponding to a sheaf of ideals $\cI$ on $X$.
  Then the *closed *subscheme $\s{Z}$ of $\s{X}$ is given by the *ideal $\s{\cI}$.
\end{cor}

\vspace{\abstand}

\begin{proof}
  According to \ref{satzNNcompat}, $\s{Z}$ is given by the *ideal
  $\image{\s{\cI}}{\cO_{\s{X}}}$.
  But $\cI\rightarrow\cO_X$ is a monomorphism and $\s{}$ is exact by \ref{thmvariety}\ref{varietycohex},
  so $\s{\cI}\hookrightarrow\cO_{\s{X}}$, and the corollary follows.
\end{proof}

\vspace{\abstand}

\begin{lemma}\label{lemmaSaturation}
  Let $A$ be a finitely generated $k$-algebra,
  let $\cI$ be an ideal of $A$,
  and let $f\in A$.
  Consider the ideals $(\cI:f^n):=\{a\in A\vert\ af^n\in\cI\}$ (for $n\in\Np$) and
  $(\cI:f^\infty):=\bigcup_{n\in\Np}(\cI:f^n)$ of $A$.
  Then
  \begin{equation}\label{eqSat1}
    \forall n\in\Np:\ (\cI:f^n)\cdot\s{A}=(\cI\cdot\s{A}:f^n)
  \end{equation}
  and
  \begin{equation}\label{eqSat2}
    (\cI:f^\infty)\cdot\s{A}=\sideset{^*}{}\bigcup_{n\in\Nps}(\cI\cdot\s{A}:f^n)=(\cI\cdot\s{A}:f^\infty).
  \end{equation}
\end{lemma}

\vspace{\abstand}

\begin{proof}
  By definition, the diagram
  \[
    \xymatrix{
      {(\cI:f^n)} \ar[r] \ar[d] \ar@{}[rd]|{\Box} & {\cI} \ar@{^{(}->}[d] \\
      A \ar[r]_{f^n} & A \\
    }
  \]
  is cartesian in the category of $A$-modules for every $n\in\Np$. Since $A\xrightarrow{\sm{A}}\s{A}$
  is (faithfully) flat by \ref{corff}, this implies \eqref{eqSat1}.
  For \eqref{eqSat2} note that $(\cI:f^\infty)$ is finitely generated, because $A$ is noetherian.
  Consequently, there is an $N\in\Np$ with
  $(\cI:f^N)=(\cI:f^{N+1})$ and hence $(\cI:f^n)=(\cI:f^N)$ for all $n\geq N$ and
  $(\cI:f^\infty)=(\cI:f^N)$. Then
  \[
    (\cI\cdot\s{A}:f^N)
    \stackrel{\eqref{eqSat1}}{=}(\cI:f^N)\cdot\s{A}
    =(\cI:f^{N+1})\cdot\s{A}
    \stackrel{\eqref{eqSat1}}{=}(\cI\cdot\s{A}:f^{N+1}),
  \]
  and hence $(\cI\cdot\s{A}:f^n)=(\cI\cdot\s{A}:f^N)$ for all $\Nps\ni n\geq N$ by transfer ---
  so \eqref{eqSat2} holds.
\end{proof}

\vspace{\abstand}

\begin{satz}\label{satzClosure}
  Let $X$ be a $k$-scheme of finite type,
  let $Y\subseteq X$ be a subscheme,
  and let $\bar{Y}\subseteq X$ be the scheme theoretic closure of $Y$ in $X$.
  Then $\s{\bar{Y}}$ is the *scheme theoretic closure of $\s{Y}$ in $\s{X}$.
\end{satz}

\vspace{\abstand}

\begin{proof}
  If $U\subseteq X$ is an open subscheme,
  then $\overline{Y\cap U}$, the closure of $Y\cap U$ in $U$, equals $\bar{Y}\cap U$.
  Therefore we can assume without loss of generality that $X=\spec{A}$ is affine
  and that $Y=\spec{A/\cI}\cap\D{f_1}\cap\ldots\cap\D{f_n}$ for an ideal $\cI\subseteq A$
  and elements $f_1,\ldots,f_n\in A$.
  Then $\bar{Y}=\spec{A/\cJ}$ with
  $\cJ=\bigcap_{i=1}^n\ker\bigl[A\xrightarrow{\mathrm{can}}A_{f_i}/\cI A_{f_i}\bigr]$.
  For any $f\in A$, we have
  \[
    \ker\bigl[A\rightarrow A_f/\cI A_f\bigr]
    =\bigl\{a\in A\bigl\vert\ \exists n\in\Np:\ f^na\in\cI\bigr\}
    =\bigcup_{n=1}^\infty(\cI:f^n)
    =(\cI:f^\infty),
  \]
  so $\cJ=\bigcup_{i=1}^n(\cI:f_i^\infty)$.
  Let $\tilde{\cJ}\subseteq\s{A}$ be the *ideal corresponding to the *schema theoretic closure of
  $\s{Y}$ in $\s{X}$. By transfer, we have
  \[
    \tilde{\cJ}
    =\bigcap_{i=1}^n(\cI\cdot\s{A}:f_i^\infty)
    \stackrel{\ref{lemmaSaturation}}{=}\bigcap_{i=1}^n(\cI:f_i^\infty)\cdot\s{A}
    =\cJ\cdot\s{A},
  \]
  so $[\s{A}]/\tilde{\cJ}=\s{[A/\cJ]}$, and we are done.
\end{proof}

\vspace{\abstand}

\begin{satz}\label{satzBlowup}
  Let $X$ be a $k$-scheme of finite type,
  let $Y\subseteq X$ be a closed subscheme,
  and let $f:Z\rightarrow X$ be the blow-up of $X$ in $Y$.
  Then $\s{f}:\s{Z}\rightarrow\s{X}$ is the *blow-up of $\s{X}$ in $\s{Y}$.
\end{satz}

\vspace{\abstand}

\begin{proof}
  First note that $\s{Y}$ is a *closed *subscheme of $\s{X}$ by \ref{satzP},
  so the statement makes sense.
  Next, by \cite[8.8.2, 8.10.5]{ega43}
  there exist a finitely generated subring $k_0$ of $k$,
  a $k_0$-scheme $X_0$ of finite type with $X=X_0\times_{k_0}k$
  and a closed subscheme $Y_0$ of $X_0$ with $Y=Y_0\times_{k_0}k$.

  Let $Z_0\rightarrow X_0$ be the blow-up of $X_0$ in $Y_0$,
  and let $W$ be the cartesian diagram of $k$-schemes
  \begin{equation}\label{eqBlowup}
    \xymatrix{
      W \ar@{^{(}.>}[r] \ar@{.>}[d] \ar@{}[rd]|{\Box} &
      {X\times_{X_0}Z_0} \ar[d]^\pi \ar[r] \ar@{}[rd]|{\Box} &
      {Z_0} \ar[d]^{\pi_0} \\
      {X\setminus Y} \ar@{^{(}->}[r] &
      X \ar[r] &
      {X_0.} \\
    }
  \end{equation}
  If $Z:=\bar{W}$ denotes the scheme theoretic closure of $W$ in $X\times_{X_0}Z_0$, then
  $f:=\pi\vert_Z:Z\rightarrow X$ is the blow-up of $X$ in $Y$ (compare \cite[IV-21]{eisenbudharris}).

  Applying the functor $\s{}$ to the left square of \eqref{eqBlowup}
  and using \ref{satzfin}\ref{finleftex} and \ref{lemmaComplement},
  we get a cartesian square of $k$-*schemes
  \[
    \xymatrix{
      {\s{W}} \ar@{^{(}.>}[r] \ar@{.>}[d] \ar@{}[rd]|{\Box} &
      {[\s{X}]\times_{\str{X_0}}\str{Z_0}} \ar[d]^{\s{\pi}} \\
      {[\s{X}]\setminus[\s{Y}]} \ar@{^{(}->}[r] &
      {\s{X},} \\
    }
  \]
  and by transfer, the *blow-up of $\s{X}$ in $\s{Y}$ is the *scheme theoretic closure of $\s{W}$
  in $[\s{X}]\times_{\str{X_0}}\str{Z_0}$.
  But according to \ref{satzClosure}, this is just $\s{\bar{W}}=\s{Z}$,
  which completes the proof.
\end{proof}

\vspace{\abstand}

\begin{defi}\label{deftangent}
  For every field $K$, every $K$-scheme $X$ and every $K$-rational point $x\in X$,
  we have the $K$-vector space $\mathrm{T}_{X,x}$, the \emph{(Zariski) tangent space of $X$ at $x$},
  defined as the $K$-dual of $\m_x/\m_x^2$.

  By transfer, for every *field $K$, every *scheme $X$ over $K$ and every $K$-valued point $x$ of $X$,
  we thus have an internal $K$-vector space $\mathrm{T}_{X,x}$ which we also call the
  \emph{(Zariski) tangent space of $X$ at $x$}.
\end{defi}

\vspace{\abstand}

\begin{satz}\label{satztangent}
  let $X$ be a $k$-scheme of finite type,
  and let $x\in X$ be a $k$-rational point.
  Then $\s{}$ induces a canonical functorial $k$-isomorphism of Zariski tangent spaces
  \[
    \s{}:
    \mathrm{T}_{X,x}\stackrel{\sim}{\longrightarrow}\mathrm{T}_{\s{X},\s{x}}.
  \]
\end{satz}

\vspace{\abstand}

\begin{proof}
  Identify $x$ with a $k$-morphism $x:\spec{k}\rightarrow X$,
  and let $e:\spec{k}\rightarrow\spec{k[\varepsilon]/\varepsilon^2}$ be the $k$-morphism induced by sending
  $\varepsilon$ to zero.

  It is well known that there is a canonical functorial isomorphism of $k$-vector spaces
  \begin{equation}\label{eqtangent}
    \mathrm{T}_{X,x}\cong
    \bigl\{t\in X(k[\varepsilon]/\varepsilon^2)\ \bigl\vert\ e^*t=x\bigr\}.
  \end{equation}
  By transfer, we get a canonical functorial isomorphism of internal $k$-vector spaces
  \begin{equation}\label{eqStangent}
    \mathrm{T}_{\s{X},\s{x}}\cong
    \bigl\{t\in (\s{X})(k\str{[\varepsilon]}/\varepsilon^2)\ \bigl\vert\ (\s{e})^*t=\s{x}\bigr\}.
  \end{equation}
  But by \ref{satzfinite} we have $k\str{[\varepsilon]}/\varepsilon^2=k[\varepsilon]/\varepsilon^2$,
  and we get the following commutative diagram of sets:
  \[
    \xymatrix@C=20mm{
      {X(k[\varepsilon]/\varepsilon^2)} \ar[r]^{e^*} \ar[d]_{\s{}}^{\wr} & {X(k)} \ar[d]^{\s{}}_{\wr} \\
      {(\s{X})(k[\varepsilon]/\varepsilon^2)} \ar[r]^{(\s{e})^*} & {(\s{X})(k),} \\
    }
  \]
  where the vertical maps are bijections because of \ref{thmBS}.
  From this, \eqref{eqtangent} and \eqref{eqStangent} the claim immediately follows.
\end{proof}

\vspace{\abstand}

\begin{cor}\label{corNonSing}
  Assume that $k$ is *algebraically closed,
  and let $X$ be a $k$-scheme of finite type.
  If $\s{X}$ is *nonsingular, then $X$ is nonsingular.
\end{cor}

\vspace{\abstand}

\begin{proof}
  Let $d:=\dim X$,
  and let $x\in X$ be a closed point.
  Since $k$ is *algebraically closed, $k$ is externally an algebraically closed field,
  and $x$ is a $k$-rational point.
  Since $\s{X}$ is *nonsingular of *dimension $d$ (by \ref{thmvariety}\ref{varietydim}),
  the tangent space $\mathrm{T}_{\s{X},\s{x}}$ has *dimension $d$,
  and the tangent space $\mathrm{T}_{X,x}$ has dimension $d$ by \ref{satztangent}.
  This shows that all tangent spaces of $X$ at closed points have dimension $d$,
  which means that $X$ is nonsingular.
\end{proof}

\vspace{\abstand}

\begin{satz}\label{satzSubschemes}
  Let $X$ be a $k$-scheme of finite type,
  and let $Y$ and $Z$ be two subschemes of $X$.
  If $\s{Y}\hookrightarrow\s{X}$ factors through $\s{Z}\hookrightarrow\s{X}$,
  then $Y\hookrightarrow X$ factors through $Z\hookrightarrow X$.
  In particular, if $\s{Y}$ and $\s{Z}$ are the same *subschemes of $\s{X}$,
  then $Y$ and $Z$ are the same subschemes of $X$.
\end{satz}

\vspace{\abstand}

\begin{proof}
  Factor $Z\hookrightarrow X$ as $Z\stackrel{i_Z}{\hookrightarrow}V\stackrel{j_Z}{\hookrightarrow}X$
  with a closed immersion $i_Z$ and an open immersion $j_Z$.
  We claim that $Y\hookrightarrow X$ factors through $j_Z$:
  Equip $Y\setminus V$ with its reduced structure and consider the cartesian diagram
  \[
    \xymatrix{
      \emptyset \ar@{^{(}->}[r] \ar@{^{(}->}[d] \ar@{}[rd]|{\Box} &
      {Y\setminus V} \ar@{^{(}->}[d] \\
      V \ar@{^{(}->}[r]_j &
      {X.} \\
    }
  \]
  Applying $\s{}$ and using \ref{satzfin}\ref{finleftex} and \ref{satzfin}\ref{finempty},
  we get a cartesian diagram
  \begin{equation}\label{eq4}
    \xymatrix{
      \emptyset \ar@{^{(}->}[r] \ar@{^{(}->}[d] \ar@{}[rd]|{\Box} &
      {\s{[Y\setminus V]}} \ar@{^{(}->}[d] \\
      \s{V} \ar@{^{(}->}[r]_{\s{j}} &
      {\s{X}.} \\
    }
  \end{equation}
  If $Y\setminus V$ was not empty, then $\s{[Y\setminus V]}$ also would not be empty by \ref{thmvariety}\ref{varietyzero}.
  But a *point of $\s{[Y\setminus V]}$ is a point of $\s{Y}$ which --- because \eqref{eq4} is cartesian ---
  is not a point of $\s{V}$, a contradiction to the fact that
  $\s{Y}\subseteq\s{Z}\subseteq\s{V}$ by assumption.

  So without loss of generality (by replacing $X$ with $V$),
  we can assume that $Z$ is a \emph{closed} subscheme of $X$.
  Factoring $Y\hookrightarrow X$ as $Y\stackrel{i_Y}{\hookrightarrow}U\stackrel{j_Y}{\hookrightarrow}X$
  with a closed immersion $i_Y$ and an open immersion $j_Y$
  and replacing $X$ with $U$ and $Z$ with $Z\cap U$, we can furthermore assume that $Y$ is also a
  closed subscheme of $X$.

  Finally, since the question is local on $X$,
  we can assume that $X=\spec{A}$ is affine and that $Y$ and $Z$ are given by ideals $\cI$ and $\cJ$ of $A$.
  By assumption, we have $\cJ\cdot\s{A}\subseteq\cI\cdot\s{A}$, and using \ref{corff},
  we conclude
  \[
    \cJ
    =A\cap[\cJ\cdot\s{A}]
    \subseteq A\cap[\cI\cdot\s{A}]
    =\cI.
  \]
\end{proof}

\vspace{\abstand}

\begin{bem}\label{bemEq}
  Let $\cC$ be a category with fibred products and a terminal object $T$,
  let $X$ and $Y$ be two objects of $\cC$,
  and let $f,g:X\rightarrow Y$ be two morphisms.
  Then the \emph{equalizer}
  \[
    \xymatrix@C=20mm{
      {\Eq{f}{g}} \ar[r]^-{\eq{f}{g}} &
      X \ar@<1mm>[r]^f \ar@<-1mm>[r]_g &
      Y \\
    }
  \]
  of $f$ and $g$ exists --- it is given by the cartesian diagram
  \begin{equation}\label{eqEq}
    \xymatrix@R=15mm@C=15mm{
      {\Eq{f}{g}} \ar@{.>}[r]^-{\eq{f}{g}} \ar[d] \ar@{}[rd]|{\Box} &
      X \ar[d]^-{(f,g)} \\
      Y \ar[r]_-{\langle\id{Y},\id{Y}\rangle} &
      {Y\times_T Y} \\
    }
  \end{equation}
\end{bem}

\vspace{\abstand}

\begin{lemma}\label{lemmaEq}
  Let $S$ be a scheme,
  let $X$ and $Y$ be two $S$-schemes,
  and let $f,g:X\rightarrow Y$ be two $S$-morphisms.
  Then the equalizer
  \[
    \xymatrix@C=20mm{
      {\Eq{f}{g}} \ar[r]^-{\eq{f}{g}} &
      X \ar@<1mm>[r]^f \ar@<-1mm>[r]_g &
      Y \\
    }
  \]
  of $f$ and $g$ exists in the category of $S$-schemes and is an immersion.
\end{lemma}

\vspace{\abstand}

\begin{proof}
  The category of $S$-schemes has fibred products and the terminal object $S$,
  so the equalizer of $f$ and $g$ exists by \ref{bemEq}.
  It is an immersion by the construction given in \eqref{eqEq},
  because
  $Y\xrightarrow{\langle\id{Y},\id{Y}\rangle}Y\times_SY$ is an immersion.
\end{proof}

\vspace{\abstand}

\begin{cor}\label{corEq}
  The functor $\s{}:\Schfp{k}\longrightarrow\sSchfp{k}$ is faithful.
\end{cor}

\vspace{\abstand}

\begin{proof}
  Let $X$ and $Y$ be $k$-schemes of finite type, and
  let $f,g:X\rightarrow Y$ be $k$-morphisms with $\s{f}=\s{g}$.
  By \ref{corEq}, $f$ and $g$ are equal if and only if
  $\Eq{f}{g}$ equals $X$ as subschemes of $X$.
  By assumption, $\Eq{\s{f}}{\s{g}}$ is the *subscheme $\s{X}$ of $\s{X}$,
  and $\Eq{\s{f}}{\s{g}}=\s{\Eq{f}{g}}$ by \ref{satzfin}\ref{finleftex},
  so the claim follows from \ref{satzSubschemes}.
\end{proof}

\vspace{\abstand}

Let $S$ be a noetherian scheme,
let $X/S$ be projective with very ample sheaf $\cO(1)$,
let $\cF$ be a coherent sheaf on $X$,
and let $P\in\Q[t]$ be a rational polynomial.
Then we have the \emph{Quot-scheme} $\Quot{P}{\cF}{X}{S}$,
projective over $S$,
which represents the contravariant functor $T\mapsto\FQuot{P}{\cF_{X\times_ST}}{X_T}{T}$
that maps a locally noetherian $S$-scheme $T$ to the set of those quotients $\cF_{X\times_ST}\twoheadrightarrow\cG$
with $\cG$ flat over $T$ and Hilbert polynomial $P$ in every fibre $t\in T$
(compare \cite[221.3]{fga}).

By transfer, for a *noetherian *scheme $S$,
a *projective $S$-*scheme $X$ with *very ample *sheaf $\cO(1)$,
a *coherent *sheaf $\cF$ on $X$
and a *polynomial $P\in\Qs\str{[t]}$,
we have a canonical *projective $S$-*scheme $\sQuot{P}{\cF}{X}{S}$
which represents the enlarged functor $T\mapsto\sFQuot{P}{\cF_{X\times_ST}}{X_T}{T}$
on *locally noetherian $S$-*schemes.

In the special case $\cF=\cO_X$,
the Quot-scheme $\Quot{P}{\cO_X}{X}{S}$ is called the \emph{Hilbert scheme} and denoted by $\Hilb{P}{X}{S}$
(its $T$-valued points correspond to closed subschemes of $X_T$ which are flat over $T$ and have Hilbert polynomial $P$
in every fibre).
--- Similarly, we call $\sHilb{P}{X}{S}:=\sQuot{P}{\cO_X}{X}{S}$ the \emph{*Hilbert scheme}.

In the following proposition, we want to show that the formation of Quot-schemes and Hilbert schemes
is compatible with the functor $\s{}$:

\begin{satz}\label{satzQuot}
  Let $X$ be a projective $k$-scheme with very ample sheaf $\cO(1)$,
  let $\cF$ be a coherent sheaf on $X$,
  and let $P\in\Q[t]$ be a rational polynomial.
  \begin{enumerate}
    \item\label{Quot1}
      We have $\s{\Quot{P}{\cF}{X}{k}}=\sQuot{P}{\s{\cF}}{\s{X}}{k}$
      and in particular $\s{\Hilb{P}{X}{k}}=\sHilb{P}{\s{X}}{k}$,
      where $P$ is considered as a *polynomial via $\Q[t]\hookrightarrow\Qs\str{[t]}$.
    \item\label{Quot2}
      Let $T$ be a $k$-scheme of finite type,
      and let $f:T\rightarrow\Quot{P}{\cF}{X}{k}$ be a $T$-valued point,
      corresponding to a quotient $\varphi:\cF_{X\times_kT}\twoheadrightarrow\cG$.
      Then $[\s{f}]$, which is a $[\s{T}]$-valued point of $\sQuot{P}{\s{\cF}}{\s{X}}{k}$
      by \ref{Quot1}, corresponds to the quotient
      $[\s{\cF}]_{[\s{X}]\times_k[\s{T}]}\twoheadrightarrow\s{G}$.

      In particular, if $g:T\rightarrow\Hilb{P}{X}{k}$ corresponds to the subscheme $Z\subseteq X\times_kT$,
      then $[\s{g}]$ corresponds to the *subscheme $[\s{Z}]\subseteq[\s{X}]\times_k[\s{T}]$.
  \end{enumerate}
\end{satz}

\vspace{\abstand}

\begin{proof}
  By \cite[8.5.2, 8.8.2, 8.10.5]{ega43},
  there exist a finitely generated subring $A_0$ of $k$,
  a projective $A_0$-scheme $X_0$ with $X=X_0\times_{A_0}k$
  and a coherent sheaf $\cF_0$ on $X_0$ with $[A_0\hookrightarrow k]^*\cF_0=\cF$.

  Then $\Quot{P}{\cF}{X}{k}=\Quot{P}{\cF_0}{X_0}{A_0}\times_{A_0}k$,
  and putting $\alpha:=\bsi{\Z}{A_0}{k}[A_0\hookrightarrow k]:\str{A_0}\rightarrow k$,
  we get
  \begin{multline*}
    \s{\Quot{P}{\cF}{X}{k}}
    \stackrel{\ref{thmS}}{=}\alpha^*\Bigl(\str{\bigl[\Quot{P}{\cF_0}{X_0}{A_0}\bigr]}\Bigr)
    =\alpha^*\Bigl[\sQuot{P}{\str{\cF_0}}{\str{X_0}}{\str{A_0}}\Bigr] \\
    =\sQuot{P}{\alpha^*[\str{\cF_0}]}{\alpha^*[\str{X_0}]}{\alpha^*[\str{A_0}]}
    \stackrel{\ref{thmS},\ref{thmScoh}}{=}\sQuot{P}{\s{\cF}}{\s{X}}{k},
  \end{multline*}
  which settles \ref{Quot1}.

  By \cite[8.8.2]{ega43}, after a possible change of $A_0$, $X_0$ and $\cF_0$,
  we find an $A_0$-scheme $T_0$ of finite type with $T=T_0\times_{A_0}k$
  and an $A_0$-morphism $f_0:T_0\rightarrow Q:=\Quot{P}{\cF_0}{X_0}{A_0}$ with
  $[A_0\hookrightarrow k]^*f_0=f$. Let $\cF_{X_0\times_{A_0}Q}\twoheadrightarrow\Guniv$
  be the universal quotient. Then $f_0$ corresponds to the quotient
  $[\cF_0]_{X_0\times_{A_0}T_0}\twoheadrightarrow[\id{X_0}\times f_0]^*\Guniv=:\cG_0$,
  and $f$ corresponds to the quotient
  $\varphi:\cF_{X\times_kT}\twoheadrightarrow[A_0\hookrightarrow k]^*\cG_0=\cG$.
  So
  \[
    \s{\cG}
    \stackrel{\ref{thmScoh}}{=}\alpha^*\bigl(\str{\cG_0}\bigr)
    =\alpha^*\Bigl(\str{\Bigl[[\id{X_0}\times f_0]^*\Guniv\Bigr]}\Bigr)
    \stackrel{\ref{thmS}}{=}\bigl[\id{\s{X}}\times[\s{f}]\bigr]^*\bigl(\alpha^*(\str{\Guniv})\bigr),
  \]
  and since $[\s{\cF}]_{[\s{X}]\times_k\sQuot{P}{\s{\cF}}{\s{X}}{k}}\twoheadrightarrow\alpha^*(\str{\Guniv})$
  obviously is the universal quotient, this proves \ref{Quot2}.
\end{proof}

\vspace{\abstand}

\begin{cor}\label{corhilb}
  Let $X$ be a projective $k$-scheme with very ample sheaf $\cO(1)$,
  and let $\cF$ be a coherent $\cO_X$-module.

  Then the Hilbert polynomial of $\cF$ (with respect to $\cO(1)$)
  coincides with the *Hilbert polynomial of $\s{\cF}$ (with respect to $\s{\cO(1)}$) in $\Qs\str{[t]}$.
\end{cor}

\vspace{\abstand}

\begin{proof}
  Denote the Hilbert polynomial of $\cF$ by $P_{\cF}\in\Q[t]\subset\Qs\str{[t]}$.
  If $\cF$ corresponds to the $k$-valued point $f$ of $\Quot{P_{\cF}}{\cF}{X}{k}$,
  then $\s{\cF}$ corresponds to the $k$-valued point $[\s{f}]$
  of $\sQuot{P_{\cF}}{\s{\cF}}{\s{X}}{k}$ according to \ref{satzQuot}\ref{Quot2}.
  But by its very definition, $\sQuot{P_{\cF}}{\s{\cF}}{\s{X}}{k}$
  parameterizes sheaves with *Hilbert polynomial $P_{\cF}\in\Qs\str{[t]}$, so we are done.
\end{proof}

\vspace{\abstand}

\begin{thm}\label{thmssheaf}
  Let $X$ be a projective $k$-scheme with very ample sheaf $\cO(1)$,
  and let $\cG$ be a *coherent *sheaf on $\s{X}$.
  Then the following two statements are equivalent:
  \begin{enumerate}
    \item\label{ssheafi}
      There is a coherent sheaf $\cH$ on $X$ with $\s{\cH}\cong\cG$.
    \item\label{ssheafii}
      There is a coherent sheaf $\cF$ on $X$,
      such that $\cG$ is  a quotient of $\s{\cF}$,
      and the *Hilbert polynomial of $\cG$ (with respect to $\s{\cO(1)}$) lies in $\Q[t]\subset\Qs\str{[t]}$.
  \end{enumerate}
\end{thm}

\vspace{\abstand}

\begin{proof}
  The implication ``\ref{ssheafi}$\Rightarrow$\ref{ssheafii}" is easy: We can simply put $\cF:=\cH$,
  and by \ref{corhilb},
  the *Hilbert polynomial of $\cG\cong\s{\cH}$ equals the Hilbert polynomial of $\cH$ and consequently
  lies in $\Q[t]$.

  For ``\ref{ssheafii}$\Rightarrow$\ref{ssheafi}", let $P\in\Q[t]\subset\Qs\str{[t]}$ be the *Hilbert polynomial
  of $\cG$. Then $\s{\cF}\twoheadrightarrow\cG$ corresponds to a $k$-valued point $g$ of
  $\sQuot{P}{\s{\cF}}{\s{X}}{k}$.
  Since
  \[
    \Bigl[\Quot{P}{\cF}{X}{k}\Bigr](k)
    \stackrel{\s{}}{\longrightarrow}\Bigl[\s{\Quot{P}{\cF}{X}{k}}\Bigr](k)
    \stackrel{\ref{satzQuot}\ref{Quot1}}{=}\Bigl[\sQuot{P}{\s{\cF}}{\s{X}}{k}\Bigr](k)
  \]
  is bijective by \ref{thmBS},
  there exists a $k$-valued point $h$ of $\Quot{P}{\cF}{X}{k}$ with $g=\s{h}$.
  If $\cF\twoheadrightarrow\cH$ is the quotient given by $h$,
  then $\s{\cH}\cong\cG$ by \ref{satzQuot}\ref{Quot2}.
\end{proof}

\vspace{\abstand}

\begin{cor}\label{corssheaf}
  Let $X$ be a projective $k$-scheme with very ample sheaf $\cO(1)$,
  and let $Z$ be a *closed *subscheme of $\s{X}$.
  Then the following two statements are equivalent:
  \begin{enumerate}
    \item
      There is a closed subscheme $W$ of $X$ with $\s{W}=Z$.
    \item
      The *Hilbert polynomial of $Z$ (with respect to $\s{\cO(1)}$) lies in $\Q[t]\subset\Qs\str{[t]}$.
  \end{enumerate}
\end{cor}

\vspace{\abstand}

\begin{proof}
  This follows immediately from \ref{thmssheaf},
  applied to the special case $\cG:=\cO_Z$ and $\cF:=\cO_X$.
\end{proof}

\vspace{\abstand}

\begin{cor}\label{corDegDim}
  Let $X$ be a projective $k$-scheme with very ample sheaf $\cO(1)$,
  and let $Z$ be a *closed *integral *subscheme (i.e. a *prime cycle) of $\s{X}$
  that has \emph{finite} *degree (with respect to $[Z\hookrightarrow\s{X}]^*\cO(1)$).
  Then there exists an integral subscheme (i.e. a prime cycle) $W$ of $X$ with $\s{W}=Z$.
\end{cor}

\vspace{\abstand}

\begin{proof}
  As $Z$ is a subscheme of $\s{X}$,
  we have $\str{\dim Z}\leq\str{\dim\s{X}}\stackrel{\ref{thmvariety}\ref{varietydim}}{=}\dim X$,
  so $Z$ is a *projective *integral *scheme of finite *degree and of \emph{finite} *dimension.
  Then transfer of \cite[XIII.6.11(i)]{sga6} shows that the *Hilbert polynomial of $Z$ has finite coefficients
  and consequently lies in $\Q[t]\subset\Qs\str{[t]}$,
  and the corollary follows from \ref{corssheaf}.
\end{proof}

\vspace{\abstand}

\begin{cor}\label{corDegDimP}
  Let $n\in\N_+$,
  and let $Z$ be a *integral *closed *subscheme of $\str{\P^n_k}$ of finite *degree.
  Then there is an integral closed subscheme $W$ of $\P^n_k$ with $\s{W}=Z$.
\end{cor}

\vspace{\abstand}

\begin{proof}
  This follows immediately from \ref{corDegDim} for $X:=\P^n_k$ and $\cO(1):=\s{\cO_{\P^n_k}(1)}$.
\end{proof}

\vspace{\abstand}

Let $S$ be a scheme,
and let $f:X\rightarrow Y$ be an $S$-morphism.
Then the \emph{graph of $f$} is the $S$-morphism $\Gamma_f:X\xrightarrow{\langle\id{X},f\rangle}X\times_SY$.
It is easy to see that the diagram
\[
  \xymatrix{
    X \ar[r]^-{\Gamma_f} \ar[d]_f \ar@{}[rd]|\Box &
    {X\times_SY} \ar[d]^{f\times\id{Y}} \\
    Y \ar[r]_-{\langle\id{Y},\id{Y}\rangle} &
    {Y\times_SY} \\
  }
\]
is cartesian, which shows that $\Gamma_f$ is an immersion (and can hence be considered as a subscheme of
$X\times_SY$, isomorphic to $X$),
which is closed if $Y/S$ is separated.

Now let $S$ be noetherian,
let $X$ and $Y$ be projective $S$-schemes with $X/S$ flat,
let $\cO(1)$ be a very ample sheaf on $X\times_SY$,
and let $P\in\Q[t]$ be a polynomial.
Consider the functor $T\mapsto\Homf{P}{S}{X}{Y}(T)$ that maps an $S$-scheme $T$ to the set of
those $T$-morphisms $f:X\times_ST\rightarrow Y\times_ST$ whose
graph $\Gamma_f\hookrightarrow X\times_SY$, a closed subscheme since $Y/S$ is separated,
has Hilbert polynomial $P$ with respect to $\cO(1)$.

It is well known (compare \cite[I.1.10]{kollar})
that this functor is represented by an open subscheme $\Homf{P}{S}{X}{Y}$ of $\Hilb{P}{X\times_SY}{S}$,
where $\Homf{P}{S}{X}{Y}\hookrightarrow\Hilb{P}{X\times_SY}{S}$ is given by sending a morphism to its graph.
Similar to the case of Quot- and Hilbert schemes, the formation of $\Homf{P}{S}{X}{Y}$ is compatible
with the functor $\s{}$ in the following sense:

\begin{satz}\label{satzHomf}
  Let $X$ and $Y$ be projective $k$-schemes,
  let $\cO(1)$ be a very ample sheaf on $X\times_kY$,
  and let $P\in\Q[t]$ be a rational polynomial.
  \begin{enumerate}
    \item\label{Homf1}
      We have $\s{\Homf{P}{k}{X}{Y}}=\sHomf{P}{k}{\s{X}}{\s{Y}}$,
      where $P$ is considered as a *polynomial via $\Q[t]\hookrightarrow\Qs\str{[t]}$.
    \item\label{Homf2}
      Let $T$ be a $k$-scheme of finite type,
      and let $f:T\rightarrow\Homf{P}{k}{X}{Y}$ be a $T$-valued point,
      corresponding to a $T$-morphism $g:X\times_kT\rightarrow Y\times_kT$.
      Then $[\s{f}]$, which is a $[\s{T}]$-valued point of $\sHomf{P}{k}{\s{X}}{\s{Y}}$
      by \ref{Homf1}, corresponds to the morphism
      $[\s{g}]:\s{X}\str{\times}_k\s{T}\rightarrow\s{Y}\str{\times}_k\s{T}$.
  \end{enumerate}
\end{satz}

\vspace{\abstand}

\begin{proof}
  This is completely analogous to the proof of \ref{satzQuot}.
\end{proof}

\vspace{\abstand}

\begin{thm}\label{thmsmorphism}
  Let $X$ and $Y$ be projective $k$-schemes,
  let $\cO(1)$ be a very ample sheaf on $X\times_kY$,
  and let $g:\s{X}\rightarrow\s{Y}$ be a morphism of $k$-*schemes.
  Then the following two statements are equivalent:
  \begin{enumerate}
    \item
      There is a $k$-morphism $f:X\rightarrow Y$ with $\s{f}=g$.
    \item
      The *Hilbert polynomial of the *graph of $g$ (with respect to $\s{\cO(1)}$) lies in $\Q[t]\subset\Qs\str{[t]}$.
  \end{enumerate}
\end{thm}

\vspace{\abstand}

\begin{proof}
  This follows from \ref{satzHomf}
  in the same way as \ref{thmssheaf} follows from \ref{satzQuot}.
\end{proof}

\vspace{\abstand}

\begin{cor}\label{corDegDimMorph}
  Let $X$ and $Y$ be projective $k$-schemes with $X$ integral,
  let $\cO(1)$ be a very ample sheaf on $X\times_kY$,
  and let $g:\s{X}\rightarrow\s{Y}$ be a morphism of $k$-*schemes
  whose *graph has \emph{finite} degree (with respect to $\s{\cO(1)}$).
  Then there exists a $k$-morphism $f:X\rightarrow Y$ with $\s{f}=g$.
\end{cor}

\vspace{\abstand}

\begin{proof}
  By transfer, the *graph $\str{\Gamma_g}$ of $g$ is isomorphic to $\s{X}$ and hence *integral.
  Then by \ref{corDegDim}, there is a closed subscheme $\Gamma$ of $X\times_kY$ with $\s{\Gamma}=\str{\Gamma_g}$,
  and it follows from \ref{corssheaf} that the *Hilbert polynomial of $\str{\Gamma_g}$ lies in $\Q[t]$.
  Then the corollary follows from \ref{thmsmorphism}.
\end{proof}

\vspace{\abstand}

\begin{cor}\label{corReflectsIsos}
  The restriction of $\s{}:\Schfp{k}\rightarrow\sSchfp{k}$ to the full subcategory of projective $k$-schemes
  reflects isomorphisms.
\end{cor}

\vspace{\abstand}

\begin{proof}
  Let $f:X\rightarrow Y$ be a morphism of projective $k$-schemes
  such that $\s{f}$ is an isomorphism with inverse $\tilde{g}:\s{Y}\rightarrow\s{X}$.
  Choose a very ample sheaf $\cO(1)$ on $X\times_kY$. If $\tau:Y\times_kX\xrightarrow{\sim}X\times_kY$
  denotes the transposition, $\tau^*\cO(1)$ is a very ample sheaf on $Y\times_kX$.
  Let $P\in\Q[t]$ be the Hilbert polynomial of $\Gamma_f$ with respect to $\cO(1)$,
  which by \ref{corhilb} is also the *Hilbert polynomial of $\str{\Gamma_{\s{f}}}$, the *graph of $\s{f}$,
  with respect to $\s{\cO(1)}$.
  If follows from transfer that the transpose $[\s{\tau}]^*[\str{\Gamma_{\s{f}}}]$ is the *graph of $\tilde{g}$
  and that its *Hilbert polynomial with respect to $\s[\tau^*\cO(1)]$ equals $P$.
  Thus by \ref{thmsmorphism}, there is a $k$-morphism $g:Y\rightarrow X$ with $\s{g}=\tilde{g}$.
  Now
  \[
    \s{[f\circ g]}=[\s{f}]\circ[\s{g}]=[\s{f}]\circ\tilde{g}=\id{\s{Y}}=\s{\id{Y}}
  \]
  and
  \[
    \s{[g\circ f]}=[\s{g}]\circ[\s{f}]=\tilde{g}\circ[\s{f}]=\id{\s{X}}=\s{\id{X}},
  \]
  so $f\circ g=\id{Y}$ and $g\circ f=\id{X}$ (because $\s{}:\Schfp{k}\rightarrow\sSchfp{k}$ is faithful by
  \ref{corEq}),
  and we see that $f$ is indeed an isomorphism (with inverse $g$).
\end{proof}

\vspace{\abstand}

\begin{lemma}\label{lemmaDegree}
  Let $\varphi:B\hookrightarrow C$ be a finite, injective morphism of integral $k$-algebras of finite type.
  Then $\s{\varphi}:\s{B}\rightarrow\s{C}$ is an injective, finite morphism of integral $k$-algebras,
  and
  \[
    [\QK{\s{C}}:\QK{\s{B}}]=[\QK{C}:\QK{B}]\in\Np.
  \]
\end{lemma}

\vspace{\abstand}

\begin{proof}
  The (internal) $k$-algebras $\s{B}$ and $\s{C}$ are integral by \ref{thmvariety}\ref{varietyinteger},
  and $\s{\varphi}$ is injective and finite, because
  \[
    \xymatrix{
      B \ar@{^{(}->}[r]^\varphi \ar[d]_{\sm{B}} & C \ar[d]^{\sm{C}} \\
      {\s{B}} \ar[r]_{\s{\varphi}} & {\s{C}} \\
    }
  \]
  is cocartesian by \ref{satzfinite} and because $\sm{B}$ is faithfully flat by \ref{corff}.
  Since $\QK{C}=C\otimes_B\QK{B}$ and
  \begin{multline*}
    \QK{\s{C}}
    =[\s{C}]\otimes_{\s{B}}\QK{\s{B}}
    \stackrel{\ref{satzfinite}}{=}C\otimes_B\QK{\s{B}} \\
    =\bigl[C\otimes_B\QK{B}\bigr]\otimes_{\QK{B}}\QK{\s{B}}
    =\QK{C}\otimes_{\QK{B}}\QK{\s{B}}.
  \end{multline*}
  Using this, we get
  \begin{multline*}
    [\QK{\s{C}}:\QK{\s{B}}]
    =\dim_{\QK{\s{B}}}\QK{\s{C}} \\
    =\dim_{\QK{\s{B}}}\bigl[\QK{C}\otimes_{\QK{B}}\QK{\s{B}}\bigr]
    =\dim_{\QK{B}}\QK{C}
    =[\QK{C}:\QK{B}],
  \end{multline*}
  and this degree is of course finite, because $\varphi$ is finite.
\end{proof}

\vspace{\abstand}

\begin{satz}\label{satzVarBirat}
  Let $f:X\rightarrow Y$ be a morphism of integral $k$-schemes of finite type.
  Then $f$ is birational if and only if $\s{f}:\s{X}\rightarrow\s{Y}$ is *birational.
\end{satz}

\vspace{\abstand}

\begin{proof}
  Assume first that $f$ is birational.
  Then by definition, there is a commutative diagram
  \[
    \xymatrix{
      & U \ar@{_{(}->}[dl]_{j_1} \ar@{^{(}->}[dr]^{j_2} \\
      X \ar[rr]_f & & Y \\
    }
  \]
  of $k$-morphisms with open immersions $j_1$ and $j_2$. So
  \[
    \xymatrix{
      & {\s{U}} \ar@{_{(}->}[dl]_{\s{j_1}} \ar@{^{(}->}[dr]^{\s{j_2}} \\
      {\s{X}} \ar[rr]_{\s{f}} & & {\s{Y}} \\
    }
  \]
  is a commutative diagram of $k$-*schemes,
  where $\s{j_1}$ and $\s{j_2}$ are *open immersions by \ref{satzP},
  which shows that $\s{f}$ is *birational.

  For the other implication, assume now that $\s{f}$ is *birational.
  Then $\s{X}$ and $\s{Y}$ have the same *dimension,
  and \ref{thmvariety}\ref{varietydim} implies that $\dim X=\dim Y$.
  Let us first show that $f$ is dominant: If it were not, there would be a non-empty open subscheme $U$ of $Y$
  and a cartesian diagram
  \[
    \xymatrix{
      \emptyset \ar[r] \ar@{}[rd]|\Box \ar@{^{(}->}[d] &
      U \ar@^{^{(}->}[d] \\
      X \ar[r]_f &
      {Y.} \\
    }
  \]
  But then by \ref{satzfin}\ref{finleftex}, \ref{satzP} and \ref{thmvariety}\ref{varietyzero},
  $\s{U}$ would be a non-empty *open *subscheme of $\s{Y}$ disjoint from $[\s{f}](\s{X})$;
  this means that $\s{f}$ would not be *dominant and consequently could not be *birational --- a contradiction.
  So $\varphi$ is indeed dominant,
  and if we denote the generic points of $X$ and $Y$ by $\xi$ respectively $\eta$,
  then $\xi$ is contained in the generic fibre $X_\eta$.
  Since we saw above that $\dim X=\dim Y$, we must have $X_\eta=\{\xi\}$ by \cite[4.1.2(i)]{ega42}.

  In particular, $X_\eta/\eta$ is of finite type and discrete,
  so by \cite[6.4.4]{ega1} it is finite.
  Then by \cite[p. 6 and 8.10.5(x)]{ega43},
  there is an affine, open, dense subset $V=\spec{B}\subseteq Y$, such that
  $f\vert_U:U\rightarrow V$ (with $U:=f^{-1}(V)$) is finite.
  Then $U=\spec{C}$ is affine, and $f^*:B\rightarrow C$ is a finite, injective morphism of integral $k$-algebras
  of finite type.
  By hypothesis we have $\QK{\s{B}}\xrightarrow{\sim}\QK{\s{C}}$,
  so \ref{lemmaDegree} implies $k(Y)=\QK{B}\xrightarrow{\sim}\QK{C}=k(X)$,
  which means that $f$ induces an isomorphism of the function fields of $X$ and $Y$
  and is therefore birational.
\end{proof}

\vspace{\abstand}


\section{The coherence theorem}

\vspace{\abstand}

For any scheme $X$, sheaf of $\cO_X$-modules $\cF$ and natural number $i\in\Nn$,
we can consider the Zariski cohomology group $\HZar{i}{X}{\cF}$.
If $X$ is an $A$-scheme for a ring $A$, then $\HZar{i}{X}{\cF}$ canonically carries the structure of an $A$-module.

If $f:X\rightarrow Y$ is a proper morphism of schemes
and if $\cF$ is a coherent $\cO_X$-module,
then we have the higher direct image $\Rf{i}{f}{\cF}$, a coherent $\cO_Y$-module
by \cite[3.2.1]{ega31}.

By transfer, if $X$ is a *scheme, $\cF$ a *finitely presented $\cO_X$-module and
$i\in\Nns$ a *natural number,
we get the \emph{*Zariski cohomology} $\HZar{i}{X}{\cF}$
which is an internal $A$-module if $X$ is an $A$-*scheme for a *ring $A$.

Similarly, if $f:X\rightarrow Y$ is a *proper morphism of *schemes
and if $\cF$ is a *coherent $\cO_X$-module,
we have the \emph{*higher direct image} $\Rf{i}{f}{\cF}$, a *coherent $\cO_Y$-module.

\vspace{\abstand}

\begin{lemma}\label{lemmaderex}
  Let $A$ be a *noetherian *ring,
  and let $f:X\rightarrow Y$ be a morphism of *schemes over $A$.
  Then the left exact functor $f_*:\sQCoh{X}\longrightarrow\sMod{Y}$
  factorizes over $\sQCoh{Y}$ and
  admits a right derived functor
  $\RF(f)_*:\Der{\sQCoh{X}}\longrightarrow\Der{\sQCoh{Y}}$.

  Furthermore, the class of flasque\footnote[2]{Note that being flasque is obviously first-order and hence
  is the same as being *flasque.}
  *quasi-coherent sheaves of $\cO_X$-*modules
  is adapted to $f_*$.
\end{lemma}

\vspace{\abstand}

\begin{proof}
  Let $B$ be a noetherian ring in $\cR$,
  and let $g:Z\rightarrow W$ be a morphism of finitely presented $B$-schemes.
  Then $Z$ and $g$ are quasi-separated (by \cite[1.2.8]{ega41}) and quasi-compact.
  It follows from \cite[B.3]{thomason_trobaugh} that $\QCoh{Z}^\cU$ has enough injective objects
  and from \cite[B.6]{thomason_trobaugh} that
  $\Rf{i}{g}{\cF}$ is quasi-coherent for all quasi-coherent $\cO_Z$-modules $\cF$ and all $i\in\Nn$.

  Furthermore, by \cite[B.4]{thomason_trobaugh}, an injective object in $\QCoh{Z}^\cU$ is
  also an injective (and hence flasque) object of $\Mod{Z}^\cU$,
  so that the class of flasque quasi-coherent $\cO_Z$-modules is adapted to $g_*$.

  Since all this is true for arbitrary $B$, $Z$, $W$ and $g$, the transferred statements are also true,
  and the lemma follows.
\end{proof}

\vspace{\abstand}

\begin{lemma}\label{lemmadercomm}
  Let $A$ be a *noetherian *ring, and let $f:X\rightarrow Y$ be a morphism of finitely presented $A$-schemes.
  Then the following diagram of exact functors between derived categories commutes
  (up to canonical isomorphism):
  \[
    \xymatrix@R=20mm@C=30mm{
      {\Der{\sQCoh{\s{X}}}} \ar[r]^{\RF[\s{f}]_*} \ar[d]_{\B{}} & {\Der{\sQCoh{\s{Y}}}} \ar[d]^{\B{}} \\
      {\Der{\Mod{X}}} \ar[r]_{\RF f_*} & {\Der{\Mod{Y}}.} \\
    }
  \]
\end{lemma}

\vspace{\abstand}

\begin{proof}
  First of all,
  note that $\RF{[\s{f}]_*}$ exists by \ref{lemmaderex}
  and that
  $\B{}:\sQCoh{\s{X}}\rightarrow\Mod{X}$ and $\B{}:\sQCoh{\s{Y}}\rightarrow\Mod{Y}$ are exact
  by \ref{satzslsBex}.

  The composition $\B{}\circ\RF{[\s{f}]_*}$ is canonically isomorphic to $\RF{[S\circ[\s{f}]_*]}$,
  because $\B{}$ is exact.
  The composition $\RF{f_*}\circ\B{}$ is canonically isomorphic to $\RF{[f_*\circ\B{}]}$,
  because $\B{}$ is exact and obviously maps flasque *sheaves
  to flasque sheaves, which are adapted to $f_*$.

  It follows immediately from the definition of $\B{}$, $f_*$ and $[\s{f}]_*$
  that $\B{}\circ[\s{f}]_*=f_*\circ\B{}$,
  so we have
  \[
    \B{}\circ\RF{[\s{f}]}_*
    \cong\RF[\B{}\circ[\s{f}]_*]
    =\RF[f_*\circ\B{}]
    \cong\RF{f_*}\circ\B{}.
  \]
\end{proof}

\vspace{\abstand}

\noindent
Let $k$ be a *field,
and let $f:X\longrightarrow Y$ be a \emph{proper} morphism of $k$-schemes of finite type.

\vspace{\abstand}

\begin{lemma}\label{lemmadercoh}
  We have a commutative diagram of exact functors
  \[
    \xymatrix{
      {\Derb{\sCoh{\s{X}}}} \ar[r]^{\RF{[\s{f}]_*}} \ar@{^{(}->}[d]_{\iota} &
      {\Derb{\sCoh{\s{Y}}}} \ar@{^{(}->}[d]^{\iota} \\
      {\Derb{\sQCoh{\s{X}}}} \ar[r]^{\RF{[\s{f}]_*}} \ar[d]_{\B{}} &
      {\Derb{\sQCoh{\s{Y}}}} \ar[d]^{\B{}} \\
      {\Derb{\Mod{X}}} \ar[r]_{\RF{f_*}} &
      {\Derb{\Mod{Y}}} \\
      {\Derb{\Coh{X}}} \ar[r]_{\RF{f_*}} \ar@{^{(}->}[u]^{\iota} &
      {\Derb{\Coh{Y}}} \ar@{^{(}->}[u]_{\iota} \\
    }
  \]
\end{lemma}

\vspace{\abstand}

\begin{proof}
  Since $X$ is finitely presented over a field, it is finite-dimensional,
  which implies that $f_*:\QCoh{X}\rightarrow\QCoh{Y}$ has finite cohomological dimension and hence
  induces $\RF{f_*}:\Derb{\QCoh{X}}\rightarrow\Derb{\QCoh{Y}}$.

  By \ref{thmvariety}\ref{varietydim} and transfer,
  $[\s{f}]_*:\sQCoh{\s{X}}\rightarrow\sQCoh{\s{Y}}$ has the same finite cohomological dimension
  and induces $\RF{[\s{f}]_*}:\Derb{\sQCoh{\s{X}}}\rightarrow\Derb{\sQCoh{\s{Y}}}$.
  So the middle square is well-defined, and it commutes by \ref{lemmadercomm}.

  The bottom square is well-defined and commutes by \cite[II.2.2]{hartshorne_residues}
  and \cite[II.2.2.2]{sga6},
  the top square is well-defined and commutes by transfer of
  \cite[II.2.2]{hartshorne_residues} and \cite[II.2.2.2.1]{sga6}.
\end{proof}

\vspace{\abstand}

\begin{satz}\label{satzcohmorph}
  There is a canonical morphism of exact functors
  \begin{equation}\label{eqcohmorph1}
    \xymatrix@C=30mm{
      {\Derb{\Coh{X}}}
        \ar@(ur,ul)[r]^{\iota\circ\RF{f_*}}
        \ar@(dr,dl)[r]_{\B{}\circ\RF{[\s{f}]_*}\circ\s{}}
        \ar@{}[r]|{\displaystyle\Downarrow} &
      {\Derb{\Mod{Y}}}
    }
  \end{equation}
  which induces a canonical morphism of $\delta$-functors
  \begin{equation}\label{eqcohmorph2}
    \xymatrix@C=30mm{
      {\Coh{X}}
        \ar@(ur,ul)[r]^{(\s{}\circ\Rf{n}{f}{})_{n\in\Nn}}
        \ar@(dr,dl)[r]_{(\Rf{n}{[\s{f}]}{}\circ\s{})_{n\in\Nn}}
        \ar@{}[r]|{\displaystyle\Downarrow} &
      {\sCoh{\s{Y}}.}
    }
  \end{equation}
\end{satz}

\vspace{\abstand}

\begin{proof}
  Morphism \eqref{eqcohmorph1} is given by the following diagram in the 2-category of triangulated categories
  \[
    \xymatrix{
      & &
      {\Derb{\Coh{Y}}} \ar[rdd]^{\iota} \ar@{}[d]|{\displaystyle\Downarrow} \\
      {\Derb{\Coh{X}}} \ar[rru]^{\RF{f_*}} \ar[rr]^{\iota} \ar[rd]^{\s{}} &
      \ar@{}[d]|{\displaystyle\Downarrow} &
      {\Derb{\Mod{X}}} \ar[rd]_{\RF{f_*}} \ar@{}[dd]|{\displaystyle\Downarrow} \\
      &
      {\Derb{\sCoh{\s{X}}}} \ar[ru]_{\B{}} \ar[rd]^{\RF{f_*}} & &
      {\Derb{\Mod{Y}}} \\
      & &
      {\Derb{\sCoh{\s{Y}}},} \ar[ru]^{\B{}} \\
    }
  \]
  where the three 2-morphisms are given by \eqref{eqIdSN} and \ref{lemmadercoh}
  (note that $\s{}:\Coh{X}\rightarrow\sCoh{\s{X}})$ is exact by \ref{thmvariety}\ref{varietycohex}).

  Applying \eqref{eqcohmorph1} to objects concentrated in degree zero (i.e. objects coming from $\Coh{X}$)
  and taking cohomology gives us a morphism of $\delta$-functors
  \[
    \xymatrix@C=30mm{
      {\Coh{X}}
        \ar@(ur,ul)[r]^{(\Rf{n}{f}{})_{n\in\Nn}}
        \ar@(dr,dl)[r]_{(\B{}\circ\Rf{n}{[\s{f}]}{}\circ\s{})_{n\in\Nn}}
        \ar@{}[r]|{\displaystyle\Downarrow\varphi} &
      {\Mod{Y}.}
    }
  \]
  Using \ref{corbsmod},
  we then get the morphism from \eqref{eqcohmorph2} for a coherent $\cO_X$-module $\cF$ and an $n\in\Nn$
  by
  \[
    \s{\Rf{n}{f}{\cF}}
    \xrightarrow{\bs{Y}{\Rf{n}{f}{\cF}}{\Rf{n}{[\s{f}]}{\s{\cF}}}^{-1}(\varphi)}
    \Rf{n}{[\s{f}]}{\s{\cF}}.
  \]
  That this is indeed a morphism of $\delta$-functors follows immediately from the exactness of $\s{}$,
  from \ref{bemtaunat}
  and from the fact that $\varphi$ is a morphism of $\delta$-functors.
\end{proof}

\vspace{\abstand}

\begin{thm}\label{thmcoh}
  The canonical morphism of functors \eqref{eqcohmorph2} is an isomorphism.
  In particular, $\s{\Rf{n}{f}{\cF}}$ is canonically isomorphic to
  $\Rf{n}{[\s{f}]}{\s{\cF}}$ for all coherent $\cO_X$-modules $\cF$ and all $n\in\Z$.
\end{thm}

\vspace{\abstand}

\begin{proof}
  Because the statement is local in $Y$, we can assume without loss of generality that
  $Y=\spec{B}$ is affine for a finitely presented $A$-algebra $B$. We split the proof in several cases:

  \vspace{\abstand}

  \textbf{First} consider the case where $f:X=\P^d_Y\rightarrow Y$ is the structural morphism
  of \emph{projective $d$-space} over $Y$.
  By \cite[2.1.15, 2.1.16]{ega31}, for any $m,n\in\Z$, we have canonical isomorphisms
  \[
    \Rf{n}{f}{\cO_X(m)}
    =\left\{\begin{array}{cl}
      \displaystyle\cO_Y[T_0,\ldots,T_d]_m & \text{if $n=0$,} \\[2mm]
      \displaystyle\cO_Y[T_0,\ldots,T_d]_{-d-1-m}^\vee & \text{if $n=d$,} \\[2mm]
      \displaystyle0 & \text{otherwise,} \\
    \end{array}\right.
  \]
  where $\cO_Y[T_0,\ldots,T_d]$ denotes
  the graded free symmetric algebra over $\cO_Y$ with generators $T_0,\ldots,T_d$
  (so that its part of degree $m$ is just the free $\cO_Y$-module with basis the homogenous monomials of degree
  $m$ in the $T_i$).
  By \ref{satzProj}, \ref{satzO1}\ref{O1proj} and transfer, we have
  \[
    \Rf{n}{[\s{f}]}{[\s{\cO_X(m)}]}
    =\left\{\begin{array}{cl}
      \displaystyle\cO_{\s{Y}}\str{[T_0,\ldots,T_d]}_m & \text{if $n=0$,} \\[2mm]
      \displaystyle\cO_{\s{Y}}\str{[T_0,\ldots,T_d]}_{-d-1-m}^\vee & \text{if $n=d$,} \\[2mm]
      \displaystyle0 & \text{otherwise.} \\
    \end{array}\right.
  \]
  Since a *monomial of degree $m$ is the same as a monomial of degree $m$,
  and since $\s{}$ respects duals by \ref{cordual},
  we see that
  $\s{\Rf{n}{f}{\cO_X(m)}}=\Rf{n}{[\s{f}]}{[\s{\cO_X(m)}]}$ for all $m$ and $n$.
  By additivity, the theorem is hence true for our special choice of $f$ and for all $\cF$
  of the form $\cO_X(m)^l$ for $l\in\Nn$ and $m\in\Z$.

  As a next step, we prove the theorem for all coherent sheaves on $\P^d_Y$ by decreasing induction on $n$
  (this part closely resembles Hartshorne's proof of the ``Theorem on Formal Functions" in \cite{hartshorne}):
  Since $\Rf{n}{f}{}$ and $\Rf{n}{[\s{f}]}{}$ both vanish for $n>d$, the theorem holds trivially in those cases.
  For the inductive step, assume that the theorem holds for all $n'>n\in\Nn$,
  and let $\cF$ be an arbitrary coherent sheaf on $X$.
  By \cite[2.2.2(iv)]{ega31},
  there exists an epimorphism $\cG:=\cO_X(m)^l\twoheadrightarrow\cF$ for suitable $l\in\Nn$ and $m\in\Z$,
  so that we have a short exact sequence
  \[
    0\longrightarrow
    \cH\longrightarrow
    \cG\longrightarrow
    \cF\longrightarrow
    0
  \]
  of coherent $\cO_X$-modules.
  By \ref{satzcohmorph}, we get an induced commutative diagram of *coherent $\cO_{\s{X}}$-modules
  with exact rows as follows:
  \[
    \xymatrix@C=5mm@R=15mm{
      {\s{\Rf{n}{f}{\cH}}} \ar[r] \ar[d]_{\alpha} &
      {\s{\Rf{n}{f}{\cG}}} \ar[r] \ar[d]_{\beta}^{\wr} &
      {\s{\Rf{n}{f}{\cF}}} \ar[r] \ar@{.>}[d]_{\gamma} &
      {\s{\Rf{n+1}{f}{\cH}}} \ar[r] \ar[d]_{\delta}^{\wr} &
      {\s{\Rf{n+1}{f}{\cG}}} \ar[d]_{\varepsilon}^{\wr} \\
      {\Rf{n}{[\s{f}]}{[\s{\cH}}]} \ar[r] &
      {\Rf{n}{[\s{f}]}{[\s{\cG}}]} \ar[r] &
      {\Rf{n}{[\s{f}]}{[\s{\cF}}]} \ar[r] &
      {\Rf{n+1}{[\s{f}]}{[\s{\cH}}]} \ar[r] &
      {\Rf{n+1}{[\s{f}]}{[\s{\cG}}]} \\
    }
  \]
  By the first part of the proof, $\beta$ and $\varepsilon$ are isomorphisms,
  and by our inductive hypothesis, $\delta$ is an isomorphism.
  Then by the five lemma, since $\beta$ and $\delta$ are epimorphisms and $\varepsilon$ is a monomorphism,
  $\gamma$ is an epimorphism.

  Since $\cF$ was chosen arbitrarily, this conclusion also applies to $\cH$, i.e. $\alpha$ is also an epimorphism.
  But then we can apply the five lemma again, using that $\alpha$ is an epimorphism and that
  $\beta$ and $\delta$ are monomorphisms,
  to conclude that $\gamma$ is a monomorphism and hence an isomorphism as desired.

  \vspace{\abstand}

  Having settled the theorem for projective space,
  we now consider the \textbf{second case} where $f:X\hookrightarrow Y$ is a \emph{closed immersion}, i.e.
  $X=\spec{B/\b}$ for an ideal $\b$ of $B$.
  Since $f_*$ and $[\s{f}]_*$ are exact in this case
  (note that $\s{f}$ is a *closed immersion by \ref{satzP}),
  we only have to show
  $\s{f_*\tilde{M}}\cong[\s{f}]_*[\s{\tilde{M}}]$
  for all $B/\b$-modules $M$ of finite type
  or --- equivalently --- that
  $\bigl[\s{f_*\tilde{M}}\bigr](\s{Y})\cong\bigl[[\s{f}]_*[\s{\tilde{M}}]\bigr](\s{Y})$.
  Now
  \[
    \bigl[\s{f_*\tilde{M}}\bigr](\s{Y})
    \stackrel{\ref{satzaff}}{\cong}[f_*\tilde{M}](Y)\otimes_B\s{B}=
    \tilde{M}(X)\otimes_B\s{B}=
    M\otimes_B\s{B}
  \]
  and (since $B\longrightarrow C:=B/\b$ is a \emph{finite} ring homomorphism)
  \begin{multline*}
    \bigl[[\s{f}]_*[\s{\tilde{M}}]\bigr](\s{Y})
    =[\s{\tilde{M}}](\s{X})
    \stackrel{\ref{satzaff}}{\cong}\tilde{M}(X)\otimes_C\s{C} \\
    =M\otimes_C\s{C}
    \stackrel{\ref{satzfinite}}{\cong}M\otimes_C\bigl(C\otimes_B\s{B}\bigr)
    \cong M\otimes_B\s{B},
  \end{multline*}
  so the theorem is true for closed immersions as well.

  \vspace{\abstand}

  As a \textbf{third case}, we take an arbitrary \emph{projective} morphism $f:X\rightarrow Y$.
  Since $Y$ is affine, it admits an ample bundle, which
  implies (see \cite[5.5.4(ii)]{ega2}) that there is a $d\in\Nn$ for which $f$ factorizes as
  $X\stackrel{i}{\hookrightarrow}\P_Y^d\xrightarrow{\pi}Y$,
  where $i$ is a closed immersion and
  $\pi$ is the structural morphism.
  Then for every coherent $\cO_X$-modules $\cF$ and every $n\in\Z$, we have
  (because $i_*$ and $[\s{i}]_*$ are exact)
  \begin{multline*}
    \s{\Rf{n}{f}{\cF}}
    =\s{\Rf{n}{[\pi i]}{\cF}}
    =\s{\Rf{n}{\pi}{i_*\cF}}
    \stackrel{\mathrm{1. case}}{\cong}\Rf{n}{[\s{\pi}]}{\s{i_*\cF}} \\
    \stackrel{\mathrm{2. case}}{\cong}\Rf{n}{[\s{\pi}]}{[\s{i}]_*\s{\cF}}
    =\Rf{n}{[\s{f}]}{\s{\cF}},
  \end{multline*}
  and the proof of this case is complete.

  \vspace{\abstand}

  Finally we consider the \textbf{general case} of an arbitrary \emph{proper} morphism $f:X\rightarrow Y$
  and imitate Grothendieck's proof of the finiteness theorem for coherent modules \cite[3.2.1]{ega31}.
  Consider the full subcategory $\cC$ of $\Coh{X}$ consisting of those coherent sheaves for which the theorem
  holds. We claim that $\cC$ has the following properties:
  \begin{enumerate}
    \item\label{serre1}
      $\cC$ is \emph{exact}, i.e.
      if $0\rightarrow\cF'\rightarrow\cF\rightarrow\cF''\rightarrow 0$ is a short exact sequence in $\Coh{X}$
      and if two of the three sheaves $\cF'$, $\cF$ and $\cF''$ belong to $\cC$,
      then so does the third (compare \cite[3.1.1]{ega31}).
    \item\label{serre2}
      If a coherent $\cO_X$-module $\cF$ belongs to $\cC$, then every direct factor of $\cF$ also belongs to $\cC$.
  \end{enumerate}
  Let $0\rightarrow\cF'\rightarrow\cF\rightarrow\cF''\rightarrow 0$ be a short exact sequence as in \ref{serre1}.
  Applying the morphism of $\delta$-functors \eqref{eqcohmorph2},
  we get the following commutative diagram with exact rows
  \[
    \xymatrix@C=3mm@R=15mm{
      {\ldots} \ar[r] &
      {\s{\Rf{n-1}{f}{\cF''}}} \ar[r]^{\delta} \ar@{.>}^{\gamma_{n-1}}[d] &
      {\s{\Rf{n}{f}{\cF'}}} \ar[r] \ar@{.>}[d]^{\alpha_n} &
      {\s{\Rf{n}{f}{\cF}}} \ar[r] \ar@{.>}[d]^{\beta_n} &
      {\s{\Rf{n}{f}{\cF''}}} \ar[r]^{\delta} \ar@{.>}[d]^{\gamma_n} &
      {\s{\Rf{n+1}{f}{\cF'}}} \ar[r] \ar@{.>}[d]^{\alpha_{n+1}} &
      {\ldots} \\
      {\ldots} \ar[r] &
      {\Rf{n-1}{[\s{f}]}{\s{\cF''}}} \ar[r]_{\delta} &
      {\Rf{n}{[\s{f}]}{\s{\cF'}}} \ar[r] &
      {\Rf{n}{[\s{f}]}{\s{\cF}}} \ar[r] &
      {\Rf{n}{[\s{f}]}{\s{\cF''}}} \ar[r]_{\delta} &
      {\Rf{n+1}{[\s{f}]}{\s{\cF'}}} \ar[r] &
      {\ldots} \\
    }
  \]
  If two of $\cF'$, $\cF$ and $\cF''$ belong to $\cC$, then for every $n$,
  two of $\alpha_n$, $\beta_n$ and $\gamma_n$ are isomorphisms.
  The five lemma shows that then all $\alpha_n$, $\beta_n$ and $\gamma_n$ are isomorphisms
  and hence $\cF'$, $\cF$ and $\cF''$ all belong to $\cC$, which proves \ref{serre1}.

  For \ref{serre2}, let $\cF$ be a coherent $\cO_X$-module in $\cC$,
  and let $\cF_1$ be a direct factor of $\cF$.
  Putting $\cF_2:=\cF/\cF_1$, we get a split short exact sequence
  \[
    0\longrightarrow
    \cF_1\longrightarrow
    \cF=\cF_1\oplus\cF_2\longrightarrow
    \cF_2\longrightarrow
    0
  \]
  and hence for any $n$ a morphism of split short exact sequences
  \[
    \xymatrix@C=5mm@R=15mm{
      0 \ar[r] &
      {\s{\Rf{n}{f}{\cF_1}}} \ar[r] \ar@{.>}[d]^{\alpha_n} &
      {\s{\Rf{n}{f}{\cF}}=\s{\Rf{n}{f}{\cF_1}}\oplus\s{\Rf{n}{f}{\cF_2}}} \ar[r] \ar[d]^{\alpha_n\oplus\beta_n}_{\wr} &
      {\s{\Rf{n}{f}{\cF_1}}} \ar[r] \ar@{.>}[d]^{\beta_n} &
      0 \\
      0 \ar[r] &
      {\Rf{n}{[\s{f}]}{\s{\cF_1}}} \ar[r] &
      {\Rf{n}{[\s{f}]}{\s{\cF}}=\Rf{n}{[\s{f}]}{\s{\cF_1}}\oplus\Rf{n}{[\s{f}]}{\s{\cF_2}}} \ar[r] &
      {\Rf{n}{[\s{f}]}{\s{\cF_1}}} \ar[r] &
      0 \\
    }
  \]
  with an isomorphism $\alpha_n\oplus\beta_n$ (because $\cF$ is in $\cC$).
  It follows immediately that $\alpha_n$ and $\beta_n$ must also be isomorphisms,
  i.e. $\cF_1$ and $\cF_2$ also belong to $\cC$, which proves \ref{serre2}.

  In order to finish the proof of the theorem, we have to show that every coherent $\cO_X$-module belongs to $\cC$,
  and we want to do so by using \emph{dévissage}: By \cite[3.1.3]{ega31},
  a full subcategory $\cC$ of $\Coh{X}$ satisfying \ref{serre1} and \ref{serre2}
  contains \emph{all} coherent $\cO_X$-modules if (and only if) for every irreducible closed subscheme $Z$ of $X$,
  there is a sheaf with support $Z$ in $\cC$.

  Let $Z\stackrel{i}{\hookrightarrow}X$ be a closed immersion with $Z$ irreducible.
  Assume that we have found a coherent sheaf $\cF_Z$ of $\cO_Z$-modules with support $Z$ such that the theorem holds for
  $\cF_Z$ and the (obviously proper) morphism $f\circ i:Z\rightarrow Y$.
  Then $\cF:=i_*\cF_Z$ is a coherent sheaf of $\cO_X$-modules with support $Z$,
  and
  \[
    \s{\Rf{n}{f}{\cF}}
    \cong\s{\Rf{n}{[fi]}{\cF_Z}}
    \cong\Rf{n}{[\s{(fi)}]}{\s{\cF_Z}}
    \cong\Rf{n}{[\s{f}]}{\bigl[(\s{i})_*\s{\cF_Z}\bigr]}
    \stackrel{\mathrm{2. case}}{\cong}\Rf{n}{[\s{f}]}{\s{\cF}},
  \]
  i.e. $\cF$ belongs to $\cC$.
  Thus without loss of generality,
  we only have to consider the case $Z=X$
  and therefore must exhibit a sheaf in $\cC$ with support $X$.

  By Chow's lemma \cite[5.6.2]{ega2},
  there is a projective and surjective morphism $g:X'\longrightarrow X$,
  with $X'$ irreducible,
  such that the composition $f\circ g:X'\longrightarrow Y$ is projective.
  Let $\cO_{X'}(1)$ be a very ample bundle for $g$. Then by \cite[2.2.1]{ega31} and \cite[3.4.7]{ega2},
  there is an $m\in\Nn$ such that $\cF:=g_*\cO_{X'}(m)$ has support $X$ and such that
  \begin{equation}\label{eqacyclic1}
    \forall n>0:\ \Rf{n}{g}{\cO_{X'}(m)}=0.
  \end{equation}
  From \eqref{eqacyclic1}, we learn two things.
  First, using the spectral sequence
  $\Rf{p}{f}{\Rf{q}{g}{\cO_{X'}(m)}}\Rightarrow\Rf{p+q}{[fg]}{\cO_{X'}(m)}$,
  we get
  \begin{equation}\label{eqacyclic2}
    \forall n\in\Z:\ \Rf{n}{f}{\cF}\cong\Rf{n}{[fg]}{\cO_{X'}(m)}.
  \end{equation}
  Second,
  applying the third case to $\cO_{X'}(m)$ and $g$, we get
  \[
    \forall n>0:\ \Rf{n}{[\s{g}]}{\s{\cO_{X'}(m)}}
    \cong\s{\Rf{n}{g}{\cO_{X'}(m)}}
    \stackrel{\eqref{eqacyclic1}}{=}0,
  \]
  and then, using the spectral sequence
  $\Rf{p}{[\s{f}]}{\Rf{q}{[\s{g}]}{\s{\cO_{X'}(m)}}}\Rightarrow\Rf{p+q}{[\s{(fg)}]}{\s{\cO_{X'}(m)}}$,
  \begin{equation}\label{eqacyclic3}
    \forall n\in\Z:\ \Rf{n}{[\s{f}]}{\s{\cF}}\cong\Rf{n}{[\s{(fg)}]}{\s{\cO_{X'}(m)}}.
  \end{equation}
  Combining these and applying the third case again, this time to $\cO_{X'}(m)$ and $fg$,
  we get
  \[
    \s{\Rf{n}{f}{\cF}}
    \stackrel{\eqref{eqacyclic2}}{\cong}\s{\Rf{n}{[fg]}{\cO_{X'}(m)}}
    \stackrel{\mathrm{3. case}}{\cong}\Rf{n}{[\s{(fg)}]}{\s{\cO_{X'}(m)}}
    \stackrel{\eqref{eqacyclic3}}{\cong}\Rf{n}{[\s{f}]}{\s{\cF}}
  \]
  for all $n\in\Z$, i.e. $\cF$ belongs to $\cC$, and the proof of the theorem is complete.

\end{proof}

\vspace{\abstand}

\begin{cor}\label{corcoh}
  If $k$ is a *field and
  if $X$ is a proper $k$-scheme,
  we have a canonical isomorphism
  \[
    \HZar{n}{X}{\cF}\stackrel{\sim}{\longrightarrow}\HZar{n}{\s{X}}{\s{\cF}}
  \]
  of finite dimensional $k$-vector spaces
  for every coherent $\cO_X$-module $\cF$ and every $n\in\N_0$,
\end{cor}

\vspace{\abstand}

\begin{proof}
  This follows immediately from \ref{thmcoh}, applied to $f:X\rightarrow\spec{k}$,
  and from \ref{satzTNlsNiso}:
  \begin{multline*}
    \HZar{n}{X}{\cF}
    =\HZar{n}{X}{\cF}\otimes_kk
    =\HZar{n}{X}{\cF}\otimes_k\s{k} \\
    \stackrel{\ref{satzTNlsNiso}}{=}\sls{\s{\Rf{n}{f}{\cF}}}
    \stackrel{\ref{thmcoh}}{=}\sls{\Rf{n}{[\s{f}]}{\s{\cF}}}
    =\HZar{n}{\s{X}}{\s{\cF}}.
  \end{multline*}
\end{proof}

\vspace{\abstand}

\begin{cor}\label{corsff}
  For a *field $k$ and a proper $k$-scheme $X$, the functor
  $\s{}:\Coh{X}\longrightarrow\sCoh{\s{X}}$ is exact and fully faithful.
\end{cor}

\vspace{\abstand}

\begin{proof}
  We already know that $\s{}$ is exact (and faithful) from \ref{thmvariety}\ref{varietycohex},
  even if $X$ is not proper over $k$.

  If $f:X\rightarrow\spec{k}$ is proper,
  and if $\cF$ and $\cG$ are coherent $\cO_X$-modules,
  we have
  \begin{multline*}
    \Hom{\cO_{\s{X}}}{\s{\cF}}{\s{\cG}}
    =\Bigl[\HOM{\s{X}}{\s{\cF}}{\s{\cG}}\Bigr](\s{X})
    \stackrel{\ref{corkhom}}{\cong}\Bigl[\s{\HOM{X}{\cF}{\cG}}\Bigr](\s{X}) \\
    =\HZar{0}{\s{X}}{\s{\HOM{X}{\cF}{\cG}}}
    \stackrel{\ref{corcoh}}{\cong}\HZar{0}{X}{\HOM{X}{\cF}{\cG}}
    =\Hom{\cO_X}{\cF}{\cG},
  \end{multline*}
  which proves fully faithfulness.
\end{proof}

\vspace{\abstand}

\begin{cor}\label{corpicinj}
  For a *field $k$ and a proper $k$-scheme $X$,
  the canonical group homomorphism $\s{}:\Pic{X}\longrightarrow\sPic{\s{X}}$ from \ref{corpic}
  is \emph{injective}.
\end{cor}

\vspace{\abstand}

\begin{proof}
  This follows immediately from the fact that $\s{}:\Modfp{X}\longrightarrow\sModfp{\s{X}}$
  is fully faithful by \ref{corsff}.
\end{proof}

\vspace{\abstand}

\begin{bsp}
  Let $k$ be a *field, and consider projective $d$-space over $k$ for a $d\in\Np$.
  Then the monomorphism $\Pic{\P^d_k}\hookrightarrow\Pic{\str{\P^d_k}}$ from \ref{corpicinj}
  is explicitly given by the following
  commutative diagram of abelian groups:
  \[
    \xymatrix{
      {m} \ar@{}[r]|{\in} \ar@{|->}[d] &
      {\Z} \ar@{^{(}->}[rr]^{*} \ar[d]_{\wr} & &
      {\Zs} \ar@{}[r]|{\ni} \ar[d]^{\wr} &
      {m} \ar@{|->}[d] \\
      {\cO_{\P^d_k}(m)} \ar@{}[r]|{\in} &
      {\Pic{\P^d_k}} \ar@{^{(}->}[rr]_{\s{}} & &
      {\sPic{\str{\P^d_k}}} \ar@{}[r]|{\ni} &
      {\cO_{\str{\P^d_k}}(m)} \\
    }
  \]
\end{bsp}

\vspace{\abstand}

\begin{cor}\label{corEP}
  Let $X$ be proper over a *field $k$,
  and let $\cF$ be a coherent $\cO_X$-module.
  Then $\chi(\cF)$, the \emph{Euler-Poincaré characteristic of $\cF$},
  equals $\chi(\s{\cF})$, the *Euler-Poincaré characteristic of $\s{\cF}$.
\end{cor}

\vspace{\abstand}

\begin{proof}
  We have
  \begin{multline*}
    \chi(\s{\cF})
    =\sideset{^*}{}\sum_{n=0}^{\str{\dim(\s{X})}}(-1)^n\cdot\str{\dim\Bigl[\HZar{n}{\s{X}}{\s{\cF}}\Bigr]} \\
    \stackrel{\ref{thmvariety}\ref{varietydim}}{=}
      \sum_{n=0}^{\dim X}(-1)^n\cdot\str{\dim\Bigl[\HZar{n}{\s{X}}{\s{\cF}}\Bigr]}
    \stackrel{\ref{corcoh}}{=}\sum_{n=0}^{\dim X}(-1)^n\cdot\dim\Bigl[\HZar{n}{X}{\cF}\Bigr]
    =\chi(\cF).
  \end{multline*}
\end{proof}

\vspace{\abstand}

\begin{cor}\label{corIdealProd}
  Let $X$ be a $k$-scheme of finite type,
  and let $\cI$ and $\cJ$ be two sheaves of ideals in $\cO_X$.
  Then $\s{[\cI\cdot\cJ]}=[\s{\cI}]\cdot[\s{\cJ}]$ as *ideals of $\cO_{\s{X}}$.
\end{cor}

\vspace{\abstand}

\begin{proof}
  Let $Z$ be the closed subscheme of $X$ given by $\cI\cdot\cJ$,
  and let $i:Z\hookrightarrow X$ be the corresponding closed immersion.
  Then we have an exact sequence of coherent $\cO_X$-modules
  \[
    \cI\otimes_{\cO_X}\cJ\longrightarrow
    \cO_X\longrightarrow
    i_*\cO_Z\longrightarrow
    0
  \]
  and hence by \ref{satzox}, \ref{cortensor}, \ref{thmvariety}\ref{varietycohex} and \ref{thmcoh}
  an exact sequence
  \[
    [\s{\cI}]\otimes_{\cO_{\s{X}}}[\s{\cJ}]\stackrel{\varphi}\longrightarrow
    \cO_{\s{X}}\longrightarrow
    [\s{i}]_*\cO_{\s{Z}}\longrightarrow
    0
  \]
  of *coherent $\cO_{\s{X}}$-*modules.
  By transfer, the image of $\varphi$ is $[\s{\cI}]\cdot[\s{\cJ}]$ and the ideal
  defining $\s{Z}$,
  which in turn is $\s{[\cI\cdot\cJ]}$ by \ref{corkNNcompat}.
\end{proof}

\vspace{\abstand}


\section{The shadow map}

\vspace{\abstand}

Let $\langle K,|.|:K\rightarrow\Rge\rangle$ be a non-trivially valued field
with locally compact completion $\langle\Kh,|.|\rangle$.
Examples of such fields are $\Q$, $\R$ and $\C$ with their usual absolute value,
$\Q$ or $\Q_p$, equipped with the $p$-adic value $|.|_p$ for a prime $p$
or --- more generally --- local fields.

Assume that $\langle\Kh,|.|\rangle$ is an element of our superstructure $\hat{M}$
(which is no restriction, since we can always choose an appropriately large $M$).

Then $\langle\Ks,|.|\rangle$
and $\langle\Ksh,|.|\rangle$ are elements of $\widehat{\str{M}}$,
where $\Ks\subseteq\Ksh$ are fields,
and $|.|:\Ks\rightarrow\str{\Rge}$ and $|.|:\Ksh\rightarrow\str{\Rge}$ are maps such that
\[
  \xymatrix{
    {\Ks} \ar@{^{(}->}[rr] \ar[rd]^{|.|} & &
    {\Ksh} \ar[ld]_{|.|} \\
    & {\str{\Rge}} \\
    & {\Rge} \ar@{^{(}->}[u]\\
    K \ar@{^{(}->}[rr] \ar@{^{(}->}[uuu]^{*} \ar[ru]_{|.|} & &
    {\Kh} \ar@{^{(}->}[uuu]_{*} \ar[lu]^{|.|} \\
  }
\]
commutes. By transfer we have

\begin{enumerate}
  \item[\Mi]
    $\forall x\in\Ksh:\ |x|=0\Longleftrightarrow x=0$,\\[0mm]
  \item[\textbf{(M2)}]
    $\forall x,y\in\Ksh:\ |x\cdot y|=|x|\cdot|y|$ and\\[0mm]
  \item[\textbf{(M3)}]
    $\forall x,y\in\Ksh:\ |x+y|\leq|x|+|y|$.\\[0mm]
\end{enumerate}
Define the set of \emph{finite} elements of $\Ks$ by
\[
  \Ksf:=\left\{x\in\Ks\ \left\vert\ \exists C\in\Rge:\ |x|<C\right.\right\}
\]
and the set of \emph{infinitesimal} elements of $\Ks$ by
\[
  \Ksi:=\left\{x\in\Ks\ \left\vert\ \forall\eps\in\Rp:\ |x|<\eps\right.\right\}.
\]

\vspace{\abstand}

\begin{satz}\label{satzvalring}
  $\Ksf\varsubsetneq\Ks$ is a valuation ring with
  maximal ideal $\Ksi$ and
  residue field canonically isomorphic to $\Kh$.
  We call the projection $\Ksf\twoheadrightarrow\Kh$ the \emph{shadow map}, denote it by $\sh$,
  and consequently get a commutative diagram of ring homomorphisms with exact row
  \begin{equation}\label{eqvalring}
    \xymatrix{
      & & & {K} \ar[ld]_{*} \ar@{^{(}->}[d] \\
      0 \ar[r] & {\Ksi} \ar[r] & {\Ksf} \ar@{^{(}->}[d] \ar[r]^{\sh} & {\Kh} \ar[r] & 0 \\
      & & {\Ks.} \\
    }
  \end{equation}
\end{satz}

\vspace{\abstand}

\begin{proof}
  \Mii\ and \Miii\ immediately imply that $\Ksf$ is a subring of $\Ks$.
  Since the value on $K$ is non-trivial, the set of values is not bounded,
  so by transfer $\Ks$ contains elements of infinite value, and $\Ksf$ is a proper subring of $\Ks$.

  If $x\in\Ks$ is not finite, it in particular satisfies $|x|>1$.
  Then $|\frac{1}{x}|<1$ (by \Mii), i.e. $\frac{1}{x}$ is finite. This proves that $\Ksf$ is indeed a valuation ring.

  For a finite $x\in\Ksf\setminus\{0\}$,
  $\frac{1}{x}$ is obviously infinite if and only if $x$ is infinitesimal, which shows that $\Ksi$ is the maximal
  ideal of $\Ksf$.

  Choose an infinite natural number $h$. We define a ring homomorphism $\alpha:\Kh\longrightarrow\Ksf/\Ksi$ by
  sending the class of a Cauchy sequence $(x_n)$ in $K$ to $x_h$.
  This is well-defined, because Cauchy sequences are bounded (so that $x_h\in\Ksf$)
  and because $\lim_{n\rightarrow\infty}x_n=0$ implies $x_h\in\Ksi$.
  Furthermore, $\alpha$ does not depend on $h$: If $h'$ is another infinite natural number,
  and if $(x_n)$ is a Cauchy sequence in $K$, then $x_h-x_{h'}$ is infinitesimal.
  Since $\Kh$ is a field, $\alpha$ is automatically injective.

  To prove that it is also surjective, we need the fact that $\Kh$ is locally compact:
  This fact implies that there exists an $\eps\in\Rp$ and a compact subset $A$ of $\Kh$
  such that
  \[
    \ball{\eps}{0}{\Kh}:=\bigl\{x\in\Kh\ \bigl\vert\ |x|<\eps\bigr\}\subseteq A.
  \]
  Now let $x$ be an arbitrary element of $\Ksf$,
  let $C\in\Rge$ with $|x|<C$,
  let $\pi\in K$ with $|\pi|>1$,
  and let $\n\in\Np$ with $|\pi^n|=|\pi|^n>=\frac{C}{\eps}$.
  Because multiplication by $\pi^n$ is a homeomorphism from $\Kh$ to itself,
  $B:=\pi^nA$ is also compact, and we have
  \[
    \ball{C}{0}{\Kh}\subseteq\ball{\eps|\pi^n|}{0}{\Kh}\subseteq B
  \]
  and hence
  \[
    x
    \in\bigl\{y\in\Ksh\ \bigl\vert\ |y|<C\bigr\}
    =\str{\ball{C}{0}{\Kh}}
    \subseteq\str{B}
    \subseteq\Ksh.
  \]
  According to the nonstandard characterization of compactness, applied to $B$,
  any element of $\str{B}$ is infinitesimally close to an element of $B$,
  so there is an $\hat{x}$ in $\Kh$ with $x-\hat{x}\in\Ksi$, i.e. $x=\alpha(\hat{x})$.
\end{proof}

\vspace{\abstand}

\begin{cor}\label{corvalproper}
  Let $X$ be a proper scheme over $\Ks$. Then the canonical map
  $X(\Ksf)\longrightarrow X(\Ks)$ is bijective.
\end{cor}

\vspace{\abstand}

\begin{proof}
  This follows immediately from \ref{satzvalring} and the valuative criterion of properness
  \cite[II.4.7]{hartshorne}.
\end{proof}

\vspace{\abstand}

\begin{cor}\label{corshadow}
  Let $X$ be a proper scheme over $K$.
  Then there is a canonical shadow map
  $\sh_X:[\str{X}](\Ks)\longrightarrow X(\Kh)$, induced by $\sh:\Ksf\longrightarrow\Kh$,
  such that the following diagram commutes:
  \[
    \xymatrix{
      & {X(K)} \ar[ldd] \ar[d] \\
      & {X(\Kh)} \\
      {X(\Ks)} \ar[r]_-{\sim}^-{\ref{thmBS}} & {[\str{X}](\Ks).} \ar@{.>}[u]_{\sh_X}  \\
    }
  \]
\end{cor}

\vspace{\abstand}

\begin{proof}
  Applying the functor $X(\_)$ to \eqref{eqvalring},
  we get the following commutative diagram,
  in which $\alpha$ is bijective by \ref{corvalproper}, so that we can define $\sh_X$ as
  $(\sh\circ\alpha^{-1}\circ\beta^{-1})$:
  \[
    \xymatrix{
      & {X(K)} \ar[ld] \ar[d] \\
      {X(\Ksf)} \ar[d]_{\wr}^{\alpha} \ar[r]^{\sh} & {X(\Kh)} \\
      {X(\Ks)} \ar[r]_-{\sim}^-{\beta} & {[\str{X}](\Ks).} \ar@{.>}[u]_{\sh_X} \\
    }
  \]
\end{proof}

\vspace{\abstand}

\begin{bsp}\label{bspshadow}
  Let $X\subseteq\P_{K}^d$ be a projective variety over $K$,
  and let $x=(x_0:\ldots:x_d)$ be a $\Ks$-valued point of $\str{X}$.
  Put $C:=\max\{|x_0|,\ldots,|x_d|\}\in\str{\Rp}$.
  Then
  \[
    \sh_X(x)=\left(\sh\left[\frac{x_0}{C}\right]:\ldots:\sh\left[\frac{x_d}{C}\right]\right)\in X(\Kh)
    \subseteq\P_{K}^d(\Kh).
  \]
\end{bsp}

\vspace{\abstand}


\section{Resolution of singularities and weak factorization}

\vspace{\abstand}

For us, a \emph{variety} over a field $k$ is an integral, separated $k$-scheme of finite type.
Similarly, if $k$ is internal, a \emph{*variety} over $k$ is a *integral, *separated *scheme in $\sSchfp{k}$.

\vspace{\abstand}

\begin{lemma}\label{lemmaVar}
  Let $k$ be a *field in $\scR$, and let $X$ be a $k$-variety. Then $\s{X}$ is a $k$-*variety.
\end{lemma}

\vspace{\abstand}

\begin{proof}
  This follows immediately from \ref{satzP} and \ref{thmvariety}\ref{varietyinteger}.
\end{proof}

\vspace{\abstand}

Let $k$ be a field, and let $X$ be a projective $k$-variety.
Then for us,
a \emph{resolution (of singularities) of $X$} is a proper, birational $k$-morphism
$X'\rightarrow X$, where $X'$ is a projective, smooth $k$-variety.

\vspace{\abstand}

\begin{satz}\label{satzResolution}
  Let $k$ be a *field in $\scR$ of \emph{external} characteristic zero,
  let $n\in\N_+$,
  and let $X$ be a *projective $k$-*variety which admits a *closed embedding into $\str{\P^n_k}$ of finite *degree.
  Then there exists a *resolution $f:X'\rightarrow X$ of $X$.
\end{satz}

\vspace{\abstand}

\begin{proof}
  By \ref{corDegDimP}, there is a projective $k$-variety $Y$ with $\s{Y}=X$,
  and by Hironaka's celebrated result on resolutions of singularities in characteristic zero,
  there exists a resolution $g:Y'\rightarrow Y$ of $Y$.

  Then $X':=\s{Y'}$ is a *projective, *smooth $k$-*variety by \ref{satzP}, \ref{satzSmooth} and \ref{lemmaVar},
  and $f:=\s{g}:X'\rightarrow X$ is *proper and *birational by \ref{satzP} and \ref{satzVarBirat}.
\end{proof}

\vspace{\abstand}

Using \ref{satzResolution},
we can now easily give a conceptual proof of the following classical result of Eklof (see \cite{eklof}):

\begin{cor}\label{corResolution}
  For any pair $(n,d)$ of natural numbers, there exists a bound $C\in\N_+$,
  such that for any field $k$ of characteristic $p\geq C$
  and any closed subvariety $X$ of $\P^n_k$ of degree $d$,
  there exists a resolution of singularities of $X$.
\end{cor}

\vspace{\abstand}

\begin{proof}
  Assume the statement is false. Then for every $i\in\N_+$,
  we find a field $k_i$ of characteristic $p_i\geq i$ and a closed subvariety $X_i$ of $\P^n_{k_i}$
  of degree $d$ which does not admit a resolution.

  We then take the full subcategory of $\Rings$ with objects $(k_i)_{i\in\N_+}$ as our base category $\cB$,
  choose an infinite $j\in\Ns$
  and get a *field $k_j$ of *characteristic $p_j\geq j$ in $\scR$
  and a *closed *subvariety $X_j$ of $\str{\P^n_k}$ of *degree $d$ which does not admit a *resolution.

  But since $p_j$ is infinite, the external characteristic of $k_j$ is zero,
  and \ref{satzResolution} states that there can be no such $X_j$.
  Thus our assumption leads to a contradiction, and the corollary is proven.
\end{proof}

\vspace{\abstand}

\begin{defi}\label{defBounded}
  Let $k$ be a field,
  let $U$ be an open subscheme of a projective $k$-variety $X$,
  and let $n\in\Nn$ be a natural number.
  We say that $U$ has \emph{complexity $n$} if $X\setminus U$, equipped with its reduced structure,
  has at most $n$ irreducible components
  and if all those components have degree at most  $n$.
\end{defi}

\vspace{\abstand}

\begin{lemma}\label{lemmaComplexity}
  Let $k$ be *field in $\scR$,
  let $X$ be a projective $k$-variety,
  and let $U'$ be a *open subscheme of $\s{X}$ of \emph{finite} *complexity.
  Then there is an open subscheme $U$ of $X$ with $\s{U}=U'$.
\end{lemma}

\vspace{\abstand}

\begin{proof}
  By definition of complexity, there is an $n\in\Nn$, such that
  $[\s{X}]\setminus U'=Z'_1\cup\ldots\cup Z'_n$ with *integral *closed *subschemes $Z'_i$ of $\s{X}$
  of *degree at most $n$,
  and by \ref{corDegDim}, there exist integral closed subschemes $Z_1,\ldots,Z_n$ of $X$ with
  $\s{Z_i}=Z'_i$ for all $i$. Put $U:=X\setminus{\bigcup_{i=1}^nZ_i}=\bigcap_{i=1}^n[X\setminus Z_i]$. Then
  \[
    \s{U}
    =\s{\left(\bigcap_{i=1}^n[X\setminus Z_i]\right)}
    \stackrel{\ref{satzfin}\ref{finleftex}}{=}\bigcap_{i=1}^n\s{[X\setminus Z_i]}
    \stackrel{\ref{lemmaComplement}}{=}\bigcap_{i=1}^n\Bigl([\s{X}]\setminus[\s{Z_i}]\Bigr)
    =\bigcap_{i=1}^n\Bigl([\s{X}]\setminus Z'_i\Bigr)
    =U'.
  \]
\end{proof}

\vspace{\abstand}

\begin{defi}\label{defWKfact}
  Let $\Phi:X\dashrightarrow Y$ be a birational map between proper nonsingular varieties over a field $k$,
  and let $U\subseteq X$ be an open subscheme where $\Phi$ is an isomorphism.
  Then a \emph{weak factorization of $\Phi$ with respect to $U$}
  is a factoring of $\Phi$ into a sequence of blow-ups and blow-downs with nonsingular irreducible
  centers disjoint from $U$.
  The \emph{length} of a weak factorization is the number of blow-ups and blow-downs in the sequence.
\end{defi}

\vspace{\abstand}

\begin{lemma}\label{lemmaWKfact}
  Let $k$ be a *field in $\scR$,
  let $\Phi:X\rightarrow Y$ be a birational morphism between proper, smooth $k$-varieties,
  and let $U\subseteq X$ be an open subscheme where $\Phi$ is an isomorphism.
  If $\Phi$ admits a weak factorization with respect to $U$ of length $n$,
  then $\s{\Phi}:\s{X}\rightarrow\s{Y}$ admits a *weak *factorization with respect to $\s{U}$ of *length $n$.
\end{lemma}

\vspace{\abstand}

\begin{proof}
  The statement makes sense, because
  $\s{X}$ and $\s{Y}$ are *proper, *nonsingular $k$-*varieties by \ref{satzP}, \ref{satzSmooth} and \ref{lemmaVar},
  $\s{\Phi}$ is *birational by \ref{satzVarBirat},
  and $[\s{\Phi}]\vert_{\s{U}}$ is trivially an isomorphism.

  Furthermore, it follows immediately from
  \ref{satzSmooth}, \ref{thmvariety}\ref{varietyinteger} and \ref{satzBlowup}
  that $\s{}$ maps any weak factorization of $\Phi$ with respect to $U$ of length $n$
  to a *weak *factorization with respect to $\s{U}$ of *length $n$.
\end{proof}

\vspace{\abstand}

\begin{satz}\label{satzWkFact}
  Let $k$ be a *algebraically closed *field in $\scR$ of \emph{external} characteristic zero,
  let $n\in\N_+$,
  let $X$ and $Y$ be *projective, *nonsingular $k$-schemes which admit a *closed embedding into $\str{\P^n_k}$
  of finite *degree,
  let $\Phi:X\rightarrow Y$ be a *birational \emph{morphism} of $k$-*schemes whose *graph has finite *degree,
  and let $U$ be a *open *subscheme of $X$ of finite *complexity where $\Phi$ is an isomorphism.
  Then $\Phi$ admits a *weak *factorization with respect to $U$ of \emph{finite} *length.
\end{satz}

\vspace{\abstand}

\begin{proof}
  By \ref{corNonSing}, \ref{corDegDimP}, \ref{corDegDimMorph}, \ref{satzVarBirat} and \ref{lemmaComplexity},
  there are projective, nonsingular $k$-varieties $X'$ and $Y'$,
  a birational morphism $\Phi':X'\rightarrow Y'$
  and an open subscheme $U'$ of $X'$,
  such that $\s{X'}=X$, $\s{Y'}=Y$, $\s{\Phi'}=\Phi$ and $\s{U'}=U$.
  Since $k$ is an algebraically closed field of characteristic zero
  and since $\Phi'\vert_{U'}$ is an isomorphism by \ref{corReflectsIsos},
  we know from \cite[0.1.1]{wkfact} that $\Phi'$ admits a weak factorization with respect to $U'$.
  The claim now follows immediately from \ref{lemmaWKfact}.
\end{proof}

\vspace{\abstand}

\begin{defi}\label{defDatum}
  Let $k$ be a field.
  A \emph{WF-datum over $k$} is a pair $\langle\Phi,U\rangle$,
  where $\Phi:X\rightarrow Y$ is a birational morphism between projective, nonsingular $k$-varieties
  and where $U$ is an open subscheme of $X$ where $\Phi$ is an isomorphism.
  A \emph{weak factorization of $\langle\Phi,U\rangle$ (of length $n$)} is a weak factorization of $\Phi$
  with respect to $U$ of length $n$.

  Let $N\in\Nn$ be a natural number.
  We say that the WF-datum $\langle\Phi,U\rangle$ has \emph{complexity $n$} if $X$ and $Y$ are (isomorphic to)
  closed subschemes of $\P^n_k$ of degree at most $n$,
  if the graph of $\Phi$ has degree at most $n$
  and if $U$ has complexity $n$.
\end{defi}

\vspace{\abstand}

\begin{cor}\label{corWKfact1}
  For any $N\in\Nn$, there exists a bound $C\in\Np$,
  such that for any algebraically closed field $k$ of characteristic $p\geq C$,
  any WF-datum of complexity $N$ has a weak factorization.
\end{cor}

\vspace{\abstand}

\begin{proof}
  This follows from \ref{satzWkFact} in the same way as \ref{corResolution} follows from \ref{satzResolution}.
\end{proof}

\vspace{\abstand}

\begin{cor}\label{corWKfact2}
  For any $N\in\Nn$, there exists a bound $D\in\Np$,
  such that for any algebraically closed field $k$ of characteristic zero,
  any WF-datum of complexity $N$ has a weak factorization of length at most $D$.
\end{cor}

\vspace{\abstand}

\begin{proof}
  This, again, follows in the same way as \ref{corResolution} and \ref{corWKfact1},
  using the fact that the *weak *factorization whose existence is proven in \ref{satzWkFact} has \emph{finite} *length.
\end{proof}

\vspace{\abstand}


\def\cprime{$'$}

\end{document}